\documentclass[10pt,a4paper,twoside]{article}
\usepackage{amsmath}
\usepackage{amssymb}
\usepackage{amsthm}
\usepackage{geometry}
\usepackage{graphicx} 
\usepackage{subfig}
\usepackage{url}

\date{}
\title{\textbf{Stability in the Energy Space of the Sum of $N$ Peakons for the Degasperis-Procesi Equation}}
\author{\textbf{Andr\'e Kabakouala}\\
L.M.P.T., U.F.R Sciences et Techniques, Universit\'e de Tours, Parc Grandmont,\\
37200 Tours, France.\\
\\
Andre.Kabakouala@lmpt.univ-tours.fr}

\begin{document}

\theoremstyle{plain}
\newtheorem{Theo}{Theorem}[section]
\newtheorem{Pro}{Proposition}[section]
\newtheorem{Lem}{Lemma}[section]
\newtheorem{Cor}{Corollary}[section]
\newtheorem{Def}{Definition}[section]
\newtheorem{Pre}{Proof}
\newtheorem{Rem}{Remark}[section]
\newtheorem{Ex}{Example}[section]
\newtheorem*{Ack}{Acknowledgements}
\newtheorem{Cost}{Construction}[section]
\numberwithin{equation}{section}

\maketitle

\begin{abstract}
The Degasperis-Procesi equation possesses well-known peaked solitary waves that are called peakons. Their stability has been established by Lin and Liu in \cite{MR2460268}. In this paper, we localize the proof (in some suitable sense detailed in Section \ref{Section 3}) of the stability of a single peakon. Thanks to this, we extend the result of stability to the sum of $N$ peakons traveling to the right with respective speeds $c_{1},\ldots,c_{N}$, such that the difference between consecutive locations of peakons is large enough.
\end{abstract}
\vspace{1cm}

\section{Introduction}\label{Section 1}
The Degasperis-Procesi (DP) equation 
\begin{equation}
u_{t}-u_{txx}+4uu_{x}=3u_{x}u_{xx}+uu_{xxx},~~(t,x)\in\mathbb{R}^{*}_{+}\times\mathbb{R}
\label{1.1}
\end{equation}
is completely integrable (see \cite{MR2001531}) and possesses, among others, the following invariants
\begin{equation}
E(u)=\int_{\mathbb{R}}yv~~\text{and}~~F(u)=\int_{\mathbb{R}}u^{3},
\label{1.2}
\end{equation}
where $y=(1-\partial^{2}_{x})u$ and $v=(4-\partial^{2}_{x})^{-1}u$. Substituting $u$ by $4v-v_{xx}$ in 
\eqref{1.2} and using integration by parts (we suppose that $u(\pm\infty)=v(\pm\infty)=v_{x}(\pm\infty)=0$), the conservation laws can be rewritten as
\begin{equation}
E(u)=\int_{\mathbb{R}}\left(4v^{2}+5v^{2}_{x}+v^{2}_{xx}\right)~~\text{and}~~F(u)=\int_{\mathbb{R}}
\left(-v^{3}_{xx}+12vv^{2}_{xx}-48v^{2}v_{xx}+64v^{3}\right)
\label{1.3}.
\end{equation}
One can see that the conservation law $E(\cdot)$ is equivalent to $\|\cdot\|^{2}_{L^{2}(\mathbb{R})}$. Indeed, using integration by parts
\begin{equation}
\|u\|^{2}_{L^{2}(\mathbb{R})}=\int_{\mathbb{R}}u^{2}=\int_{\mathbb{R}}(4v-v_{xx})^{2}=\int_{\mathbb{R}}\left(16v^{2}+8v^{2}_{x}+
v^{2}_{xx}\right)\le 4E(u),
\label{1.4}
\end{equation}
and applying Plancherel-Parseval identity 
\begin{equation}
E(u)=\int_{\mathbb{R}}yv=\int_{\mathbb{R}}\frac{1+\omega^{2}}{4+\omega^{2}}|\widehat{u}(\omega)|^{2}\le\int_{\mathbb{R}}|\widehat{u}(\omega)|^{2}=\|u\|^{2}_{L^{2}(\mathbb{R})},
\label{1.5}
\end{equation}
where $\widehat{u}$ denotes the Fourier transform of $u$. In the sequel we will denote 
\begin{equation}\label{n1}
\|u\|_{\mathcal{H}}=\sqrt{E(u)}.
\end{equation}
Note that, by reversing the operator $(1-\partial^{2}_{x})(\cdot)$ in \eqref{1.1}, the DP equation can be rewritten in conservation form as 
\begin{equation}
u_{t}+\frac{1}{2}\partial_{x}u^{2}+\frac{3}{2}(1-\partial^{2}_{x})^{-1}\partial_{x}u^{2}=0,
~~(t,x)\in\mathbb{R}^{*}_{+}\times\mathbb{R}
\label{1.6}.
\end{equation}

The DP equation possesses solitary waves called \textit{peakons} (see Fig. \ref{fig1}) and defined by 
\begin{equation}
u(t,x)=\varphi_{c}(x-ct)=c\varphi(x-ct)=ce^{-|x-ct|},~~c\in\mathbb{R}^{*},
~~(t,x)\in\mathbb{R}^{*}_{+}\times\mathbb{R},
\label{1.7}
\end{equation}
but they are not smooth since $\varphi_{c}\notin C^{1}(\mathbb{R})$ (see Fig. \ref{fig2}). The peakons are only global weak solutions of \eqref{1.6}. It means, for any smooth test function $\phi\in C^{\infty}(\mathbb{R}_{+}\times\mathbb{R})$, it holds
\begin{align*}
\int_{0}^{+\infty}\int_{\mathbb{R}}&\varphi_{c}(x-ct)\phi_{t}(t,x)dtdx+\frac{1}{2}\int_{0}^{+\infty}\int_{\mathbb{R}}\varphi^{2}_{c}(x-ct)\phi_{x}(t,x)dtdx\\
&+\frac{3}{2}\int_{0}^{+\infty}\int_{\mathbb{R}}(1-\partial^{2}_{x})^{-1}\varphi^{2}_{c}(x-ct)\phi_{x}(t,x)dtdx+\int_{\mathbb{R}}\varphi_{c}(x)\phi(0,x)dx=0.
\end{align*}
The goal of our work is to prove that ordered trains of peakons are stable under small perturbations in the energy space $ {\mathcal H }$ (equivalent to $L^{2}$).

\begin{Def}[Stability]
Let $c>0$ be given.  The peakon $\varphi_{c}$ is said stable in $\mathcal{H}$, if for all $\varepsilon >0$, there exists $\delta>0$ such that if 
\begin{equation}
\|u_{0}-\varphi_{c}\|_{\mathcal{H}}\le\delta,
\label{1.8}
\end{equation}
then for all $t\ge 0$, there exists $\xi(t)$ such that 
\begin{equation}
\|u(t,\cdot)-\varphi_{c}\left(\cdot-\xi(t)\right)\|_{\mathcal{H}}\le\varepsilon,
\label{1.9}
\end{equation}
where $u(t)$ is the solution to  \eqref{1.1} emanating from $u_{0}$.
\end{Def}

Lin and Liu proved in \cite{MR2460268} the stability of a single peakon under the additional condition that $(1-\partial^{2}_{x})u_{0} \in {\mathcal M}^{+}(\mathbb{R})$. Using this result and the general strategy introduced by Martel, Merle  and Tsai in \cite{MR1946336} for the generalized Korteweg-de Vries (gKdV) equation and adapted by El Dika and Molinet in \cite{MR2542735} and \cite{MR2525156} for the Camassa-Holm (CH) equation, we prove here the stability of the sum of $N$ peakons for the DP equation.

Before stating the main result we introduce the function space where will live our class of solutions to the equation. For $I$ a finite or infinite time interval of $\mathbb{R}_{+}$, we denote by $\mathcal{X}(I)$ the function space \footnote{$ W^{1,1}(\mathbb R) $ is the space of $ L^1(\mathbb R) $ functions with derivatives in $ L^1(\mathbb R) $ and $ BV(\mathbb R) $ is the space of function with bounded variation.}
\begin{equation}
\mathcal{X}(I)=\left\lbrace u\in C\left(I;H^{1}(\mathbb{R})\right)\cap L^{\infty}\left(I;W^{1,1}(\mathbb{R})\right),~ u_{x}\in L^{\infty}\left(I;BV(\mathbb{R})\right)\right\rbrace.
\label{1.10}
\end{equation}

The main result of the present paper is the following theorem.

\begin{Theo}[Stability of the Sum of $N$ Peakons]\label{Theoreme 1.1}
Let be given $N$ velocities $c_{1},\ldots,c_{N}$ such that $0<c_{1}<\ldots<c_{N}$. Let $u\in\mathcal{X}([0,T[)$, with $0<T\le +\infty$, be a solution of the DP equation. There exist $C>0$, $L_{0}>0$ and $\varepsilon_{0}>0$ only depending on the speeds $(c_{i})_{i=1}^{N}$,  such that if 
\begin{equation}
y_{0}=(1-\partial^{2}_{x})u_{0}\in\mathcal{M}^{+}(\mathbb{R})
\label{1.11}
\end{equation}
and
\begin{equation}
\left\|u_{0}-\sum_{i=1}^{N}\varphi_{c_{i}}(\cdot-z^{0}_{i})\right\|_{\mathcal{H}}\le\varepsilon^{2},~~\text{with}~~0<\varepsilon<\varepsilon_{0},
\label{1.12}
\end{equation}
for some $z^{0}_{1},\ldots,z^{0}_{N}$ satisfying 
\begin{equation}
z^{0}_{1}<\ldots<z^{0}_{N}~~\text{and}~~
z^{0}_{i}-z^{0}_{i-1}\ge L,~~\text{with}~~L>L_{0}>0,~~i=2,\ldots,N,
\label{1.13}
\end{equation}
then there exist $\xi^{1}_{1}(t),\ldots,\xi^{N}_{1}(t)$ such that 
\begin{equation}
\left\|u(t)-\sum_{i=1}^{N}\varphi_{c_{i}}(\cdot-\xi^{i}_{1}(t))\right\|_{\mathcal{H}}\le C(\sqrt{\varepsilon}+L^{-1/8}),~~\forall t\in[0,T[
\label{1.14}
\end{equation}
and 
\begin{equation}
\xi^{i}_{1}(t)-\xi^{i-1}_{1}(t)>\frac{L}{2},~~\forall t\in[0,T[,~~i=2,\ldots,N,
\label{1.15}
\end{equation}
where $\xi^{1}_{1}(t),\ldots,\xi^{N}_{1}(t)$ are defined in Subsection \ref{Subsection 4.1}.
\end{Theo}

\section{Preliminaries}\label{Section 2}
In this section, we briefly recall the global well-posedness results for the DP equation and its consequences 
(see \cite{MR2271927} and \cite{MR2249792} for details).

\begin{Theo}[Global Weak Solution; See \cite{MR2271927} and \cite{MR2249792}]\label{Theoreme 2.1}
Assume that  $u_{0}\in L^{2}(\mathbb{R})$ with  $y_{0}=(1-\partial^{2}_{x})u_{0}\in\mathcal{M}^{+}(\mathbb{R})$. Then the DP equation has a unique global weak solution $u\in\mathcal{X}(\mathbb{R}_{+})$ such that 
\begin{equation}
y(t,\cdot)=(1-\partial^{2}_{x})u(t,\cdot)\in \mathcal{M}^{+}(\mathbb{R}),~~\forall t\in\mathbb{R}_+.
\label{2.1}
\end{equation}
Moreover $ E(\cdot) $ and $ F(\cdot) $ are conserved by the flow.
\end{Theo}

\begin{Rem}[Control of $L^{\infty}$ Norm by $L^{2}$ Norm]\label{Remark 2.1}
\normalfont
From \eqref{2.1}, it holds
$$u(x)=\frac{e^{-x}}{2}\int_{-\infty}^{x}e^{x'}y(x')dx'+\frac{e^{x}}{2}\int_{x}^{+\infty}e^{-x'}y(x')dx'$$
and 
$$u_{x}(x)=-\frac{e^{-x}}{2}\int_{-\infty}^{x}e^{x'}y(x')dx'+\frac{e^{x}}{2}\int_{x}^{+\infty}e^{-x'}y(x')dx',$$
which lead to 
\begin{equation}
|u_{x}(x)|\le u(x),~~\forall x\in\mathbb{R}.
\label{2.b}
\end{equation}
Then, using the Sobolev embedding of $H^{1}(\mathbb{R})$ into $L^{\infty}(\mathbb{R})$ and \eqref{2.b}, we infer that there exists a constant $C_{S}>0$ such that 
\begin{equation}
\|u\|_{L^{\infty}(\mathbb{R})}\le C_{S}\|u\|_{H^{1}(\mathbb{R})}\le 2C_{S}\|u\|_{L^{2}(\mathbb{R})}.
\label{2.3}
\end{equation}
\end{Rem}

\begin{Lem}[Positivity; See \cite{MR2249792}]\label{Lemma 2.1}
Let $u\in H^{1}(\mathbb{R})$ with $y=(1-\partial^{2}_{x})u\in\mathcal{M}^{+}(\mathbb{R})$. If $k_{1}\ge 1$, then we have
\begin{equation}
(k_{1}\pm\partial_{x})u(x)\ge 0,~~\forall x\in\mathbb{R}.
\label{2.4}
\end{equation}
\end{Lem}

\begin{Lem}[Positivity; See \cite{MR2249792}]\label{Lemma 2.2}
Let $w(x)=(k_{1}\pm\partial_{x})u(x)$. Assume that $u\in H^{1}(\mathbb{R})$ with $y=(1-\partial^{2}_{x})u\in\mathcal{M}^{+}(\mathbb{R})$. If $k_{1}\ge 1$ and $k_{2}\ge 2$, then we have
\begin{equation}
(k_{2}\pm\partial_{x})(4-\partial^{2}_{x})^{-1}w(x)\ge 0,~~\forall x\in\mathbb{R}.
\label{2.5}
\end{equation}
\end{Lem}

\section{Stability of a single peakon}\label{Section 3}
The proof  of Lin and Liu in \cite{MR2460268} is not entirely suitable for our work, because it involves all local extrema of the function $v=(4-\partial^{2}_{x})^{-1}u$ on $\mathbb{R}$, and thus  is not  local. For our work, we have  to localize the  estimates. Therefore, we need to modify  a little the proof of Lin and Liu. We do this first for a single peakon.

\begin{Theo}[Stability of Peakons]\label{Theorem 3.1}
Let $u\in \mathcal{X}([0,T[)$, with $0<T\le+\infty$, be a solution of the DP equation and $\varphi_{c}$ be the peakon defined in \eqref{1.7}, traveling to the right at the speed $c>0$. There exist $C>0$ and $\varepsilon_{0}>0$ only depending  on the speed $c>0$, such that if
\begin{equation}
y_{0}=(1-\partial^{2}_{x})u_{0}\in\mathcal{M}^{+}(\mathbb{R})
\label{3.1}
\end{equation}
and
\begin{equation}
\|u_{0}-\varphi_{c}\|_{\mathcal{H}}\le\varepsilon^{2},~~\text{with}~~0<\varepsilon<\varepsilon_{0},
\label{3.2}
\end{equation}
then 
\begin{equation}
\|u(t,\cdot)-\varphi_{c}(\cdot-\xi_{1}(t))\|_{\mathcal{H}}\le C\sqrt{\varepsilon},~~\forall t\in[0,T[,
\label{3.3}
\end{equation}
where $\xi_{1}(t)\in\mathbb{R}$ is any point where the function $v(t,\cdot)=(4-\partial^{2}_{x})^{-1}u(t,\cdot)$ attains its maximum.
\end{Theo}
To prove this theorem we first need the following lemma that enables to control the distance of $E(u)$ and $F(u)$ to respectively $E(\varphi_c)$ and $F(\varphi_c)$.

\begin{Lem}[Control of Distances Between Energies]\label{Lemma 3.1}
Let $u\in H^{1}(\mathbb{R})$ with $y=(1-\partial^{2}_{x})u\in\mathcal{M}^{+}(\mathbb{R})$. If 
$\|u-\varphi_{c}\|_{\mathcal{H}}\le\varepsilon^2 $, then
\begin{equation}
|E(u)-E(\varphi_{c})|\le O(\varepsilon^2)
\label{3.4}
\end{equation}
and 
\begin{equation}
|F(u)-F(\varphi_{c})|\le O(\varepsilon^2),
\label{3.5}
\end{equation}
where $O(\cdot)$ only depends on the speed $c$.
\end{Lem}

\textbf{Proof.} For the first estimate, applying triangular inequality, and using that 
$\|u-\varphi_{c}\|_{\mathcal{H}}\le\varepsilon^2$ and $\|\varphi_{c}\|_{\mathcal{H}}=c/\sqrt{3}$, we have  
\begin{align*}
|E(u)-E(\varphi_{c})|&=|\left\|u\|_{\mathcal{H}}-\|\varphi_{c}\|_{\mathcal{H}}\right|\left(\|u\|_{\mathcal{H}}+\|\varphi_{c}\|_{\mathcal{H}}\right)\\
&\le \|u-\varphi_{c}\|_{\mathcal{H}}\left(\|u-\varphi_{c}\|_{\mathcal{H}}+2\|\varphi_{c}\|_{\mathcal{H}}\right)\\
&\le\varepsilon^2 \left(\varepsilon^2+\frac{2c}{\sqrt{3}}\right)\\
&\le O(\varepsilon^2).
\end{align*}

For the second estimate, applying the H\"older inequality, and using that $\|u-\varphi_{c}\|_{\mathcal{H}}\le\varepsilon^2$ and \eqref{2.3}, we have
\begin{align*}
|F(u)-F(\varphi_{c})|& \le\int_{\mathbb{R}}\left|u^{3}-\varphi^{3}_{c}\right|\\
&\le\int_{\mathbb{R}}|u-\varphi_{c}|(u^{2}+u\varphi_{c}+\varphi^{2}_{c})\\
&\le\|u-\varphi_{c}\|_{L^{2}(\mathbb{R})}\left(\int_{\mathbb{R}}(u^{2}+u\varphi_{c}+\varphi^{2}_{c})^{2}\right)^{1/2}\\
&=\|u-\varphi_{c}\|_{L^{2}(\mathbb{R})}\left(\int_{\mathbb{R}}(u^{4}+2u^{3}\varphi_{c}+3u^{2}\varphi^{2}_{c}+2u\varphi^{3}_{c}+\varphi^{4}_{c})\right)^{1/2}\\
&\le\|u-\varphi_{c}\|_{L^{2}(\mathbb{R})}\\
&\hspace{1cm}\cdot\left(4C^{2}_{S}\|u\|^{4}_{L^{2}(\mathbb{R})}+4cC_{S}\|u\|^{2}_{L^{2}(\mathbb{R})}
+3c^{2}\|u\|^{2}_{L^{2}(\mathbb{R})}+\frac{8}{3}c^{3}C_{S}\|u\|_{L^{2}(\mathbb{R})}+\frac{1}{2}c^{4}\right)^{1/2}\\
&\le O(\varepsilon^2),
\end{align*}
where we also use that the $L^{2}$ norm of $u$ is bounded and the following measures of peakon:
$$\|\varphi_{c}\|_{L^{\infty}(\mathbb{R})}=c,~~\|\varphi_{c}\|_{L^{3}(\mathbb{R})}=\sqrt[3]{\frac{2}{3}}c~~
\text{and}~~\|\varphi_{c}\|_{L^{4}(\mathbb{R})}=\frac{1}{\sqrt[4]{2}}c.$$
This proves the lemma.
\hfill $ \square $ \vspace*{2mm}

Now, to prove Theorem \ref{Theorem 3.1}, by the conservation of $ E(\cdot) $, $ F(\cdot) $ and the continuity of the map $t\mapsto u(t)$ from $[0,T[$ to $H^1(\mathbb{R}) \hookrightarrow \mathcal{H}$ (since $\mathcal{H}\simeq L^{2}$ and $\|u\|_{L^{2}(\mathbb{R})}\le\|u\|_{H^{1}(\mathbb{R})}$), it suffices to prove that for any function $u\in H^{1}(\mathbb{R})$ satisfying $y=(1-\partial^{2}_{x})u\in\mathcal{M}^{+}(\mathbb{R})$, \eqref{3.4} and \eqref{3.5}, if 
\begin{equation}
\inf_{\xi\in\mathbb{R}}\|u-\varphi_{c}(\cdot-\xi)\|_{\mathcal{H}}\le\varepsilon^{1/4},
\label{3.6}
\end{equation}
then 
\begin{equation}
\|u-\varphi_{c}(\cdot-\xi_1)\|_{\mathcal{H}}\le C\sqrt{\varepsilon},
\label{3.7}
\end{equation}
where $\xi_1\in\mathbb{R}$ is any point where the function $v=(4-\partial^{2}_{x})^{-1}u$ attains its maximum.

\begin{figure}[h]
\centering
\subfloat[$\varphi(x)$ and $\rho(x)$ profiles.]
{\includegraphics[width=8cm, height=6cm]{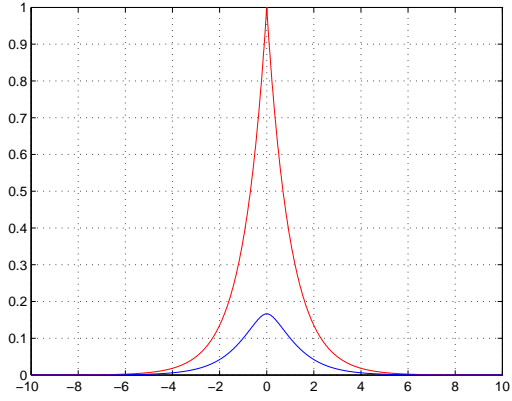}\label{fig1}} \\
\subfloat[$\varphi'(x)$ and  $\rho'(x)$ profiles.]
{\includegraphics[width=8cm, height=6cm]{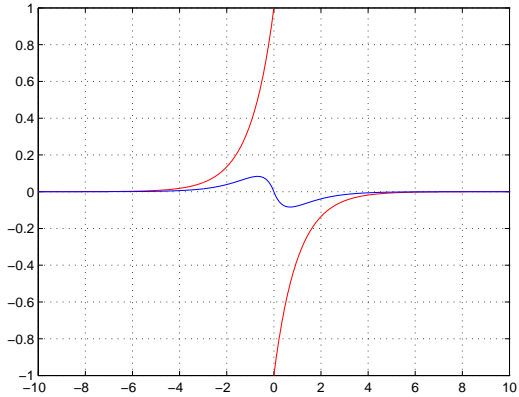}\label{fig2}}
\subfloat[$(1-\partial^{2}_{x})\varphi(x)=2\delta_{0}$ and $(1-\partial^{2}_{x})\rho(x)=(1/2)e^{-2|x|}$ profiles]
{\includegraphics[width=8cm, height=6cm]{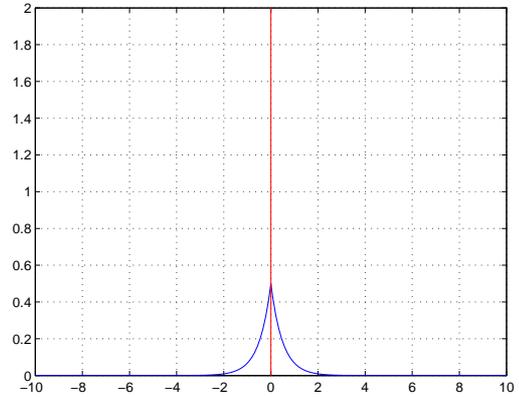}\label{fig3}}
\caption{Variation of peakon and smooth-peakon at initial time with the speed $c=1$.}
\end{figure}

We divide the proof of Theorem \ref{Theorem 3.1} into a sequence of lemmas. In the sequel, we will need to introduce the following \textit{smooth-peakons} defined for all $x\in\mathbb{R}$ by:
\begin{equation}
\rho_{c}(x)=c\rho(x)=(4-\partial^{2}_{x})^{-1}\varphi_{c}(x)=\frac{c}{3}e^{-|x|}-\frac{c}{6}e^{-2|x|}.
\label{3.8}
\end{equation}
One can check that $\rho_{c}\in H^{3}(\mathbb{R})\hookrightarrow C^{2}(\mathbb{R})$ (by the Sobolev embedding) since $\varphi_{c}\in H^{1}(\mathbb{R})$. Indeed, we have 
\begin{equation}
\|\rho_{c}\|^{2}_{H^{3}(\mathbb{R})}=\int_{\mathbb{R}}\frac{(1+\omega^{2})^{3}}{(4+\omega^{2})^{2}}
|\widehat{\varphi}_{c}(\omega)|^{2}
\le\int_{\mathbb{R}}(1+\omega^{2})|\widehat{\varphi}_{c}(\omega)|^{2}
\le \|\varphi_{c}\|^{2}_{H^{1}(\mathbb{R})}
=2c^{2}.
\label{a.a}
\end{equation}
Moreover, $\rho_{c}$ is a positive even function which decays to $0$ at infinity, and admits a single maximum $c/6$ at point $0$ (see Fig. \ref{fig1}-\ref{fig3}).

\begin{Lem}[Uniform Estimates]\label{Lemma 3.2}
Let $u\in H^{1}(\mathbb{R})$ with $y=(1-\partial^{2}_{x})u\in\mathcal{M}^{+}(\mathbb{R})$,
and $\xi\in\mathbb{R}$. If $\|u-\varphi_{c}(\cdot-\xi)\|_{\mathcal{H}}\le \varepsilon^{1/4}$,
then
\begin{equation}
\|u-\varphi_{c}(\cdot-\xi)\|_{L^{\infty}(\mathbb{R})}\le O(\varepsilon^{1/8})
\label{3.10}
\end{equation}
and 
\begin{equation}
\|v-\rho_{c}(\cdot-\xi)\|_{L^{\infty}(\mathbb{R})}\le O(\varepsilon^{1/4}),
\label{3.11}
\end{equation}
where $v=(4-\partial^{2}_{x})^{-1}u$ and $\rho_{c}$ is defined in \eqref{3.8}.
\end{Lem}
\textbf{Proof.}
For the second estimate, applying the H\"older inequality and using assumption, we get for all $x\in\mathbb{R}$,
\begin{align*}
\left|v(x)-\rho_{c}(x-\xi)\right|&\le\frac{1}{4}\int_{\mathbb{R}}e^{-2|x'|}
\left|u(x-x')-\varphi_{c}\left[(x-x')-\xi\right]\right|dx'\\
&\le\frac{1}{4}\left(\int_{\mathbb{R}}e^{-4|x'|}dx'\right)^{1/2}
\left(\int_{\mathbb{R}}\left|u(x')-\varphi_{c}(x'-\xi)\right|^{2}dx'\right)^{1/2}\\
&\le\frac{1}{2\sqrt{2}}\|u-\varphi_{c}(\cdot-\xi)\|_{\mathcal{H}}\\
&\le O(\varepsilon^{1/4}).  
\end{align*}

For the first estimate,  note that the assumption $y=(1-\partial^{2}_{x})u\ge 0$ implies that $u=(1-\partial^{2}_{x})^{-1}y\ge 0$ and satisfies \eqref{2.b}. Then, applying triangular inequality, and using that $|\varphi'_{c}|=\varphi_{c}$ on $\mathbb{R}$ and \eqref{2.3}, we have 
\begin{align*}
\|u-\varphi_{c}(\cdot-\xi)\|_{H^{1}(\mathbb{R})}&\le\|u\|_{H^{1}(\mathbb{R})}
+\|\varphi_{c}\|_{H^{1}(\mathbb{R})}\\
&\le 2\|u\|_{L^{2}(\mathbb{R})}+2\|\varphi_{c}\|_{L^{2}(\mathbb{R})}\\
&\le 2\|u-\varphi_{c}(\cdot-\xi)\|_{L^{2}(\mathbb{R})}+4\|\varphi_{c}\|_{L^{2}(\mathbb{R})}\\
&\le O(\varepsilon^{1/4})+O(1).
\end{align*}
Now, applying the Gagliardo-Nirenberg inequality and using assumption, we obtain
\begin{align*}
\|u-\varphi_{c}(\cdot-\xi)\|_{L^{\infty}(\mathbb{R})}&\le 
C_{G}\|u-\varphi_{c}(\cdot-\xi)\|^{1/2}_{L^{2}(\mathbb{R})}
\|u-\varphi_{c}(\cdot-\xi)\|^{1/2}_{H^{1}(\mathbb{R})}\\
&\le O(\varepsilon^{1/8})\left(O(\varepsilon^{1/8})+O(1)\right)\\
&\le O(\varepsilon^{1/8}).
\end{align*}
This proves the lemma.
\hfill $ \square $ \vspace*{2mm}

\begin{Lem}[Quadratic Identity; See \cite{MR2460268}]\label{Lemma 3.3}
For any $u\in L^{2}(\mathbb{R})$ and $\xi\in\mathbb{R}$, it holds
\begin{equation}
E(u)-E(\varphi_{c})=\|u-\varphi_{c}(\cdot-\xi)\|^{2}_{\mathcal{H}}+4c\left(v(\xi)-\frac{c}{6}\right),
\label{3.12}
\end{equation}
where $v=(4-\partial^{2}_{x})^{-1}u$ and $c/6=\rho_{c}(0)=\max_{x\in\mathbb{R}}\rho_{c}(x-\xi)$.
\end{Lem}

\textbf{Sketch of proof.} The proof follows by direct computation, with the aid of two integration by parts, and using that $(1-\partial^{2}_{x})\varphi_{c}(\cdot-\xi)=2c\delta_{\xi}$, where $\delta_{\xi}$ denotes the Dirac mass applied at point $\xi$.
\hfill $ \square $ \vspace*{2mm}

Let $u\in H^{1}(\mathbb{R})$ with $y = (1-\partial^{2}_{x})u\in\mathcal{M}^{+}(\mathbb{R})$, and assume that there exists $\xi\in\mathbb{R}$ such that \eqref{3.6} holds for some $\xi\in\mathbb{R}$. We consider now the interval in which the peakon (respectively the smooth-peakon) is concentrated, and we will decompose this interval according to the variation of $v=(4-\partial^{2}_{x})^{-1}u$ in the following way: we set 
\begin{equation}
\alpha=\sup\left\lbrace x<\xi,~v(x)=\frac{c}{2400}\right\rbrace~~\text{and}~~
\beta=\inf\left\lbrace x>\xi,~v(x)=\frac{c}{2400}\right\rbrace.
\label{3.13}
\end{equation}
According to Lemma \ref{Lemma 3.2}, we know that $v$ is close to $\rho_{c}(\cdot-\xi)$ in $L^{\infty}$ norm with  $\rho_{c}(0)=c/6$. Therefore $v$ must have at least one local maximum on $[\alpha,\beta]$. Assume that on $[\alpha,\beta]$ the function $v$ admits $k+1$ points $(\xi_{j})_{j=1}^{k+1}$ with local maximal values for some integer $k\ge 0$, where $\xi_{1}$ is the first local maximum point and $\xi_{k+1}$ the last local maximum point\footnote{In the case of an infinite countable number of local maximal values, the proof is exactly the same.}. Then between $\xi_{1}$ and $\xi_{k+1}$, the function $v$ admits  $k$ points $(\eta_{j})_{j=1}^{k}$ with local minimal values. We rename $\alpha=\eta_{0}$ and $\beta=\eta_{k+1}$ so that it holds 
\begin{equation}
\eta_{0}<\xi_{1}<\eta_{1}<\ldots<\xi_{j}<\eta_{j}<\xi_{j+1}<\eta_{j+1}
<\ldots<\eta_{k}<\xi_{k+1}<\eta_{k+1}.
\label{3.14}
\end{equation}
Let 
\begin{equation}
M_{j}=v(\xi_{j}),~~j=1,\ldots,k+1,~~\text{and}~~m_{j}=v(\eta_{j}),~~j=1,\ldots,k.
\label{3.15}
\end{equation}
By construction
\begin{equation}
v_{x}(x)\ge 0,~~\forall x\in[\eta_{j-1},\xi_{j}],~~j=1,\ldots,k
\label{3.16}
\end{equation}
and
\begin{equation}
v_{x}(x)\le 0,~~\forall x\in[\xi_{j},\eta_{j}],~~j=1,\ldots,k+1.
\label{3.17}
\end{equation}
We claim that 
\begin{equation}
v(x)\le\frac{c}{300},~~\forall x\in\mathbb{R}\setminus[\eta_{0},\eta_{k+1}],
\label{3.18}
\end{equation}
\begin{equation}
u(x)\le\frac{c}{300},~~\forall x\in\mathbb{R}\setminus[\eta_{0},\eta_{k+1}],
\label{3.19}
\end{equation}
and there exits $C_{0}>0$ such that 
\begin{equation}
[\eta_{0},\eta_{k+1}]\subset[\xi-C_{0},\xi+C_{0}].
\label{3.20}
\end{equation}
Indeed, for some $0<\varepsilon\ll 1$ fixed, using \eqref{3.11} we have 
$$\rho_{c}(\eta_{0}-\xi)=v(\eta_{0})+O(\varepsilon^{1/4})=\frac{c}{2400}+O(\varepsilon^{1/4})\le\frac{c}{350}.$$
Please note that, we abuse notation by writing that the difference between $v$ and $\rho_{c}(\cdot-\xi)$  is equal to $O(\varepsilon^{1/4})$. Therefore, using that $\rho_{c}(\cdot-\xi)$ is increasing on $]-\infty,\xi]$, it holds for all $x\in]-\infty,\eta_{0}[$,
$$v(x)=\rho_{c}(x-\xi)+O(\varepsilon^{1/4})
\le\frac{c}{350}+O(\varepsilon^{1/4})
\le\frac{c}{300}.$$
Proceeding in the same way for $x\in]\eta_{k+1},+\infty[$, we obtain \eqref{3.18}.

One can remark that for all $x\in\mathbb{R}$,
\begin{equation}
\varphi_{c}(x)-6\rho_{c}(x)=ce^{-|x|}-6\left(\frac{c}{3}e^{-|x|}-\frac{c}{6}e^{-2|x|}\right)
=-ce^{-|x|}+ce^{-2|x|}
\le 0.
\label{3.21}
\end{equation}
Thus, combining \eqref{3.10}, \eqref{3.21} and proceeding as for the estimate \eqref{3.18}, we infer \eqref{3.19}.

Finally, from \eqref{3.11} we have 
$$\rho_{c}(\eta_{0}-\xi)=v(\eta_{0})+O(\varepsilon^{1/4})
=\frac{c}{2400}+O(\varepsilon^{1/4})
\ge\frac{c}{3000}.$$
Therefore, since $\rho_{c}=(c/3)e^{-|\cdot|}-(c/6)e^{-2|\cdot|}$ and that  $x\mapsto (1/3)e^{-|x|}-(1/6)e^{-2|x|}$ is a positive even function decreasing to $0$ on $\mathbb{R}_{+}$ (see Fig. \ref{fig1}), there exists a universal constant $C_{0}>0$ such that \eqref{3.20} holds.

We now are ready to establish the connection between the conservation laws. Please note that, we will change the order of the extrema of $v=(4-\partial^{2}_{x})^{-1}u$ while keeping the same notations as in \eqref{3.15}.

\begin{Lem}[Connection Between $E(\cdot)$ and the Local Extrema of $v$]\label{Lemma 3.4}
Let $u\in H^{1}(\mathbb{R})$ and $v=(4-\partial^{2}_{x})^{-1}u\in H^{3}(\mathbb{R})$. Define the function $g$ by
\begin{equation}
  g(x)=\left\{
    \begin{aligned}
     &2v+v_{xx}-3v_{x},~~x<\xi_{1},\\
     &2v+v_{xx}+3v_{x},~~\xi_{j}<x<\eta_{j},\\
     &2v+v_{xx}-3v_{x},~~\eta_{j}<x<\xi_{j+1},\\
     &2v+v_{xx}+3v_{x},~~x>\xi_{k+1},\\
    \end{aligned}
  \right.~~j=1,\ldots,k.
  \label{3.22}
\end{equation}
Then it holds
\begin{equation}
\int_{\mathbb{R}}g^{2}(x)dx=E(u)-12\left(\sum_{j=0}^{k}M^{2}_{j+1}-\sum_{j=1}^{k}m^{2}_{j}\right).
\label{3.23}
\end{equation}
\end{Lem}

\textbf{Proof.} We have 
\begin{equation}
\int_{\mathbb{R}}g^{2}(x)dx=\int_{-\infty}^{\xi_{1}}g^{2}(x)dx
+\sum_{j=1}^{k}\int_{\xi_{j}}^{\xi_{j+1}}g^{2}(x)dx
+\int_{\xi_{k+1}}^{+\infty}g^{2}(x)dx.
\label{3.24}
\end{equation}
For $j=1,\ldots,k$,
\begin{align*}
\int_{\xi_{j}}^{\xi_{j+1}} g^{2}(x)(x)dx&=\int_{\xi_{j}}^{\eta_{j}}\left( 2v+v_{xx}+3v_{x}\right)^{2}+\int_{\eta_{j}}^{\xi_{j+1}}\left(2v+v_{xx}-3v_{x}\right)^{2}\\
&=J+I.
\end{align*}
Let us compute $I$,
\begin{align*}
I &=\int_{\eta_{j}}^{\xi_{j+1}}\left(4v^{2}+v^{2}_{xx}+9v^{2}_{x}+4vv_{xx}
-12vv_{x}-6v_{x}v_{xx}\right)\\
&=\int_{\eta_{j}}^{\xi_{j+1}}\left(4v^{2}+v^{2}_{xx}+9v^{2}_{x}\right)+
4\int_{\eta_{j}}^{\xi_{j+1}}vv_{xx}-12\int_{\eta_{j}}^{\xi_{j+1}}vv_{x}
-6\int_{\eta_{j}}^{\xi_{j+1}}v_{x}v_{xx}\\
&=\int_{\eta_{j}}^{\xi_{j+1}}\left(4v^{2}+v^{2}_{xx}+9v^{2}_{x}\right)+I_{1}+I_{2}
+I_{3}.
\end{align*}
Applying integration by parts and using that $v_{x}(\xi_{j})=v_{x}(\eta_{j})=0$, we get
$$I_{1}=-4\int_{\eta_{j}}^{\xi_{j+1}}v^{2}_{x},~~I_{2}=-6\int_{\eta_{j}}^{\xi_{j+1}}\partial_{x}(v^{2})\nonumber=-6v^{2}(\xi_{j+1})+6v^{2}(\eta_{j})~~\text{and}~~I_{3}=-3\int_{\eta_{j}}^{\xi_{j+1}}\partial_{x}(v^{2}_{x})=0.$$
Therefore
\begin{equation}\label{0a1}
I=\int_{\eta_{j}}^{\xi_{j+1}}\left(4v^{2}+5v^{2}_{x}+v^{2}_{xx}\right)-6v^{2}(\xi_{j+1})+6v^{2}(\eta_{j}).
\end{equation}
Similar computations lead to 
\begin{equation}
J=\int_{\xi_{j}}^{\eta_{j}}\left(4v^{2}+5v^{2}_{x}+v^{2}_{xx}\right)-6v^{2}(\xi_{j})+6v^{2}(\eta_{j}),~~
\int_{-\infty}^{\xi_1} g^2(x) \, dx =\int_{-\infty}^{\xi_1}\left(4v^{2}+5v^{2}_{x}+v^{2}_{xx}\right)-6v^{2}(\xi_{1}) 
\label{3.26}
\end{equation}
and 
\begin{equation}
 \int_{\xi_{k+1}}^{+\infty} g^2(x) \, dx =\int_{\xi_{k+1}}^{+\infty} \left(4v^{2}+5v^{2}_{x}+v^{2}_{xx}\right)-6v^{2}(\xi_{k+1}).
\label{3.2666}
\end{equation}
 Adding $I$ and $J$, and summing over $j\in\lbrace 1,\ldots, k\rbrace$, we obtain
\begin{equation}
\int_{\xi_{1}}^{\xi_{k+1}} g^{2}(x)dx=
\int_{\xi_{1}}^{\xi_{k+1}}\left(4v^{2}+5v^{2}_{x}+v^{2}_{xx}\right)
-6\sum_{j=1}^{k}v^{2}(\xi_{j+1})
-6\sum_{j=1}^{k}v^{2}(\xi_{j})
+12\sum_{j=1}^{k}v^{2}(\eta_{j}).
\label{3.27}
\end{equation}
The lemma follows by combining \eqref{3.24} and \eqref{3.26}-\eqref{3.27}.
\hfill $ \square $ \vspace*{2mm}

\begin{Lem}[Connection Between $F(\cdot)$ and the Local Extrema of $v$]\label{Lemma 3.5}
Let $u\in H^{1}(\mathbb{R})$ and $v=(4-\partial^{2}_{x})^{-1}u\in H^{3}(\mathbb{R})$. Define the function $h$ by
\begin{equation}
  h(x)=\left\{
    \begin{aligned}
     &-v_{xx}-6v_{x}+16v,~~x<\xi_{1},\\
     &-v_{xx}+6v_{x}+16v,~~\xi_{j}<x<\eta_{j},\\
     &-v_{xx}-6v_{x}+16v,~~\eta_{j}<x<\xi_{j+1},\\
     &-v_{xx}+6v_{x}+16v,~~x>\xi_{k+1},\\
    \end{aligned}
  \right.~~j=1,\ldots,k.
  \label{3.28}
\end{equation}
Then it holds
\begin{equation}
\int_{\mathbb{R}}h(x)g^{2}(x)dx=F(u)-144\left(\sum_{j=0}^{k}M^{3}_{j+1}-\sum_{j=1}^{k}m^{3}_{j}\right).
\label{3.29}
\end{equation}
\end{Lem}

\textbf{Proof.} We have 
\begin{equation}
\int_{\mathbb{R}}h(x)g^{2}(x)dx=\int_{-\infty}^{\xi_{1}}
h(x)g^{2}(x)dx
+\sum_{j=1}^{k}\int_{\xi_{j}}^{\xi_{j+1}}h(x)g^{2}(x)dx+\int_{\xi_{k+1}}^{+\infty}h(x)g^{2}(x)dx.
\label{3.30}
\end{equation}
For $j=1,\ldots,k$, 
\begin{align*}
\int_{\xi_{j}}^{\xi_{j+1}}h(x)g^{2}(x)dx&=\int_{\xi_{j}}^{\eta_{j}}\left(-v_{xx}-6v_{x}+16v\right)\left(2v+v_{xx}
-3v_{x}\right)^{2}\\
&\hspace{1cm}+\int_{\eta_{j}}^{\xi_{j+1}}\left(-v_{xx}+6v_{x}+16v\right)\left(2v+v_{xx}
+3v_{x}\right)^{2}\\
&=J+I.
\end{align*}
Let us compute $I$,
\begin{align*}
I&=\int_{\eta_{j}}^{\xi_{j+1}}\left(-v^{3}_{xx}+12vv^{2}_{xx}
+64v^{3}+60v^{2}v_{xx}\right)-54\int_{\eta_{j}}^{\xi_{j+1}}v^{3}_{x}
+27\int_{\eta_{j}}^{\xi_{j+1}}v^{2}_{x}v_{xx}\\
&\hspace{1cm}-108\int_{\eta_{j}}^{\xi_{j+1}}vv_{x}v_{xx}
-216\int_{\eta_{j}}^{\xi_{j+1}}v^{2}v_{x}
+216\int_{\eta_{j}}^{\xi_{j+1}}vv^{2}_{x}\\
&=\int_{\eta_{j}}^{\xi_{j+1}}\left(-v^{3}_{xx}+12vv^{2}_{xx}
+64v^{3}+60v^{2}v_{xx}\right)\\
&\hspace{1cm}-54\int_{\eta_{j}}^{\xi_{j+1}}v^{3}_{x}
+I_{1}+I_{2}+I_{3}+I_{4}.
\end{align*}
Applying integration by parts and using that $v_{x}(\xi_{j})=v_{x}(\eta_{j})=0$, we get
$$I_{1}=9\int_{\eta_{j}}^{\xi_{j+1}}\partial_{x}(v^{3}_{x})=0,~~
I_{2}=54\int_{\eta_{j}}^{\xi_{j+1}}v^{3}_{x},~~
I_{3}=-72\int_{\eta_{j}}^{\xi_{j+1}}\partial_{x}(v^{3})=-72v^{3}(\xi_{j+1})+72v^{3}(\eta_{j})$$
and
\begin{align*}
I_{4}&=108\int_{\eta_{j}}^{\xi_{j+1}}\partial_{x}(v^{2})v_{x}=-108\int_{\eta_{j}}^{\xi_{j+1}}v^{2}v_{xx}.
\end{align*}
Therefore
\begin{equation}
I=\int_{\eta_{j}}^{\xi_{j+1}}\left(-v^{3}_{xx}+12vv^{2}_{xx}
+64v^{3}+60v^{2}v_{xx}\right)-72v^{3}(\xi_{j+1})+72v^{3}(\eta_{j}).
\label{0a2}
\end{equation}
Similar computations lead to 
\begin{equation}
J=\int_{\xi_{j}}^{\eta_{j}}\left(-v^{3}_{x}+12vv^{2}_{xx}
+64v^{3}+60v^{2}v_{xx}\right)-72v^{3}(\xi_{j})
+72v^{3}(\eta_{j}),
\label{0a3}
\end{equation}
\begin{equation}
\int_{-\infty}^{\xi_1} h(x) g^2(x) \\ dx =\int_{-\infty}^{\xi_1}\left(-v^{3}_{x}+12vv^{2}_{xx}
+64v^{3}+60v^{2}v_{xx}\right)-72 v^{3}(\xi_{1})
\label{3.322}
\end{equation}
and 
\begin{equation}
 \int_{\xi_{k+1}}^{+\infty} h(x) g^2(x) \, dx =\int_{\xi_{k+1}}^{+\infty} \left(-v^{3}_{x}+12vv^{2}_{xx}
+64v^{3}+60v^{2}v_{xx}\right)-72 v^{3}(\xi_{k+1}).
\label{3.3222}
\end{equation}
Adding \eqref{0a1} and \eqref{0a3}, and summing over $j\in\lbrace 1,\ldots,k\rbrace$, we obtain
\begin{align}
\int_{\xi_{1}}^{\xi_{k+1}}h(x)g^{2}(x)dx&=
\int_{\xi_{1}}^{\xi_{k+1}}\left(-v^{3}_{x}+12vv^{2}_{xx}
+64v^{3}+60v^{2}v_{xx}\right)-72\sum_{j=1}^{k}v^{3}(\xi_{j+1})\nonumber\\
&\hspace{1cm}-72\sum_{j=1}^{k}v^{3}(\xi_{j})
+72\sum_{j=1}^{k}v^{3}(\eta_{j}).
\label{3.33}
\end{align}
The lemma follows by combining \eqref{3.30} and \eqref{3.322}-\eqref{3.33}.
\hfill $ \square $ \vspace*{2mm}

\begin{Lem}[Connection Between $E(\cdot)$ and $F(\cdot)$]\label{Lemma 3.6}
Let $u\in H^{1}(\mathbb{R})$, with  $y=(1-\partial^{2}_{x})u \in {\mathcal M}^+({\mathbb R}) $, that satisfies  \eqref{3.6} for some $\xi\in\mathbb{R}$. Assume that 
$v=(4-\partial^{2}_{x})^{-1}u$ satisfies \eqref{3.13}-\eqref{3.20}, with local extrema on $[\eta_{0},\eta_{k+1}]$ arranged in decreasing order in the following way:
\begin{equation}
M_{1}\ge M_{2}\ge\ldots\ge M_{k+1}\ge 0,~~m_{1}\ge m_{2}\ge\ldots\ge m_{k}\ge 0,~~M_{j+1}\ge m_{j},~j=1,\ldots,k.
\label{3.34}
\end{equation}
There exists $\varepsilon_{0}>0$ only depending on the speed $c$, such that if $0<\varepsilon<\varepsilon_{0}$, then it holds
\begin{equation}
M^{3}_{1}-\frac{1}{4}E(u)M_{1}+\frac{1}{72}F(u)\le 0.
\label{3.35}
\end{equation}
\end{Lem}

\textbf{Proof.}
The key is to show that $h\le 18 M_{1}$ on $\mathbb{R}$. Note that by \eqref{3.6} we know that $18M_{1}\ge c/4$.
We rewrite the function $h$ as
\begin{equation*}
  h(x)=\left\{
    \begin{aligned}
     &-v_{xx}-6v_{x}+16v,~~x<\eta_{0},\\
     &-\left(\partial^{2}_{x}+3\partial_{x}+2\right)v-3v_{x}+18v,~~\eta_{0}<x<\xi_{1},\\
     &-\left(\partial^{2}_{x}-3\partial_{x}+2\right)v+3v_{x}+18v,~~\xi_{j}<x<\eta_{j},\\
     &-\left(\partial^{2}_{x}+3\partial_{x}+2\right)v-3v_{x}+18v,~~\eta_{j}<x<\xi_{j+1},\\
     &-\left(\partial^{2}_{x}-3\partial_{x}+2\right)v+3v_{x}+18v,~~\xi_{k+1}<x<\eta_{k+1},\\
     &-v_{xx}+6v_{x}+16v,~~x>\eta_{k+1},\\
    \end{aligned}
  \right.~~j=1,\ldots,k.
  \end{equation*}
First, one can remark that for all $x\in\mathbb{R}$,
$$v(x)=\frac{e^{-2x}}{4}\int_{-\infty}^{x}e^{2x'}u(x')dx'+\frac{e^{2x}}{4}\int_{x}^{+\infty}e^{-2x'}u(x')dx'$$
and 
$$v_{x}(x)=-\frac{e^{-2x}}{2}\int_{-\infty}^{x}e^{2x'}u(x')dx'+\frac{e^{2x}}{2}\int_{x}^{+\infty}e^{-2x'}u(x')dx'.$$
Then using that $u=(1-\partial^{2}_{x})^{-1}y\ge 0$ on $\mathbb{R}$, we get  
\begin{equation}
|v_{x}(x)|\le 2v(x),~~\forall x\in\mathbb{R}.
\label{3.36}
\end{equation}
Next, if $x\in\mathbb{R}\setminus[\eta_{0},\eta_{k+1}]$, using that $v_{xx}=4v-u$, \eqref{3.18}, \eqref{3.19} and \eqref{3.36}, it holds 
\begin{align*}
h&\le|v_{xx}|+6|v_{x}|+16v\le u+32v\le\frac{c}{9}.
\end{align*}
If $\eta_{0}<x<\xi_{1}$, then $v_{x}\ge 0$, and using that $y=(1-\partial^{2}_{x})u\ge 0$, it follows from  Lemma \ref{Lemma 2.2} that 
\begin{align*}
h&=-(\partial^{2}_{x}+3\partial_{x}+2)v-3v_{x}+18v\\
&=-(2+\partial_{x})(4-\partial^{2}_{x})^{-1}(1+\partial_{x})u-3v_{x}+18v\\
&\le 18v.
\end{align*}
If $\xi_{j}<x<\eta_{j}$, then $v_{x}\le 0$, and similarly using that $y=(1-\partial^{2}_{x})u\ge 0$, it follows from  Lemma \ref{Lemma 2.2} that 
\begin{align*}
h&=-(\partial^{2}_{x}-3\partial_{x}+2)v+3v_{x}+18v\\
&=-(2-\partial_{x})(4-\partial^{2}_{x})^{-1}(1-\partial_{x})u+3v_{x}+18v\\
&\le 18v.
\end{align*}
Therefore, it holds 
\begin{equation}
h(x)\le 18\max_{x\in\mathbb{R}}v(x)=18M_{1},~~\forall x\in\mathbb{R}.
\label{3.37}
\end{equation}

Now, combining \eqref{3.23}, \eqref{3.29} and \eqref{3.37}, we get 
\begin{align*}
F(u)-144&\left(\sum_{j=0}^{k}M^{3}_{j+1}-\sum_{j=1}^{k}m^{3}_{j}\right)\\
&=\int_{\mathbb{R}}h(x)g^{2}(x)dx\\
&\le\|h\|_{L^{\infty}(\mathbb{R})}\int_{\mathbb{R}}g^{2}(x)dx\\
&\le 18M_{1}\left[E(u)-12\left(\sum_{j=0}^{k}M^{2}_{j+1}-\sum_{j=1}^{k}m^{2}_{j}\right)\right].
\end{align*}
For $j=1,\ldots,k$, we set
$$A_{j}=M^{3}_{j+1}-m^{3}_{j}~~\text{and}~~B_{j}=M^{2}_{j+1}-m^{2}_{j},$$
and our inequality becomes 
\begin{equation}
M^{3}_{1}-\frac{1}{4}E(u)M_{1}+\frac{1}{72}F(u)\le 2\sum_{j=1}^{k}\left(A_{j}-\frac{3}{2}M_{1}B_{j}\right).
\label{3.38}
\end{equation}
On the other hand, using that $M_{j+1}\ge m_{j}$, we have
\begin{equation}
A_{j}-\frac{3}{2}M_{1}B_{j}=-\frac{1}{2}\left(M_{j+1}-m_{j}\right)
\left(3M_{1}m_{j}+3M_{1}M_{j+1}-2M_{j+1}m_{j}-2M^{2}_{j+1}-2m^{2}_{j}\right)\le 0.
\label{3.39}
\end{equation}
Finally, combining \eqref{3.38} and \eqref{3.39}, we obtain the lemma.
\hfill $ \square $ \vspace*{2mm}

\textbf{Proof of Theorem \ref{Theorem 3.1}.}
We argue as El Dika and Molinet in \cite{MR2542735}. As noticed after the statement of the theorem, it suffices to prove \eqref{3.7} assuming that $u \in H^1(\mathbb{R}) $ satisfies \eqref{3.1}, \eqref{3.2} and 
\eqref{3.4}-\eqref{3.6}.  We set $M_{1}=v(\xi_{1})=\max_{x\in\mathbb{R}}v(x)$ and $\delta=c/6-M_{1}$. 
We first remark that if $\delta\le 0$, combining \eqref{3.4} and \eqref{3.12}, it holds
\begin{equation*}
\|u-\varphi_{c}(\cdot-\xi_{1})\|_{\mathcal{H}}\le|E(u_{0})-E(\varphi_{c})|^{1/2}
\le O(\varepsilon),
\end{equation*}
that yields the desired result. 
Now suppose that $\delta>0$, that is the maximum of the function $v$ is less than the maximum of $\rho_{c}$. 
 Combining \eqref{3.4}, \eqref{3.5} and \eqref{3.35}, we get
$$M^{3}_{1}-\frac{1}{4}E(\varphi_{c})M_{1}+\frac{1}{72}F(\varphi_{c})\le O(\varepsilon^{2}).$$
Using that $E(\varphi_{c})=c^{2}/3$ and $F(\varphi_{c})=2c^{3}/3$, our inequality becomes
$$\left(M_{1}-\frac{c}{6}\right)^{2}\left(M_{1}+\frac{c}{3}\right)\le  O(\varepsilon^{2}).$$
Substituting $M_{1}$ by $c/6-\delta$ and using that $(M_{1}+c/3)^{-1}<3/c$, it holds
\begin{equation}
\delta^{2}\le O(\varepsilon^{2})\Rightarrow\delta\le O(\varepsilon).
\label{3.40}
\end{equation}
Finally, combining \eqref{3.4}, \eqref{3.12} and \eqref{3.40}, we obtain
$$\|u-\varphi_{c}(\cdot-\xi_{1})\|_{\mathcal{H}}\le C\sqrt{\varepsilon},$$
where $C>0$ only depends on the speed $c$. This completes the proof of  the stability of a single peakon.

\section{Stability of the trains of peakons}\label{Section 4}
For $\gamma>0$ and $L>0$, we define the following neighborhood of all the sums of $N$ peakons of speed $c_{1},...,c_{N}$ with spatial shifts $z_{i}$ that satisfied $z_{i}-z_{i-1}\ge L$,
\begin{equation}
U(\gamma,L)=\left\lbrace u\in H^{1}(\mathbb{R});~\inf_{z_{i}-z_{i-1}>L}\left\|u-\sum_{i=1}^{N}\varphi_{c_{i}}(\cdot-z_{i})\right\|_{\mathcal{H}}\le\gamma\right\rbrace.
\label{4.1}
\end{equation}
By the continuity of the map $t\mapsto u(t)$ from $[0,T[$ into $H^{1}(\mathbb{R})\hookrightarrow\mathcal{H}$, to prove Theorem $1.1$ it suffices to prove that there exist $A>0$, $\varepsilon_{0}>0$ and $L_{0}>0$ such that for all $ L>L_{0}$ and $0<\varepsilon<\varepsilon_{0}$, if $u_{0}$ satisfies \eqref{1.11}-\eqref{1.13}, and if for some $0<t_{0}<T$,
\begin{equation}
u(t)\in U\left(A(\sqrt{\varepsilon}+L^{-1/8}),\frac{L}{2}\right),~~\forall t\in[0,t_{0}],
\label{4.2}
\end{equation}
then
\begin{equation}
u(t_{0})\in U\left(\frac{A}{2}(\sqrt{\varepsilon}+L^{-1/8}),\frac{2L}{3}\right).
\label{4.3}
\end{equation}
Therefore, in the sequel of this section we will assume \eqref{4.2} for some $0<\varepsilon<\varepsilon_{0}$ and $L>L_{0}$, with $A$, $\varepsilon_{0}$ and $L_{0}$ to be specified later, and we will prove \eqref{4.3}.

\begin{Rem}[Distance Between $v$ and the Sum of $N$ Smooth-peakons]\label{Remark 4.1}
\normalfont
From the definition of $E(\cdot)$ and $\mathcal{H}$ (see respectively \eqref{1.3} and \eqref{n1}), one can clearly see that $\|u\|_{\mathcal{H}}$ is equivalent to $\|v\|_{H^{2}(\mathbb{R})}$, where 
$v=(4-\partial^{2}_{x})^{-1}u$. Let $t_{1}\in[0,t_{0}]$ fixed, if $u(t_{1})\in U(\gamma,L/2)$, 
then there exists $\tilde{Z}=(\tilde{z}_{i})_{i=1}^{N}$ with $\tilde{z}_{i}-\tilde{z}_{i-1}\ge L/2$, such that 
$v(t_{1})\in H^{3}(\mathbb{R})\hookrightarrow C^{2}(\mathbb{R})$  stays close to $\sum_{i=1}^{N}\rho_{c_{i}}
(\cdot-\tilde{z}_{i})$ in the $H^{2}$ norm, where $\rho_{c_{i}}$ is defined in 
\eqref{3.8}. 
\end{Rem}

\subsection{Control of the distance between the peakons}\label{Subsection 4.1}
In this subsection, we want to prove that the different bumps of $u$ (respectively of $v$) that are individually close to a peakon 
(respectively a smooth-peakon) get away from each others as time is increasing. This is crucial in our analysis since we do not know how to manage strong interactions.

\begin{Lem}[Decomposition of the Solution Around $\varphi_{c}$]\label{Lemma 4.1}
Let $u_{0}$ satisfying \eqref{1.11}-\eqref{1.13}. There exist $\gamma_{0}>0$, $L_{0}>0$ and $C_{0}>0$ such that for all $0<\gamma<\gamma_{0}$ and $0<L_{0}<L$, if $u(t)\in U(\gamma,L/2)$ on $[0,t_{0}]$ for some $0<t_{0}<T$, then there exist $N~C^{1}$ functions $\tilde{x}_{1},\ldots,\tilde{x}_{N}$ defined on $[0,t_{0}]$ such that 
\begin{equation}
\left\|u(t)-\sum_{i=1}^{N}\varphi_{c_{i}}(\cdot-\tilde{x}_{i}(t))\right\|_{\mathcal{H}}\le O(\gamma),
\label{4.5}
\end{equation}
\begin{equation}
\left\|v(t)-\sum_{i=1}^{N}\rho_{c_{i}}(\cdot-\tilde{x}_{i}(t))\right\|_{C^{1}(\mathbb{R})}
\le O(\gamma),
\label{4.6}
\end{equation}
\begin{equation}
\left|\dot{\tilde{x}}_{i}(t)-c_{i}\right|\le c_1^{-2}\left(O(\gamma)+O(e^{-L/4})\right) ,~~i=1,\ldots,N,
\label{4.7}
\end{equation}
and
\begin{equation}
\tilde{x}_{i}(t)-\tilde{x}_{i-1}(t)\ge\frac{3L}{4}+\frac{(c_{i}-c_{i-1})t}{2},~~i=2,\ldots,N.
\label{4.8}
\end{equation}
Moreover, for $i=1,\ldots,N$, setting $J_{i}=[y_{i}(t),y_{i+1}(t)]$, with
\begin{equation}
  \left\{
    \begin{aligned}
     &y_{1}=-\infty,\\
     &y_{i}(t)=\frac{\tilde{x}_{i-1}(t)+\tilde{x}_{i}(t)}{2},\\
     &y_{N+1}=+\infty,\\
    \end{aligned}
  \right.~~i=2,\ldots,N,
  \label{4.9}
\end{equation}
it holds
\begin{equation}
\left|\xi^{i}_{1}(t)-\tilde{x}_{i}(t)\right|\le\frac{L}{12},~~i=1,\ldots,N,
\label{4.10}
\end{equation}
where $\xi^{1}_{1}(t),\ldots,\xi^{N}_{1}(t)$ are any point such that
\begin{equation}
v\left(t,\xi^{i}_{1}(t)\right)=\max_{x\in J_{i}}v(t,x),~~i=1,\ldots,N,
\label{4.11}
\end{equation}
and where $v=(4-\partial^{2}_{x})^{-1}u$ and $O(\cdot)$ only depends on the speeds $(c_{i})_{i=1}^{N}$.
\end{Lem}

\textbf{Proof.}
We will slightly modify the construction done by El Dika and Molinet in \cite{MR2542735}. One can remark that the peakons $\varphi_{c_{i}}(\cdot-c_{i}t)$ and the smooth-peakons $\rho_{c_{i}}(\cdot-c_{i}t)$ travel at the same speed $c_{i}$, thanks to this, we will do our construction with $v=(4-\partial^{2}_{x})^{-1}u$ instead of $u$. We do that because the  $\mathcal{H}$ 
(equivalent to $L^{2}$) approximation \eqref{4.2} does not permit us to construct a $C^{1}$ function, which is crucial for application of the Implicit Function Theorem. We note that the same approach can also be used for the CH equation.

\begin{figure}[ht]
\centering
\subfloat[$\left(\int_{0}^{x}|\varphi(x-10)-\varphi(x-10-10^{-9})|^{2}\right)^{1/2}\approx 10^{-9}$]
{\includegraphics[width=8cm, height=6cm]{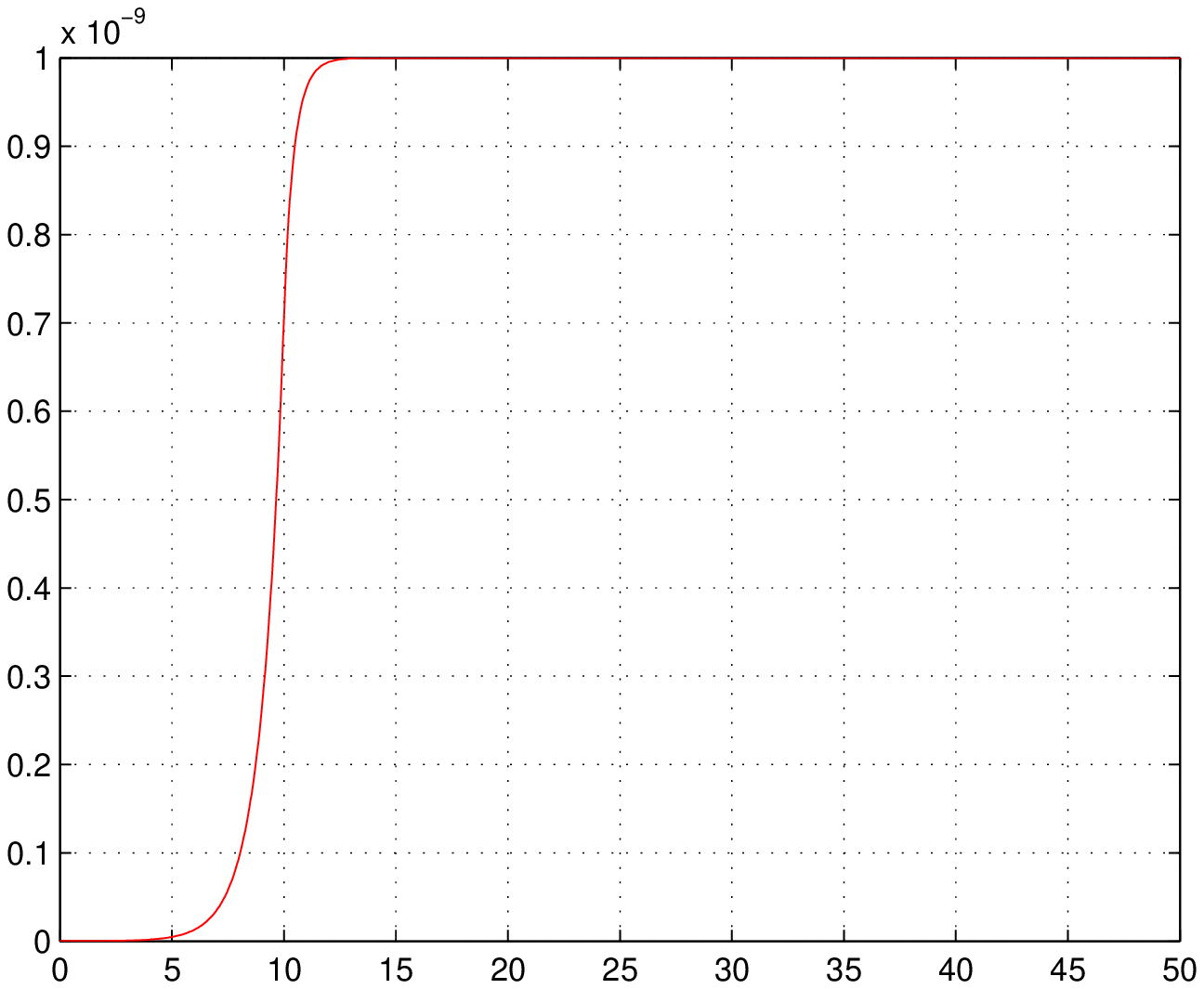}\label{aa1}} 
\subfloat[$\log\left(\int_{0}^{x}|\varphi(x-10)-\varphi(x-10-10^{-9})|^{2}\right)^{1/2}\approx-20.72$]
{\includegraphics[width=8cm, height=6cm]{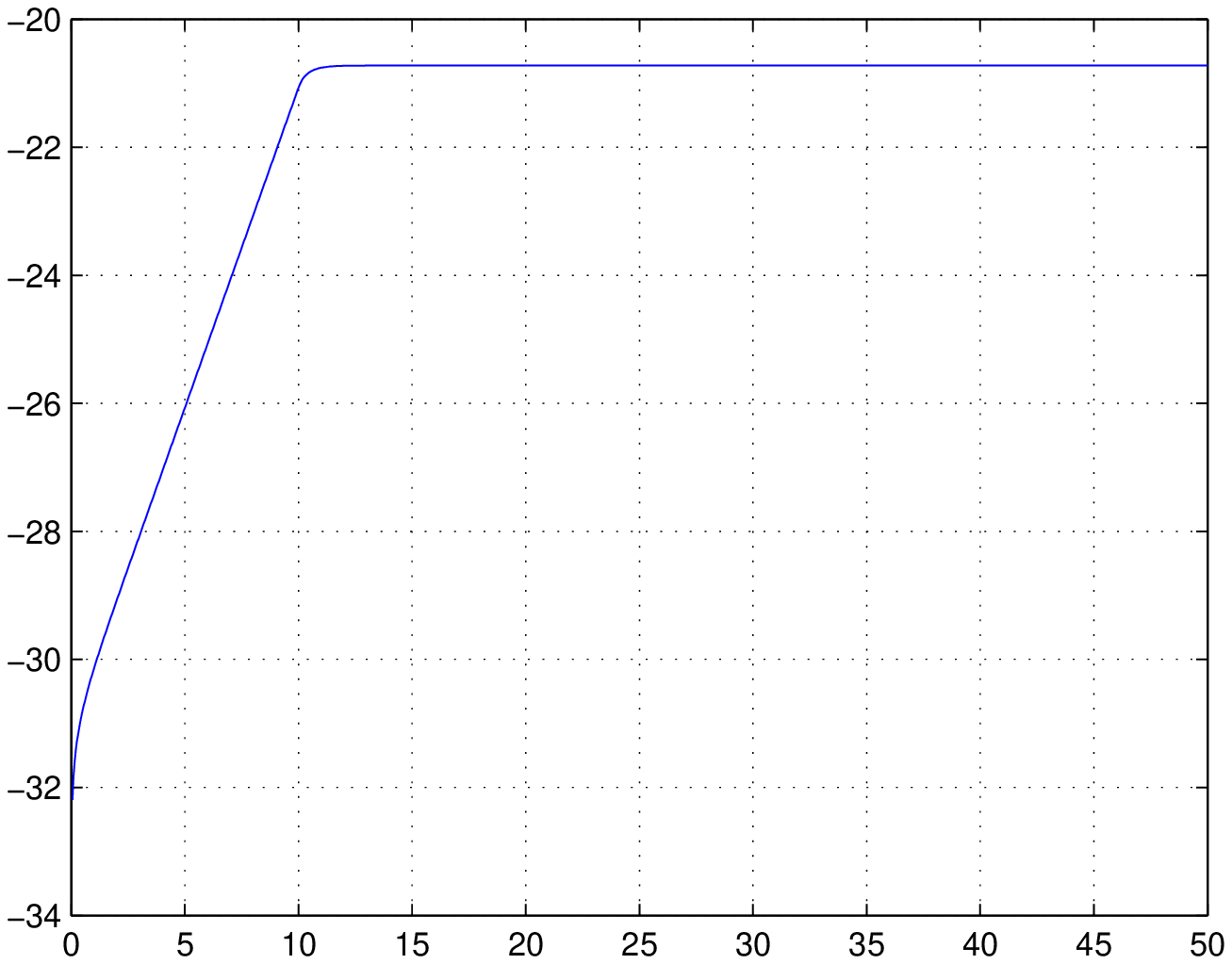}\label{aa2}}
\caption{Distance between two very close peakons (at time $t=10$ with respective speeds $1$ and $1+10^{-10})$.}
\end{figure}

For $Z=(z_{1},\ldots,z_{N})\in\mathbb{R}^{N}$ fixed such that $|z_{i}-z_{i-1}|>L/2$, we set 
\begin{equation}
R_{Z}(\cdot)=\sum_{i=1}^{N}\rho_{c_{i}}(\cdot-z_{i})~~\text{and}~~S_{Z}(\cdot)=\sum_{i=1}^{N}\varphi_{c_{i}}
(\cdot-z_{i}).
\label{4.12}
\end{equation}
For $0<\gamma<\gamma_{0}$, we define the function 
\begin{align*}
\mathcal{Y}:(-\gamma,\gamma)^{N}&\times B_{H^{2}}(R_{Z},\gamma)\rightarrow\mathbb{R}^{N},\\
&(y_{1},\ldots,y_{N},v)\mapsto\left(\mathcal{Y}^{1}(y_{1},\ldots,y_{N},v),\ldots,\mathcal{Y}^{N}(y_{1},\ldots,
y_{N},v)\right)
\end{align*}
with
$$\mathcal{Y}^{i}(y_{1},\ldots,y_{N},v)=\int_{\mathbb{R}}\left(v-\sum_{j=1}^{N}\rho_{c_{j}}(\cdot-z_{j}-y_{j})\right)\partial_{x}\rho_{c_{i}}(\cdot-z_{i}-y_{i}).$$
$Y$ is clearly of class $C^{1}$. For $i=1,\ldots,N$,
\begin{equation*}
\frac{\partial \mathcal{Y}^{i}}{\partial y_{i}}(y_{1},\ldots,y_{N},v)=-\int_{\mathbb{R}}\left(v-\sum\limits_{\underset{j \neq i}{1\le j\le N}}\rho_{c_{j}}(\cdot-z_{j}-y_{j})\right)\partial^{2}_{x}\rho_{c_{i}}(\cdot-z_{i}-y_{i})
\end{equation*}
and for $j\neq i$,
$$\frac{\partial\mathcal{Y}^{i}}{\partial y_{j}}(y_{1},\ldots,y_{N},v)=\int_{\mathbb{R}}\partial_{x}\rho_{c_{j}}(\cdot-z_{j}-y_{j})\partial_{x}\rho_{c_{i}}(\cdot-z_{i}-y_{i}).$$
Hence
\begin{equation*}
\frac{\partial\mathcal{Y}^{i}}{\partial y_{i}}(0,\ldots,0,R_{Z})=\|\partial_{x}\rho_{c_{i}}\|^{2}_{L^{2}(\mathbb{R})}=\frac{c^{2}_{i}}{54}\ge\frac{c^{2}_{1}}{54}
\end{equation*}
and for $j\neq i$, using the exponential decay of $\varphi_{c_{i}}$ and that $|z_{i}-z_{i-1}|>L/2$, for $L>L_{0}>0$ with $L_{0}\gg 1$, it holds 
\begin{align*}
&\left|\frac{\partial \mathcal{Y}^{i}}{\partial y_{j}}(0,\ldots,0,R_{Z})\right|\\
&\hspace{1cm}=\left| \int_{\mathbb{R}}\partial_{x}\rho_{c_{j},\alpha}(\cdot-z_{j})\partial_{x}\rho_{c_{i},\alpha}(\cdot-z_{i})\right|\\
&\hspace{1cm}=\left| \int_{\mathbb{R}}\rho_{c_{j},\alpha}(\cdot-z_{j})\partial^{2}_{x}\rho_{c_{i},\alpha}(\cdot-z_{i})\right|\\
&\hspace{1cm}\le\frac{1}{9}(c_{i}-\alpha)(c_{j}-\alpha)\left\lbrace\int_{\mathbb{R}}e^{-|x-z_{j}|-|x-z_{i}|}dx+2\int_{\mathbb{R}}e^{-|x-z_{j}|-2|x-z_{i}|}dx\right.\\
&\hspace{2cm}\left.+\frac{1}{2}\int_{\mathbb{R}}e^{-2|x-z_{j}|-|x-z_{i}|}dx
+\int_{\mathbb{R}}e^{-2|x-z_{j}|-2|x-z_{i}|}dx\right\rbrace\\
&\hspace{1cm}\le O(e^{-L/4}).
\end{align*}
We deduce that, for $L>0$ large enough, $D_{(y_{1},\ldots,y_{N})}\mathcal{Y}(0,\ldots,0,R_{Z})=D+P$ where $D$ is an invertible diagonal matrix with $\|D^{-1}\|\le (c_{1}/3\sqrt{6})^{-2}$ and $\|P\|\le O(e^{-L/4})$. Hence there exists $L_{0}>0$ such that for $L>L_{0}$, $D_{(y_{1},\ldots,y_{N})}\mathcal{Y}(0,\ldots,0,R_{Z})$ is invertible with an inverse matrix of norm smaller than $2(c_{1}/3\sqrt{6})^{-2}$. From the Implicit Function Theorem we deduce that there exists $\beta_{0}>0$ and $C^{1}$ functions $(y_{1},\ldots,y_{N})$ from $B_{H^{2}}(R_{Z},\beta_{0})$ to a neighborhood of $(0,\ldots,0)$ which are uniquely determined such that
$$\mathcal{Y}(y_{1}(v),\ldots,y_{N}(v),v)=0,~~\forall v\in B_{H^{2}}(R_{Z},\beta_{0}).$$
In particular, there exists $C_{0}>0$ such that if $v\in B_{H^{2}}(R_{Z},\beta)$, with $0<\beta\le\beta_{0}$, then
\begin{equation}
\sum_{i=1}^{N}|y_{i}(v)|\le C_{0}\beta.
\label{4.13}
\end{equation}
Note that $\beta_{0}$ and $C_{0}$ only depend on $c_{1}$ and $L_{0}$ and not on the point $(z_{1},\ldots,z_{N})$. For $v\in B_{H^{2}}(R_{Z},\beta_{0})$ we set $\tilde{x}_{i}(v)=z_{i}+y_{i}(v)$. Assuming that $\beta_{0}\le L_{0}/8C_{0}$, $(\tilde{x}_{1}(v),\ldots,\tilde{x}_{N}(v))$ are thus $C^{1}$ functions on $B_{H^{2}}(R_{Z},\beta)$ satisfying 
\begin{equation}
\tilde{x}_{i}(v)-\tilde{x}_{i-1}(v)=z_{i}-z_{i-1}+y_{i}(v)-y_{i-1}(v)>\frac{L}{2}-2C_{0}\beta>\frac{L}{4}.
\label{4.14}
\end{equation}
For $L>L_{0}$ and $0<\gamma<\gamma_{0}<\beta_{0}/2$ to be chosen later, we define the modulation of $v$ in the following way: we cover the trajectory of $v$ by a finite number of open balls in the following way:
$$\lbrace v(t), t\in[0,t_{0}]\rbrace\subset\bigcup_{k=1,\ldots,M}B_{H^{2}}(R_{Z^{k}},2\gamma).$$
This is possible thanks to Remark \ref{Remark 4.1}.
It is worth noticing that, since $0<\gamma<\gamma_{0}<\beta_{0}/2$, the functions $\tilde{x}_{i}(v)$ are uniquely determined for $v\in B_{H^{2}}(R_{Z^{k}},2\gamma)\cap B_{H^{2}}(R_{Z^{k'}},2\gamma)$. We can thus define the functions $t\mapsto\tilde{x}_{i}(t)$ on $[0,t_{0}]$ by setting $\tilde{x}_{i}(t)=\tilde{x}_{i}(v(t))$. By construction
\begin{equation}
\int_{\mathbb{R}}\left(v(t,\cdot)-\sum_{j=1}^{N}\rho_{c_{j}}(\cdot-\tilde{x}_{j}(t))\right)\partial_{x}
\rho_{c_{i}}(\cdot-\tilde{x}_{i}(t))=0.
\label{4.15}
\end{equation}
For $0<\gamma<\gamma_{0}$, with $\gamma_{0}\ll 1$, using that $u\in U(\gamma,L/2)$ and \eqref{4.13}, we have 
\begin{align*}
\|u(t)&-S_{\tilde{X}(t)}\|_{\mathcal{H}}\\
&\le\|u(t)-S_{Z}(t)\|_{\mathcal{H}}+\sum_{i=1}^{N}\|\varphi_{c_{i}}(\cdot-z_{i})-\varphi_{c_{i}}
(\cdot-z_{i}-y_{i}(v(t)))\|_{L^{2}(\mathbb{R})}\\
&\le\gamma+\sqrt{2}\sum_{i=1}^{N}\left(\int_{\mathbb{R}}\varphi^{2}_{c_{i}}(x)dx-\int_{\mathbb{R}}\varphi_{c_{i}}(x-z_{i})\varphi_{c_{i}}(x-z_{i}-y_{i}(v(t)))dx\right)^{1/2}\\
&=\gamma+\sqrt{2}\sum_{i=1}^{N}c_{i}\left(1-e^{-|y_{i}(v(t))|}-|y_{i}(v(t))|e^{-|y_{i}(v(t))|}\right)^{1/2}\\
&\le\gamma+\sum_{i=1}^{N}O(|y_{i}(v(t))|)\\
&\le O(\gamma),
\end{align*}
where we apply two time the mean value theorem with the function $\varphi$ on $[0,|y_{i}(v(t))|]$ for substituting $(1-e^{-|y_{i}(v(t))|})$ by $|y_{i}(v(t))|e^{-\theta}$, with $\theta\in]0,|y_{i}(v(t))|[$, and this proves \eqref{4.5} (see Fig. \ref{aa1}-\ref{aa2}).

The estimate \eqref{4.6} follows directly by using \eqref{4.5}, Remark \ref{Remark 4.1} and the Sobolev embedding of $H^{2}(\mathbb{R})$ into $C^{1}(\mathbb{R})$.

To prove that the speed of $\tilde{x}_{i}(\cdot)$ stays close to $c_{i}$, we set 
$$S_{j}(t)=\varphi_{c_{j}}(\cdot-\tilde{x}_{j}(t)),~~\varepsilon_{1}(t)=u(t)-\sum_{j=1}^{N}S_{j}(t)$$
and
$$R_{j}(t)=\rho_{c_{j}}(\cdot-\tilde{x}_{j}(t)),~~\varepsilon_{2}(t)=v(t)-\sum_{j=1}^{N}R_{j}(t).$$
One can notice that 
\begin{equation}
\partial^{2}_{x}R_{i}=4R_{i}-S_{i},
\label{4.16}
\end{equation}
and using the Fourier transformation 
\begin{align}
(1-\partial^{2}_{x})^{-1}(4-\partial^{2}_{x})^{-1}(\cdot)&=\mathcal{F}^{-1}\left[\frac{1}{(1+\omega^{2})(4+\omega^{2})}\right](\cdot)\nonumber\\
&=\mathcal{F}^{-1}\left[\frac{1}{3(1+\omega^{2})}-\frac{1}{3(4+\omega^{2})}\right](\cdot)\nonumber\\
&=\frac{1}{3}(1-\partial^{2}_{x})^{-1}(\cdot)-\frac{1}{3}(4-\partial^{2}_{x})^{-1}(\cdot).
\label{4.17}
\end{align}
Differentiating \eqref{4.15} with respect to time and using \eqref{4.16}, we get 
$$\int_{\mathbb{R}}\partial_{t}\varepsilon_{2}\partial_{x}R_{i}=\dot{\tilde{x}}_{i}(t)\left(4\int_{\mathbb{R}}\varepsilon_{2}R_{i}-\int_{\mathbb{R}}\varepsilon_{2}S_{i}\right)$$
and thus
\begin{align}
\left|\int_{\mathbb{R}}\partial_{t}\varepsilon_{2}\partial_{x}R_{i}\right|&\le|\dot{\tilde{x}}_{i}(t)|
\left(4\|\varepsilon_{2}\|_{L^{\infty}(\mathbb{R})}\|R_{i}\|_{L^{1}(\mathbb{R})}
+\|\varepsilon_{2}\|_{L^{\infty}(\mathbb{R})}\|S_{i}\|_{L^{1}(\mathbb{R})}\right)\nonumber\\
&\le|\dot{\tilde{x}}_{i}(t)-c_{i}|O(\gamma)+O(\gamma),
\label{4.18}
\end{align}
we point out that $\left\|S_{i}\right\|_{L^{1}(\mathbb{R})}=c$ and $\left\|R_{i}\right\|_{L^{1}(\mathbb{R})}=c/2$. 
Substituting $u$ by $\varepsilon_{1}+\sum_{j=1}^{N}S_{j}$ in \eqref{1.6} and using that $S_{j}$ satisfies 
$$\partial_{t}S_{j}=-(\dot{\tilde{x}}_{j}(t)-c_{j})\partial_{x}S_{j}-\frac{1}{2}\partial_{x}S_{j}^{2}-\frac{3}{2}(1-\partial^{2}_{x})^{-1}\partial_{x}S^{2}_{j},$$
we infer that $\varepsilon_{1}$ satisfies on $[0,t_{0}]$,
\begin{align*}
\partial_{t}\varepsilon_{1}&-\sum_{j=1}^{N}(\dot{\tilde{x}}_{j}(t)-c_{j})\partial_{x}S_{j}\\
&=-\frac{1}{2}\partial_{x}\left[\left(\varepsilon_{1}+
\sum_{j=1}^{N}S_{j}\right)^{2}-\sum_{j=1}^{N}S^{2}_{j}\right]\\
&\hspace{1cm}-\frac{3}{2}\partial_{x}(1-\partial^{2}_{x})^{-1}\left[\left(\varepsilon_{1}
+\sum_{j=1}^{N}S_{j}\right)^{2}-\sum_{j=1}^{N}S^{2}_{j}\right].
\end{align*}
Multiplying by $(4-\partial^{2}_{x})^{-1}(\cdot)$ and using \eqref{4.17}, we get 
$$\partial_{t}\varepsilon_{2}-\sum_{j=1}^{N}(\dot{\tilde{x}}_{j}(t)-c_{j})\partial_{x}R_{j}
=-\frac{1}{2}\partial_{x}(1-\partial^{2}_{x})^{-1}
\left[\left(\varepsilon_{1}+\sum_{j=1}^{N}S_{j}\right)^{2}-\sum_{j=1}^{N}S^{2}_{j}\right].$$
Taking the $L^{2}$ scalar product with $\partial_{x}R_{i}$, integrating by parts, we find 
\begin{align}
-(\dot{\tilde{x}}_{i}(t)&-c_{i})\int_{\mathbb{R}}(\partial_{x}R_{i})^{2}\nonumber\\
&=-\int_{\mathbb{R}}\partial_{t}\varepsilon_{2}\partial_{x}R_{i}+\sum\limits_{\underset{j \neq i}{1\le j\le N}}(\dot{\tilde{x}}_{j}(t)-c_{j})\int_{\mathbb{R}}(\partial_{x}R_{i})(\partial_{x}R_{j})\nonumber\\
&\hspace{1cm}+\frac{1}{2}\int_{\mathbb{R}}(1-\partial^{2}_{x})^{-1}
\left[\left(\varepsilon_{1}+\sum_{j=1}^{N}S_{j}\right)^{2}-\sum_{j=1}^{N}S^{2}_{j}\right]\partial^{2}_{x}
R_{i}.
\label{4.19}
\end{align}
We set
\begin{equation*}
Q=\left(\varepsilon_{1}+\sum_{j=1}^{N}S_{j}\right)^{2}-\sum_{j=1}^{N}S^{2}_{j}\\
=\varepsilon_{1}^{2}+2\varepsilon_{1}\left(\sum_{j=1}^{N}S_{j}\right)
+\sum\limits_{\underset{j \neq i}{1\le i, j\le N}}S_{i}S_{j},
\end{equation*}
then 
\begin{align*}
(1-\partial^{2}_{x})^{-1}Q&=(1-\partial^{2}_{x})^{-1}\varepsilon_{1}^{2}
+2\sum_{j=1}^{N}(1-\partial^{2}_{x})^{-1}\left(\varepsilon_{1}S_{j}\right)
+\sum\limits_{\underset{j \neq i}{1\le i, j\le N}}(1-\partial^{2}_{x})^{-1}\left(S_{i}S_{j}\right)\\
&=I+J+K.
\end{align*}
We have the following estimates 
\begin{equation*}
I=\frac{1}{2}\int_{\mathbb{R}}e^{-|x-x'|}\varepsilon^{2}_{1}(x')
\le\frac{1}{2}\|e^{-|\cdot|}\|_{L^{\infty}(\mathbb{R})}\|\varepsilon_{1}\|^{2}_{L^{2}(\mathbb{R})}
=\frac{1}{2}\|\varepsilon_{1}\|^{2}_{L^{2}(\mathbb{R})},
\end{equation*}
\begin{equation*}
J=\sum_{j=1}^{N}\int_{\mathbb{R}}e^{-|x-x'|}\varepsilon_{1}(x')S_{j}(x')
\le\|e^{-|\cdot|}\|_{L^{\infty}(\mathbb{R})}\|\varepsilon_{1}\|_{L^{2}(\mathbb{R})}
\sum_{j=1}^{N}\|S_{j}\|_{L^{2}(\mathbb{R})}
=\left(\sum_{j=1}^{N}c_{j}\right)\|\varepsilon_{1}\|_{L^{2}(\mathbb{R})}
\end{equation*}
and 
\begin{equation*}
K=\frac{1}{2}\sum\limits_{\underset{j \neq i}{1\le i, j\le N}}\int_{\mathbb{R}}
e^{-|x-x'|}S_{j}(x')S_{i}(x')
\le\frac{1}{2}\sum\limits_{\underset{j \neq i}{1\le i, j\le N}}
\int_{\mathbb{R}}S_{j}(x')S_{i}(x').
\end{equation*}
Thus, using \eqref{4.5} and the exponential decay of $S_{j}$, it holds
\begin{equation*}
\|(1-\partial^{2}_{x})^{-1}Q\|_{L^{\infty}(\mathbb{R})}\le O(\gamma)+O(e^{-L/4})
\end{equation*}
and then
\begin{align}\label{papa1}
\left|\frac{1}{2}\int_{\mathbb{R}}[(1-\partial^{2}_{x})^{-1}Q]\partial^{2}_{x}R_{i}\right|&
\le\frac{1}{2}\|(1-\partial^{2}_{x})^{-1}Q\|_{L^{\infty}(\mathbb{R})}\|\partial^{2}_{x}R_{i}\|_{L^{1}(\mathbb{R})}\nonumber\\
&\le O(\gamma)+O(e^{-L/4}),
\end{align}
where $\left\|\partial^{2}_{x}R_{i}\right\|_{L^{1}(\mathbb{R})}=c/3$.
Now, combining \eqref{4.18}, \eqref{papa1}, and using the exponential decay of $\partial_{x}R_{i}$, it holds
$$\left|\dot{\tilde{x}}_{i}(t)-c_{i}\right|\|\partial_{x}R_{i}\|^{2}_{L^{2}(\mathbb{R})}\le\left|\dot{\tilde{x}}_{i}(t)-c_{i}\right|O(\gamma)+O(\gamma)+O(e^{-L/4}),$$
then 
$$\left|\dot{\tilde{x}}_{i}(t)-c_{i}\right|\left(\frac{c_i^2}{54}-O(\gamma)\right)
\le O(\gamma)+O(e^{-L/4}),$$
which yields \eqref{4.7}.

Taking $0<\gamma<\gamma_{0}$ and $L>L_{0}>0$ with $\gamma_{0}\ll 1$ and $L_{0}\gg 1$, combining \eqref{1.11}-\eqref{1.13}, \eqref{4.7} and \eqref{4.14}, we deduce that 
\begin{align*}
\tilde{x}_{i}(t)-\tilde{x}_{i-1}(t)&=\tilde{x}_{i}(0)-\tilde{x}_{i-1}(0)+(c_{i}-c_{i-1})t\\
&\ge L-2C_{0}\gamma_{0}+\frac{(c_{i}-c_{i-1})t}{2}\\
&\ge \frac{3L}{4}+\frac{(c_{i}-c_{i-1})t}{2},
\end{align*}
this proves \eqref{4.8}.

From \eqref{4.6}, we infer that 
$$v(x)=\sum_{j=1}^{N}\rho_{c_{j}}(x-\tilde{x}_{j})+O(\gamma),~~\forall x\in\mathbb{R},$$
please note that we abuse notation by writing $\varepsilon_{2}(x)=O(\gamma)$. 
Applying this formula with $x=\xi^{i}_{1}$ and $v(\xi^{i}_{1})=\max_{x\in J_{i}}v(x)$, and using \eqref{4.8}, it holds
\begin{align*}
v(\xi^{i}_{1})&=\max_{x\in J_{i}}\left\lbrace\sum_{j=1}^{N}\rho_{c_{j}}(x-\tilde{x}_{j})\right\rbrace+O(\gamma)\\
&=\frac{c_{i}}{6}+O(e^{-L/4})+O(\gamma)\\
&\ge \frac{c_{i}}{7}.
\end{align*}
On the other hand, for $x\in J_{i}\setminus[\tilde{x}_{i}(t)-L/12,~\tilde{x}_{i}(t)+L/12]$, we get
$$v(x)\le\frac{c_{i}}{3}e^{-L/12}+O(e^{-L/4})+O(\gamma)\\
\le \frac{c_{i}}{8}.$$
This ensures that $\xi^{i}_{1}\in[\tilde{x}_{i}(t)-L/12,~\tilde{x}_{i}(t)+L/12]$, and this concluded the proof of the lemma.
\hfill $ \square $ \vspace*{2mm}

\subsection{Monotonicity property}\label{Subsection 4.2}
\begin{center}
\begin{figure}[ht]
\centering
\includegraphics[width=10cm, height=6cm]{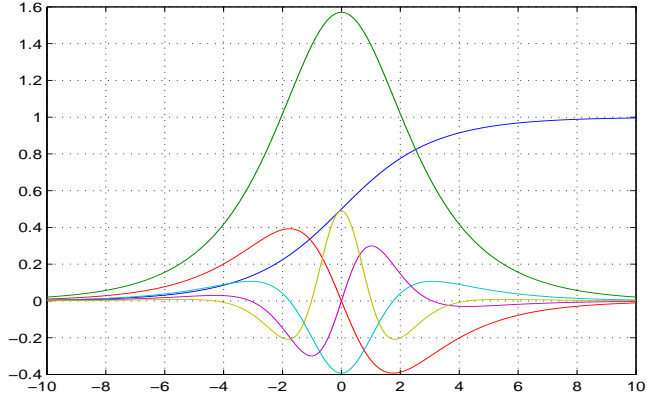}
\caption{$\psi^{(q)}(x)$, $q=0,1,2,3,4,5$, profiles. Note that $\psi^{(5)}$ will not be used. \label{fig4}}
\end{figure}
\end{center}

Thanks to the preceding lemma, for $\varepsilon_{0}>0$ small enough and $L_{0}>0$ large enough, one can construct $N~C^{1}$ functions $\tilde{x}_{1},\ldots,\tilde{x}_{N}$ defined on $[0,t_{0}]$ such that \eqref{4.5}-\eqref{4.9} are satisfied. In this subsection, we state the almost monotonicity of functionals that are very close to the energy at the right of $i$th bump, $i=1,\ldots,N-1$ of $u$ (respectively of $v$). Let $\psi$ be a 
$C^{\infty}$ test-function (see Fig. \ref{fig4}) such that 
\begin{equation}
  \left\{
    \begin{aligned}
     &0<\psi(x)<1,~\psi'(x)>0,&x&\in\mathbb{R},\\
     &|\psi^{(q)}(x)|\le 10\psi'(x),~q=2,3,4,5, &x&\in[-10,10],\\
    \end{aligned}
  \right.
  \label{4.24}
\end{equation}
and
\begin{equation}
  \psi(x)=\left\{
    \begin{aligned}
     &e^{-|x|},&x&<-10,\\
     &1-e^{-|x|},&x&>10.\\
    \end{aligned}
  \right.
  \label{4.25}
\end{equation}

Setting $\psi_{K}=\psi(\cdot/K)$, we introduce for $i=2,\ldots,N$,
\begin{equation}
\mathcal{J}_{i,K}(t)=\int_{\mathbb{R}}\left(4v^{2}+5v^{2}_{x}+v^{2}_{xx}\right)\psi_{i,K}(t),
\label{4.26}
\end{equation}
where $\psi_{i,K}(t,x)=\psi_{K}(x-y_{i}(t))$ with $y_{i}$'s as in \eqref{4.9}. Note that $\mathcal{J}_{i,K}(t)$ is close to $\|u(t)\|^{2}_{\mathcal{H}(x>y_{i}(t))}$ (respectively to $\|v(t)\|^{2}_{H^{2}(x>y_{i}(t))}$) and thus measures the energy at the right of the $(i-1)$th bump of $u$ (respectively of $v$). Finally, we set 
\begin{equation}
\sigma_{0}=\frac{1}{4}\min\lbrace c_{1},c_{2}-c_{1},\ldots,c_{N}-c_{N-1}\rbrace.
\label{4.27}
\end{equation}

We have the following monotonicity result.

\begin{Pro}[Exponential Decay of the Functional $\mathcal{J}_{i,K}(t)$]\label{Proposition 4.1}
Let $u\in \mathcal{X}([0,T[)$, with $0<T\le +\infty$, be a solution of equation \eqref{1.1} that satisfies \eqref{1.11}-\eqref{1.13} and \eqref{4.5}-\eqref{4.6}. There exist $\gamma_{0}>0$ and $L_{0}>0$ only depending on $c_{1}$ such that if $0<\gamma<\gamma_{0}$ and $L>L_{0}>0$, then for any $4\le K\lesssim \sqrt{L}$,
\begin{equation}
\mathcal{J}_{i,K}(t)-\mathcal{J}_{i,K}(0)\le O(e^{-\frac{L}{8K}}),~~\forall t\in[0,t_{0}],~~i=2,\ldots,N.
\label{4.28}
\end{equation}
\end{Pro}

The proof of Proposition \ref{Proposition 4.1} relies on the following Virial type identity.

\begin{Lem}[Virial Type Identity]\label{Lemma 4.2}
Let $u\in \mathcal{X}([0,T[)$, with $0<T\le +\infty$, be a solution of equation \eqref{1.1} that satisfies \eqref{1.11}-\eqref{1.13}. For any smooth space function $g :\mathbb{R}\mapsto\mathbb{R}$, it holds
\begin{align}
\frac{d}{dt}\int_{\mathbb{R}}&\left(4v^{2}+5v^{2}_{x}+v^{2}_{xx}\right)g\nonumber\\
&=\frac{2}{3}\int_{\mathbb{R}}u^{3}g'-4\int_{\mathbb{R}}u^{2}vg'
-\frac{1}{2}\int_{\mathbb{R}}u^{2}vg'''
+\frac{1}{2}\int_{\mathbb{R}}u^{2}v_{x}g''
+\int_{\mathbb{R}}uhg'\nonumber\\
&\hspace{1cm}+\frac{1}{2}\int_{\mathbb{R}}uh_{x}g''
-\frac{5}{2}\int_{\mathbb{R}}vh_{x}g''-2\int_{\mathbb{R}}v_{x}hg''+\frac{1}{2}\int_{\mathbb{R}}vh_{x}g^{(4)}
\label{4.29}
\end{align}
where $y=(1-\partial^{2}_{x})u$, $v=(4-\partial^{2}_{x})^{-1}u$ and $h=(1-\partial^{2}_{x})^{-1}u^{2}$.
\end{Lem}
The full proof of Lemma \ref{Lemma 4.2} is given in the Appendix \ref{Appendix}.

\textbf{Proof of Proposition 4.1.}
We first note that, combining \eqref{4.7} and \eqref{4.9}, it holds for $i=2,\ldots,N$,
\begin{align}
\dot{y}_{i}(t)&=\frac{\dot{\tilde{x}}_{i-1}(t)+\dot{\tilde{x}}_{i}(t)}{2}\nonumber\\
&=\frac{c_{i-1}+c_{i}}{2}+O(\gamma)\nonumber\\
&\ge c_{i-1}+O(\gamma)\nonumber\\
&\ge\frac{c_{1}}{2}.
\label{4.49}
\end{align}
Recall that the assumption \eqref{1.11} ensures that $u\ge 0$ and $v\ge 0$ on $\mathbb{R}$. Now,  applying the Virial type identity \eqref{4.29} with $g=\psi_{i,K}$ and using \eqref{4.49}, we get
\begin{align}
\frac{d}{dt}\mathcal{J}_{i,K}(t)
&=-\dot{y}_{i}\int_{\mathbb{R}}\left(4v^{2}+5v^{2}_{x}+v^{2}_{xx}\right)\psi'_{i,K}
+\frac{2}{3}\int_{\mathbb{R}}u^{3}\psi_{i,K}'-4\int_{\mathbb{R}}u^{2}v\psi_{i,K}'\nonumber\\
&\hspace{1cm}-\frac{1}{2}\int_{\mathbb{R}}u^{2}v\psi_{i,K}'''
+\frac{1}{2}\int_{\mathbb{R}}u^{2}v_{x}\psi''_{i,K}
+\int_{\mathbb{R}}uh\psi_{i,K}'+\frac{1}{2}\int_{\mathbb{R}}uh_{x}\psi_{i,K}''\nonumber\\
&\hspace{2cm}-\frac{5}{2}\int_{\mathbb{R}}vh_{x}\psi''_{i,K}
-2\int_{\mathbb{R}}v_{x}h\psi''_{i,K}+\frac{1}{2}\int_{\mathbb{R}}vh_{x}\psi^{(4)}_{i,K}\nonumber\\
& \le -\dot{y}_{i}\int_{\mathbb{R}}\left(4v^{2}+5v^{2}_{x}+v^{2}_{xx}\right)\psi'_{i,K}
+\frac{2}{3}\int_{\mathbb{R}}u^{3}\psi_{i,K}'-\frac{1}{2}\int_{\mathbb{R}}u^{2}v\psi'''_{i,K}\nonumber\\
&\hspace{1cm}+\frac{1}{2}\int_{\mathbb{R}}u^{2}v_{x}\psi_{i,K}''
+\int_{\mathbb{R}}uh\psi_{i,K}'+\frac{1}{2}\int_{\mathbb{R}}uh_{x}\psi_{i,K}''-\frac{5}{2}\int_{\mathbb{R}}vh_{x}\psi''_{i,K}\nonumber\\
&\hspace{2cm}-2\int_{\mathbb{R}}v_{x}h\psi''_{i,K}
+\frac{1}{2}\int_{\mathbb{R}}vh_{x}\psi^{(4)}_{i,K}\nonumber\\
&\le -\frac{c_{1}}{2}\int_{\mathbb{R}}\left(4v^{2}+5v^{2}_{x}+v^{2}_{xx}\right)\psi'_{i,K}+\sum_{k=1}^{8}J_{k}.
\label{4.50}
\end{align}
We claim that for $k=1,\ldots,8$, it holds
\begin{equation}
J_{k}\le\frac{c_{1}}{20}\int_{\mathbb{R}}\left(4v^{2}+5v^{2}_{x}+v^{2}_{xx}\right)\psi'_{i,K}+\frac{C}{K}\|u_{0}\|^{3}_{\mathcal{H}}e^{-\frac{1}{K}(\sigma_{0}t+L/8)}.
\label{4.51}
\end{equation}
We divide $\mathbb{R}$ into two regions $D_{i}$ and $D^{c}_{i}$ with
$$D_{i}=\left[\tilde{x}_{i-1}(t)+\frac{L}{4},\tilde{x}_{i}(t)-\frac{L}{4}\right],~~i=2,\ldots,N.$$
Combining \eqref{4.8} and \eqref{4.9}, one can check that for $x\in D^{c}_{i}$,
\begin{align}
|x-y_{i}(t)|&\ge\frac{\tilde{x}_{i}(t)-\tilde{x}_{i-1}(t)}{2}-\frac{L}{4}\nonumber\\
&\ge\frac{c_{i}-c_{i-1}}{4}t+\frac{L}{8}\nonumber\\
&\ge\sigma_{0}t+\frac{L}{8}.
\label{4.52}
\end{align}
Let us begin by an estimate of $J_{1}$. Using \eqref{2.3}, \eqref{4.52} and the exponential decay of $\psi'_{i,K}$ on $D^{c}_{i}$, we get 
\begin{align}
\frac{2}{3}\int_{\mathbb{R}}u^{3}\psi'_{i,K}&=\frac{2}{3}\int_{D_{i}}u^{3}\psi_{i,K}'+\frac{2}{3}\int_{D^{c}_{i}}u^{3}\psi'_{i,K}\nonumber\\
&\le\frac{2}{3}\|u\|_{L^{\infty}(D_{i})}\int_{\mathbb{R}}
u^{2}\psi'_{i,K}+\frac{2}{3}\|\psi'_{i,K}\|_{L^{\infty}(D^{c}_{i})}\|u\|_{L^{\infty}(\mathbb{R})}\|u\|^{2}_{L^{2}(\mathbb{R})}\nonumber\\
&\le\frac{2}{3}\|u\|_{L^{\infty}(D_{i})}\int_{\mathbb{R}}
u^{2}\psi'_{i,K}+
\frac{C}{K}\|u\|^{3}_{L^{2}(\mathbb{R})}e^{-\frac{1}{K}(\sigma_{0}t+L/8)}.
\label{4.53}
\end{align}
Note that, using the exponential decay of $|\psi'''_{j,K}|$ on $D^{c}_{i}$, 
$|\psi'''_{i,K}|\le(10/K^{2})\psi'_{i,K}$ on $D_{i}$, with $K\ge 4$, and that $\|u\|_{\mathcal{H}}=\|u_{0}\|_{\mathcal{H}}$,  we have
\begin{align}
\int_{\mathbb{R}}u^{2}\psi'_{i,K}&=\int_{\mathbb{R}}(4v-v_{xx})^{2}\psi'_{i,K}\nonumber\\
&=16\int_{\mathbb{R}}v^{2}\psi'_{i,K}+\int_{\mathbb{R}}v^{2}_{xx}\psi'_{i,K}-8\int_{\mathbb{R}}vv_{xx}\psi'_{i,K}\nonumber\\
&=16\int_{\mathbb{R}}v^{2}\psi'_{i,K}+\int_{\mathbb{R}}v^{2}_{xx}\psi'_{i,K}+8\int_{\mathbb{R}}
v^{2}_{x}\psi'_{i,K}+4\int_{\mathbb{R}}\partial_{x}(v^{2})\psi''_{i,K}\nonumber\\
&=16\int_{\mathbb{R}}v^{2}\psi'_{i,K}+\int_{\mathbb{R}}v^{2}_{xx}\psi'_{i,K}+8\int_{\mathbb{R}}
v^{2}_{x}\psi'_{i,K}-4\int_{\mathbb{R}}v^{2}\psi'''_{i,K}\nonumber\\
&\le\int_{\mathbb{R}}\left(16v^{2}+8v^{2}_{x}+v^{2}_{xx}\right)\psi'_{i,K}
+4\int_{D_{i}}v^{2}|\psi'''_{j,K}|+4\int_{D^{c}_{i}}v^{2}|\psi'''_{j,K}|\nonumber\\
&\le\int_{\mathbb{R}}\left(16v^{2}+8v^{2}_{x}+v^{2}_{xx}\right)\psi'_{i,K}
+\frac{40}{K^{2}}\int_{\mathbb{R}}v^{2}\psi'_{i,K}
+\frac{C}{K^{3}}\|u_{0}\|^{2}_{\mathcal{H}}e^{-\frac{1}{K}(\sigma_{0}t+L/8)}\nonumber\\
&\le 5\int_{\mathbb{R}}\left(4v^{2}+5v^{2}_{x}+v^{2}_{xx}\right)\psi'_{i,K}
+\frac{C}{K^{3}}\|u_{0}\|^{2}_{\mathcal{H}}e^{-\frac{1}{K}(\sigma_{0}t+L/8)}.
\label{4.54}
\end{align}
Now, using the exponential decay of $\varphi_{c_{i}}$ on $D_{i}$, \eqref{4.5}, and proceeding as for the estimate \eqref{3.10} (see Lemma \ref{Lemma 3.2}), it holds
\begin{align}
\|u\|_{L^{\infty}(D_{i})} &\le\left\|u-\sum_{j=1}^{N}\varphi_{c_{j}}(\cdot-\tilde{x}_{j}(t))\right\|_{L^{\infty}(D_{i})}+\sum_{j=1}^{N}\left\|\varphi_{c_{j}}(\cdot-\tilde{x}_{j}(t))\right\|_{L^{\infty}(D_{i})}\nonumber\\
&\le O(\gamma^{1/2})+O(e^{-L/8}).
\label{4.55}
\end{align}
Therefore, for $0<\gamma<\gamma_{0}$ and $L>L_{0}>0$, with $\gamma_{0}\ll 1$ and $L_{0}\gg 1$, combining \eqref{4.53}-\eqref{4.54}, we obtain
$$J_{1}\le\frac{c_{1}}{20}\int_{\mathbb{R}}\left(4v^{2}+5v^{2}_{x}+v^{2}_{xx}\right)\psi'_{i,K}+\frac{C}{K}\|u\|^{3}_{L^{2}(\mathbb{R})}e^{-\frac{1}{K}(\sigma_{0}t+L/8)}.$$

Next, the estimate of $J_{2}$ on $D^{c}_{i}$ gives us
$$\left|\frac{1}{2}\int_{D^{c}_{i}}u^{2}v\psi'''_{i,K}\right| \le\frac{1}{2}\|\psi'''_{i,K}\|_{L^{\infty}(D^{c}_{i})}\|v\|_{L^{\infty}(\mathbb{R})}\|u\|^{2}_{L^{2}(\mathbb{R})}.$$
Note that, applying the H\"older inequality, we have for all $x\in\mathbb{R}$,
\begin{align}
v(x)&=\frac{1}{4}\int_{\mathbb{R}}e^{-2|x-x'|}u(x')dx'\nonumber\\
&\le\frac{1}{4}\left(\int_{\mathbb{R}}e^{-4|x-x'|}dx'\right)^{1/2}\left(\int_{\mathbb{R}}|u(x')|^{2}dx'\right)^{1/2}\nonumber\\
&=\frac{1}{4\sqrt{2}}\|u\|_{L^{2}(\mathbb{R})}
\label{4.56}
\end{align}
and thus, using \eqref{4.56} and the exponential decay of $|\psi'''_{i,K}|$ on $D^{c}_{i}$, it holds
\begin{equation}
\left|\frac{1}{2}\int_{D^{c}_{i}}u^{2}v\psi'''_{i,K}\right| \le\frac{1}{8\sqrt{2}K^{3}}\|u\|^{3}_{L^{2}(\mathbb{R})}e^{-\frac{1}{K}(\sigma_{0}t+L/8)}.
\label{4.57}
\end{equation}
Using that $|\psi'''_{i,K}|\le(10/K^{2})\psi'_{i,K}$ on $D_{i}$, the estimate of $J_{2}$ on $D_{i}$ leads to
\begin{equation}
\left|\frac{1}{2}\int_{D_{i}}u^{2}v\psi'''_{i,K}\right|\le\frac{5}{K^{2}}\|u\|_{L^{\infty}(D_{i})}\int_{\mathbb{R}}uv\psi'_{i,K}.
\label{4.58}
\end{equation}
Also, one can notice that, using the exponential decay of $|\psi'''_{j,K}|$ on $D^{c}_{i}$, 
$|\psi'''_{i,K}|\le(10/K^{2})\psi'_{i,K}$ on $\mathbb{R}$, with $K\ge 4$, and that $\|u\|_{\mathcal{H}}=\|u_{0}\|_{\mathcal{H}}$,  we have
\begin{align}
\int_{\mathbb{R}}uv\psi'_{i,K} &=\int_{\mathbb{R}}(4v-v_{xx})v\psi'_{i,K}\nonumber\\
&=4\int_{\mathbb{R}}v^{2}\psi'_{i,K}+\int_{\mathbb{R}}\partial_{x}(v\psi'_{i,K})v_{x}\nonumber\\
&=4\int_{\mathbb{R}}v^{2}\psi'_{j,K}+\int_{\mathbb{R}}v^{2}_{x}\psi'_{i,K}
+\int_{\mathbb{R}}vv_{x}\psi''_{i,K}\nonumber\\
&=\int_{\mathbb{R}}\left(4v^{2}+v^{2}_{x}\right)\psi'_{i,K}-\frac{1}{2}\int_{\mathbb{R}}v^{2}\psi'''_{i,K}\nonumber\\
&\le \int_{\mathbb{R}}\left(4v^{2}+v^{2}_{x}\right)\psi'_{i,K}+\frac{1}{2}\int_{D_{i}}v^{2}|\psi'''_{i,K}|+\frac{1}{2}\int_{D^{c}_{i}}v^{2}|\psi'''_{i,K}|\nonumber\\
&\le\int_{\mathbb{R}}\left(4v^{2}+v^{2}_{x}\right)\psi'_{i,K}
+\frac{5}{K^{2}}\int_{\mathbb{R}}v^{2}\psi'_{i,K}
+\frac{C}{K^{3}}\|u_{0}\|^{2}_{\mathcal{H}}e^{-\frac{1}{K}(\sigma_{0}t+L/8)}\nonumber\\
&\le 2\int_{\mathbb{R}}\left(4v^{2}+5v^{2}_{x}+v^{2}_{xx}\right)\psi'_{i,K}
+\frac{C}{K^{3}}\|u_{0}\|^{2}_{\mathcal{H}}e^{-\frac{1}{K}(\sigma_{0}t+L/8)}.
\label{4.59}
\end{align}
Therefore, for $0<\gamma<\gamma_{0}$ and $L>L_{0}>0$, with $\gamma_{0}\ll 1$ and $L_{0}\gg 1$, combining \eqref{4.55}, \eqref{4.57}-\eqref{4.59}, it holds
$$J_{2}\le\frac{c_{1}}{20}\int_{\mathbb{R}}\left(4v^{2}+5v^{2}_{x}+v^{2}_{xx}\right)\psi'_{i,K}+\frac{C}{K^{3}}\|u\|^{3}_{L^{2}(\mathbb{R})}e^{-\frac{1}{K}(\sigma_{0}t+L/8)}.$$
In the same way, using that $|v_{x}|\le 2v$ on $\mathbb{R}$ (see \eqref{3.36}), and the definition of 
$\psi_{i,K}$ (see \eqref{4.24} and \eqref{4.25}), we deduce the estimate of $J_{3}$.

Let us tackle now the estimate of $J_{4}$.  On $D^{c}_{i}$ we have 
\begin{align*}
\int_{D^{c}_{i}}uh\psi_{i,K}'&\le\|\psi'_{i,K}\|_{L^{\infty}(D^{c}_{i})}\int_{\mathbb{R}}uh\\
&=\|\psi'_{i,K}\|_{L^{\infty}(D^{c}_{i})}\int_{\mathbb{R}}u[(1-\partial^{2}_{x})^{-1}u^{2}]\\
&=\|\psi'_{i,K}\|_{L^{\infty}(D^{c}_{i})}\int_{\mathbb{R}}u^{2}[(1-\partial^{2}_{x})^{-1}u]\\
&\le \|\psi'_{i,K}\|_{L^{\infty}(D^{c}_{j})}\|(1-\partial^{2}_{x})^{-1}u\|_{L^{\infty}(\mathbb{R})}\|u\|^{2}_{L^{2}(\mathbb{R})}.
\end{align*}
Remark that, applying the H\"older inequality, we have for all $x\in\mathbb{R}$,
\begin{align}
(1-\partial^{2}_{x})^{-1}u(x)&\le\frac{1}{2}\int_{\mathbb{R}}e^{-|x-x'|}u(x')dx'\nonumber\\
&\le\frac{1}{2}\left(\int_{\mathbb{R}}e^{-2|x-x'|}dx'\right)^{1/2}\left(\int_{\mathbb{R}}|u(x')|^{2}dx'\right)^{1/2}\nonumber\\
&=\frac{1}{2}\|u\|_{L^{2}(\mathbb{R})}
\label{4.60}
\end{align}
and thus, using \eqref{4.60} and the exponential decay of $\psi'_{i,K}$ on $D^{c}_{i}$, it holds
\begin{equation}
\int_{D^{c}_{i}}uh\psi'_{i,K}\le\frac{1}{2K}\|u\|^{3}_{L^{2}(\mathbb{R})}e^{-\frac{1}{K}(\sigma_{0}t+L/8)}.
\label{4.61}
\end{equation}
The estimate of $J_{4}$ on $D_{i}$ leads to
\begin{align}
\int_{D_{i}}uh\psi'_{i,K} &\le\|u\|_{L^{\infty}(D_{i})}\int_{\mathbb{R}}\psi'_{i,K}[(1-\partial^{2}_{x})^{-1}u^{2}]\nonumber\\
&=\|u\|_{L^{\infty}(D_{i})}\int_{\mathbb{R}}u^{2}[(1-\partial^{2}_{x})^{-1}\psi'_{i,K}].
\label{4.62}
\end{align}
On the other hand, using that $|\psi'''_{i,K}|\le(10/K^{2})\psi'_{i,K}$ on $ \mathbb{R}$, we have
\begin{equation*}
(1-\partial^{2}_{x})\psi'_{i,K}(x)=\psi'_{i,K}(x)-\psi'''_{i,K}(x)
\ge\left(1-\frac{10}{K^{2}}\right)\psi'_{i,K}(x),~~\forall x\in\mathbb{R},
\end{equation*}
and since $K\ge 4$, it holds
\begin{equation}
(1-\partial^{2}_{x})^{-1}\psi'_{i,K}(x)\le\left(1-\frac{10}{K^{2}}\right)^{-1}\psi'_{i,K}(x),
~~\forall x\in\mathbb{R}.
\label{4.63}
\end{equation}
Therefore, for $0<\gamma<\gamma_{0}$ and $L>L_{0}>0$, with $\gamma_{0}\ll 1$ and $L_{0}\gg 1$, combining \eqref{4.55}, \eqref{4.61}-\eqref{4.63}, it holds
$$J_{4}\le \frac{c_{1}}{20}\int_{\mathbb{R}}\left(4v^{2}+5v^{2}_{x}+v^{2}_{xx}\right)\psi'_{i,K}+\frac{C}{K}\|u\|^{3}_{L^{2}(\mathbb{R})}e^{-\frac{1}{K}(\sigma_{0}t+L/8)}.$$

Noticing that for all $x\in\mathbb{R}$,
$$h(x)=\frac{e^{-x}}{2}\int_{-\infty}^{x}e^{x'}u^{2}(x')dx'+\frac{e^{x}}{2}\int_{x}^{+\infty}e^{-x'}u^{2}(x')dx'$$
and 
$$h_{x}(x)=-\frac{e^{-x}}{2}\int_{-\infty}^{x}e^{x'}u^{2}(x')dx'+\frac{e^{x}}{2}\int_{x}^{+\infty}e^{-x'}u^{2}(x')dx',$$
we infer that 
\begin{equation}
|h_{x}(x)|\le h(x),~~\forall x\in\mathbb{R}.
\label{4.64}
\end{equation}
Then, combining \eqref{4.24}, \eqref{4.64}, and proceeding as for the estimate of $J_{4}$, we deduce the estimate of $J_{5}$.

Now, combining \eqref{4.56} and \eqref{4.60}, we have for all $x\in\mathbb{R}$,
\begin{align}
(1-\partial^{2}_{x})^{-1}v(x)&=\frac{1}{3}(1-\partial^{2}_{x})^{-1}u(x)-\frac{1}{3}v(x)\nonumber\\
&\le\frac{1}{3}\left\|(1-\partial^{2}_{x})^{-1}u\right\|_{L^{\infty}(\mathbb{R})}+\frac{1}{3}\|v\|_{L^{\infty}(\mathbb{R})}\nonumber\\
&\le\frac{4+\sqrt{2}}{24}\|u\|_{L^{2}(\mathbb{R})},
\label{4.65}
\end{align}
and using the exponential decay of $\rho_{c_{i}}$ on $D_{i}$ and \eqref{4.6}, it holds
\begin{align}
\|v\|_{L^{\infty}(D_{i})} &\le\left\|v-\sum_{j=1}^{N}\rho_{c_{j}}(\cdot
-\tilde{x}_{j}(t))\right\|_{L^{\infty}(D_{i})}+\sum_{j=1}^{N}\left\|\rho_{c_{j}}(\cdot-\tilde{x}_{j}(t))\right\|_{L^{\infty}(D_{i})}\nonumber\\
&\le O(\gamma)+O(e^{-L/8}).
\label{4.66}
\end{align}
Therefore, combining \eqref{3.36}, \eqref{4.24}, \eqref{4.64}-\eqref{4.66}, and proceeding as for the estimate of $J_{4}$, we deduce the estimates of the remaining terms.

Finally, combining \eqref{4.50}, \eqref{4.51} and using that $\|u\|_{L^{2}(\mathbb{R})}\sim\|u_{0}\|_{\mathcal{H}}$, it holds actually
$$\frac{d}{dt}\mathcal{J}_{i,K}(t)\le \frac{C}{K}\|u_{0}\|^{3}_{\mathcal{H}}e^{-\frac{1}{K}(\sigma_{0}t+L/8)}.$$
Integrating between $0$ and $t$, we obtain
\begin{align*}
\mathcal{J}_{i,K}(t)-\mathcal{J}_{i,K}(0) &\le \frac{C}{K}\|u_{0}\|^{3}_{\mathcal{H}}\left(-\frac{K}{\sigma_{0}}e^{-\frac{1}{K}(\sigma_{0}t+L/8)}+\frac{K}{\sigma_{0}}e^{-\frac{L}{8K}}\right)\\
&\le \frac{C}{\sigma_{0}}\|u_{0}\|^{3}_{\mathcal{H}}e^{-\frac{L}{8K}},
\end{align*}
and this proves the proposition for smooth initial solutions.

For $ u\in\mathcal{X}([0,T[) $, we will use that for any $T_0>0 $ and any sequence 
$(u_{0,n})_{n\ge 1}\subset L^{2}(\mathbb{R})$ such that $(u_{0,n}-\partial^{2}_{x}u_{0,n})_{n\ge 1}\subset\mathcal{M}^{+}(\mathbb{R})$ and $u_{0,n}\rightarrow u_{0}$ in $L^{2}(\mathbb{R})$,  the sequence of emanating global weak solutions $(u_n)_{n\ge 1} $ to the DP  equation  satisfies
\begin{equation}
u_n\underset{n\to+\infty}{\longrightarrow} u~~\text{in}~~C\left([0,T_0];L^{2}(\mathbb{R})\right),
\label{R1}
\end{equation}
where $ u$ is the global weak solution emanating from $ u_0 $.
This fact can be easily deduced from the proof of the existence of the global weak solutions in \cite{MR2271927}. Indeed, by the same arguments developed in this proof, we obtain that, up to a subsequence,  $(u_{n})_{n\ge 1}$ converges in $ C\left([0,T_0];L^{2}(\mathbb{R})\right)$ towards  a solution of the DP equation emanating  from $ u_0$. 
 \eqref{R1} then follows by the uniqueness result. Combining \eqref{R1} and Remark \ref{Remark 4.1}, it follows that 
\begin{equation}
v_n\underset{n\to+\infty}{\longrightarrow} v~~\text{in}~~C\left([0,T_0];H^{2}(\mathbb{R})\right),
\label{R2}
\end{equation}
where $v=(4-\partial^{2}_{x})^{-1}u$.
For all $t\in[0,T[$, we set 
\begin{equation}
\mathcal{J}^{n}_{i,K}(t)=\mathcal{J}_{i,K}(u_{n}(t))=\int_{\mathbb{R}}\left(4v^{2}_{n}+5v^{2}_{n,x}+v^{2}_{n,xx}\right)\psi_{i,K}(t),
\label{R3}
\end{equation}
and we claim that 
\begin{equation}
\lim_{n\to +\infty}\sup_{0\le t<T}\left|\mathcal{J}^{n}_{i,K}(t)-\mathcal{J}_{i,K}(t)\right|=0.
\label{R4}
\end{equation}
Let $t\in[0,T[$ be fixed, we compute 
\begin{align}
\mathcal{J}^{n}_{i,K}(t)-\mathcal{J}_{i,K}(t)&=4\int_{\mathbb{R}}(v^{2}_{n}-v^{2})\psi_{i,K}(t)
+5\int_{\mathbb{R}}(v^{2}_{n,x}-v^{2}_{x})\psi_{i,K}(t)
+\int_{\mathbb{R}}(v^{2}_{n,xx}-v^{2}_{xx})\psi_{i,K}(t)\nonumber\\
&=K^{n}_{1}(t)+K^{n}_{2}(t)+K^{n}_{3}(t).
\label{R5}
\end{align}
Then it is easy to check that 
\begin{align}
|K^{n}_{1}(t)|&\le 4\int_{\mathbb{R}}|v_{n}-v|(v_{n}+v)\psi_{i,K}(t)\nonumber\\
&=4\|v_{n}-v\|_{L^{2}(\mathbb{R})}\|v_{n}+v\|_{L^{2}(\mathbb{R})}
\|\psi_{i,K}\|_{L^{\infty}(\mathbb{R})}\nonumber\\
&\le O(\|v_{n}-v\|_{L^{2}(\mathbb{R})})\nonumber\\
&\rightarrow 0~~\text{as}~~n\rightarrow+\infty,
\label{R6}
\end{align}
and
\begin{align}
|K^{n}_{2}(t)|&\le 5\int_{\mathbb{R}}|v_{n,x}-v_{x}|\left|v_{n,x}+v_{x}\right|\psi_{i,K}(t)\nonumber\\
&=4\|v_{n,x}-v_{x}\|_{L^{2}(\mathbb{R})}\|v_{n,x}+v_{x}\|_{L^{2}(\mathbb{R})}
\|\psi_{i,K}\|_{L^{\infty}(\mathbb{R})}\nonumber\\
&\le O(\|v_{n,x}-v_{x}\|_{L^{2}(\mathbb{R})})\nonumber\\
&\rightarrow 0~~\text{as}~~n\rightarrow+\infty.
\label{R7}
\end{align}
Recalling that $v_{xx}=4v-u$ and thus  $v^{2}_{xx}=16v^{2}+u^{2}-8uv$, we also get
\begin{align}
|K^{n}_{3}(t)|&\le 16\int_{\mathbb{R}}|v_{n}-v|(v_{n}+v)\psi_{i,K}(t)+\int_{\mathbb{R}}|u_{n}-u|(u_{n}+u)\psi_{i,K}(t)\nonumber\\
&\hspace{1cm}+8\int_{\mathbb{R}}u_{n}|v_{n}-v|\psi_{i,K}(t)+8\int_{\mathbb{R}}v|u_{n}-u|\psi_{i,K}(t)\nonumber\\
&\le 16\|v_{n}-v\|_{L^{2}(\mathbb{R})}\|v_{n}+v\|_{L^{2}(\mathbb{R})}\|\psi_{i,K}\|_{L^{\infty}(\mathbb{R})}
\nonumber\\
&\hspace{1cm}+\|u_{n}-u\|_{L^{2}(\mathbb{R})}\|u_{n}+u\|_{L^{2}(\mathbb{R})}\|\psi_{i,K}\|_{L^{\infty}(\mathbb{R})}\nonumber\\
&\hspace{2cm}+8\|u_{n}\|_{L^{2}(\mathbb{R})}\|v_{n}-v\|_{L^{2}(\mathbb{R})}
\|\psi_{i,K}\|_{L^{\infty}(\mathbb{R})}\nonumber\\
&\hspace{3cm}+8\|v\|_{L^{2}(\mathbb{R}}\|u_{n}-u\|_{L^{2}(\mathbb{R})}
\|\psi_{i,K}\|_{L^{\infty}(\mathbb{R})}\nonumber\\
&\le O\left(\|u_{n}-u\|_{L^{2}(\mathbb{R})})+O(\|v_{n}-v\|_{L^{2}(\mathbb{R})}\right)\nonumber\\
&\rightarrow 0~~\text{as}~~n\rightarrow+\infty.
\label{R8}
\end{align}
Combining \eqref{R5}-\eqref{R8}, we obtain \eqref{R4}.

Thanks to \eqref{R4}, the monotonicity formula \eqref{4.28} holds for any $u\in \mathcal{X}([0,T[)$, with $0<T\le+\infty$.
\hfill $ \square $ \vspace*{2mm}

\subsection{A localized and a global estimate}\label{Subsection 4.3}
\begin{figure}[ht]
\centering
\includegraphics[width=10cm, height=6cm]{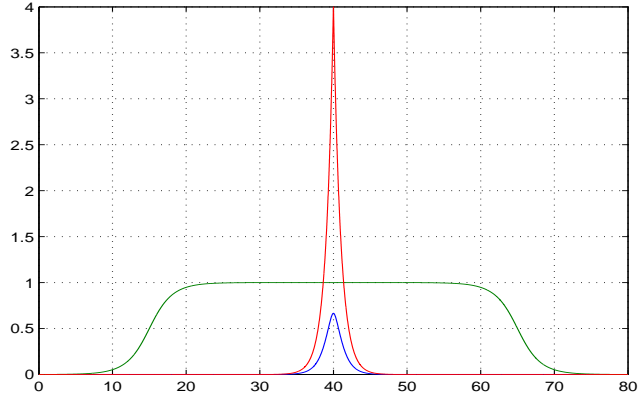}
\caption{Localization-function $\phi_{green}(x)=\psi(x-15)-\psi(x-65)$ (at time $t=10$) profile. Also, the peakon $4\varphi(x-40)$ and the smooth-peakon $4\rho(x-40)$ (at time $t=10$ with speed $c=4$) profiles. In this example, one can see that $\phi_{green}$ is close to $1$ on $]25,55[$, and decays exponentially to $0$ on $\mathbb{R}\setminus ]10,70[$.\label{de1}}
\end{figure}

Let $K=\sqrt{L}/8$ and define the function $\phi_{i}=\phi_{i}(t,x)$ (see Fig. \ref{de1}) by 
\begin{equation}
  \left\{
    \begin{aligned}
     &\phi_{1}=1-\psi_{2,K}=1-\psi_{K}(\cdot-y_{2}(t)),\\
     &\phi_{i}=\psi_{i,K}-\psi_{i+1,K}=\psi_{K}(\cdot-y_{i}(t))-\psi_{K}(\cdot-y_{i+1}(t)),\\
     &\phi_{N}=\psi_{N,K}=\psi_{K}(\cdot-y_{N}(t)),\\
    \end{aligned}
  \right.~~i=2,\ldots,N-1,
  \label{4.a1}
\end{equation}
where $\psi_{i,K}$'s and $y_{i}$'s are defined in Subsection \ref{Subsection 4.2}. One can see that the 
$\phi_{i}$'s are positive functions and that $\sum_{i=1}^{N}
\phi_{i}=1$. We take  $L/K>4$ so that $\phi_{i}$ satisfies for $i=1,\ldots,N$,
\begin{equation}
|1-\phi_{i}|\le 2e^{-\frac{L}{8K}}~~\text{on}~~\left]y_{i}+\frac{L}{8},y_{i+1}-\frac{L}{8}\right[
\label{4.67}
\end{equation}
and 
\begin{equation}
|\phi_{i}|\le 2e^{-\frac{L}{8K}}~~\text{on}~~\mathbb{R}\setminus\left]y_{i}-\frac{L}{8},y_{i+1}+\frac{L}{8}\right[.
\label{4.68}
\end{equation}

We will use the following localized version of the conservation laws defined for $i=1,\ldots,N$ by 
\begin{equation}
E_{i}(t)=\int_{\mathbb{R}}\left(4v^{2}+5v^{2}_{x}+v^{2}_{xx}\right)\phi_{i}(t)~~\text{and}~~
F_{i}(t)=\int_{\mathbb{R}}
\left(-v^{3}_{xx}+12vv^{2}_{xx}-48v^{2}v_{xx}+64v^{3}\right)\phi_{i}(t).
\label{4.69}
\end{equation}
One can remark that the functional $E_{i}(\cdot)$ and $F_{i}(\cdot)$ do not depend on time in the statement below since we fix $-\infty=y_{1}<y_{2}<\ldots<y_{N}<y_{N+1}=+\infty$.

For $i=1,\ldots,N$, we set $\Omega_{i}=]y_{i}-L/8,y_{i+1}+L/8[$. First, one can notice that  
\begin{equation}
\sum_{j=1}^{N}\rho_{c_{j}}(x-\tilde{x}_{j})=\rho_{c_{i}}(x-\tilde{x}_{i})+O(e^{-L/4}),~~
\forall x\in\Omega_{i},
\label{A1}
\end{equation}
we abuse notation by writing $\rho_{c_{i}}(x-\tilde{x}_{i})= O(e^{-L/4})$ for all $x\in\mathbb{R}\setminus\Omega_{i}$. We will now decompose this interval according to the variation of $v=(4-\partial^{2}_{x})^{-1}u$ in the same way as in Section \ref{Section 3}. We set 
\begin{equation}
\alpha_{i}=\sup\left\lbrace x<\tilde{x}_{i},~v(x)=\frac{c_{i}}{2400}\right\rbrace~~\text{and}~~
\beta_{i}=\inf\left\lbrace x>\tilde{x}_{i},~v(x)=\frac{c_{i}}{2400}\right\rbrace.
\label{4.70}
\end{equation}
According to Lemma \ref{Lemma 4.1}, we know that $v$ is close to $\sum_{i=1}^{N}\rho_{c}(\cdot-\tilde{x}_{i})$ in $L^{\infty}$ norm with $\rho_{c_{i}}(0)=c_{i}/6$. Therefore $v$ must have at least one local maximum on $[\alpha_{i},\beta_{i}]$. Assume that on $[\alpha_{i},\beta_{i}]$ the function $v$ admits $k_{i}+1$ points $(\xi^{i}_{j})_{j=1}^{k_{i}+1}$ with local maximal values for some integer $k_{i}\ge 0$, where $\xi^{i}_{1}$ is the first local maximum point and $\xi^{i}_{k_{i}+1}$ the last local maximum point\footnote{In the case of an infinite countable number of local maximal values, the proof is exactly the same.}. Then between $\xi^{i}_{1}$ and $\xi^{i}_{k_{i}+1}$, the function $v$ admits  $k_{i}$ points $(\eta^{i}_{j})_{j=1}^{k_{i}}$ with local minimal values.  We rename $\alpha_{i}=\eta^{i}_{0}$ and $\beta_{i}=\eta^{i}_{k+1}$ so that it holds 
\begin{equation}
\eta^{i}_{0}<\xi^{i}_{1}<\eta^{i}_{1}<\ldots<\xi^{i}_{j}<\eta^{i}_{j}<\xi^{i}_{j+1}<\eta^{i}_{j+1}
<\ldots<\eta^{i}_{k_{i}}<\xi^{i}_{k_{i}+1}<\eta^{i}_{k_{i}+1}.
\label{4.71}
\end{equation}
Let
\begin{equation}
M^{i}_{j}=v(\xi^{i}_{j}),~~j=1,\ldots,k_{i}+1,~~\text{and}~~m^{i}_{j}=v(\eta^{i}_{j}),~~j=1,\ldots,k_{i}.
\label{4.72}
\end{equation}
By construction
\begin{equation}
v_{x}(x)\ge 0,~~\forall x\in[\eta^{i}_{j-1},\xi^{i}_{j}],~~j=1,\ldots,k_{i}
\label{4.73}
\end{equation}
and
\begin{equation}
v_{x}(x)\le 0,~~\forall x\in[\xi^{i}_{j},\eta^{i}_{j}],~~j=1,\ldots,k_{i}+1.
\label{4.74}
\end{equation}
Proceeding as for \eqref{3.18}-\eqref{3.20}, we also have
\begin{equation}
v(x)\le\frac{c_{i}}{300},~~\forall x\in\Omega_{i}\setminus[\eta^{i}_{0},\eta^{i}_{k_{i}+1}],
\label{4.75}
\end{equation}
\begin{equation}
u(x)\le\frac{c_{i}}{300},~~\forall x\in\Omega_{i}\setminus[\eta^{i}_{0},\eta^{i}_{k_{i}+1}],
\label{4.76}
\end{equation}
and taking $L>L_{0}>8C_{0}$, it holds
\begin{equation}
[\eta^{i}_{0},\eta^{i}_{k_{i}+1}]\subset[\tilde{x}_{i}-C_{0},\tilde{x}_{i}+C_{0}]
\subset\left]y_{i}+\frac{L}{8},y_{i+1}-\frac{L}{8}\right[,
\label{4.77}
\end{equation}
where $C_{0}>0$ is the universal constant appearing in \eqref{3.20}.

We now derive versions of Lemma \ref{Lemma 3.4}, Lemma \ref{Lemma 3.5} and Lemma \ref{Lemma 3.6} where the global functional $ E(\cdot)$ and $F(\cdot) $ are replaced by their localized versions
 $E_i(\cdot) $ and $F_i(\cdot) $. Please note that, we will change the order of the extrema of $v=(4-\partial^{2}_{x})^{-1}u$ while keeping the same notations as in \eqref{4.72}.

\begin{Lem}[Connection Between $E_{i}(\cdot)$ and the Local Extrema of $v$]\label{Lemma 4.3}
Let $u\in H^{1}(\mathbb{R})$ and $v=(4-\partial^{2}_{x})^{-1}u\in H^{3}(\mathbb{R})$. For $i=1,\ldots,N$, define the function $g_{i}$ by
\begin{equation}
  g_{i}(x)=\left\{
    \begin{aligned}
     &2v+v_{xx}-3v_{x},~~x<\xi^{i}_{1},\\
     &2v+v_{xx}+3v_{x},~~\xi^{i}_{j}<x<\eta^{i}_{j},\\
     &2v+v_{xx}-3v_{x},~~\eta^{i}_{j}<x<\xi^{i}_{j+1},\\
     &2v+v_{xx}+3v_{x},~~x>\xi^{i}_{k_{i}+1},\\
    \end{aligned}
  \right.~~j=1,\ldots,k_{i}.
  \label{4.78}
\end{equation}
Then it holds
\begin{align}
\int_{\mathbb{R}}g^{2}_{i}(x)\phi_{i}(x)&=E_{i}(u)-12\left(\sum_{j=0}^{k_{i}}(M^{i}_{j+1})^{2}\phi_{i}(\xi^{i}_{j+1})-\sum_{j=1}^{k_{i}}(m^{i}_{j})^{2}\phi_{i}(\eta_{j})\right)
+\|u\|^{2}_{\mathcal{H}}O(L^{-1/2}).
\label{4.79}
\end{align}
\end{Lem}

\textbf{Proof.} We have 
\begin{equation}
\int_{\mathbb{R}}g^{2}_{i}(x)\phi_{i}(x)dx=\int_{-\infty}^{\xi^{i}_{1}}g^{2}_{i}(x)\phi_{i}(x)dx
+\sum_{j=1}^{k_{i}}\int_{\xi^{i}_{j}}^{\xi^{i}_{j+1}} g^{2}_{i}(x)\phi_{i}(x)dx
+\int_{\xi^{i}_{k_{i}+1}}^{+\infty}g^{2}_{i}(x)\phi_{i}(x)dx.
\label{4.80}
\end{equation}
For $j=1,\ldots,k_{i}$,
\begin{align*}
\int_{\xi^{i}_{j}}^{\xi^{i}_{j+1}} g^{2}_{i}(x)\phi_{i}(x)dx&=\int_{\xi^{i}_{j}}^{\eta^{i}_{j}}\left( 2v+v_{xx}+3v_{x}\right)^{2}\phi_{i}(x)dx+\int_{\eta^{i}_{j}}^{\xi^{i}_{j+1}}\left(2v+v_{xx}-3v_{x}\right)^{2}\phi_{i}(x)dx\\
&=J+I.
\end{align*}
Computing $I$ , we obtain
\begin{align}
I &=\int_{\eta^{i}_{j}}^{\xi^{i}_{j+1}}\left(4v^{2}+v^{2}_{xx}+9v^{2}_{x}+4vv_{xx}
-12vv_{x}-6v_{x}v_{xx}\right)\phi_{i}\nonumber\\
&=\int_{\eta^{i}_{j}}^{\xi^{i}_{j+1}}\left(4v^{2}+v^{2}_{xx}+9v^{2}_{x}\right)\phi_{i}+
4\int_{\eta^{i}_{j}}^{\xi^{i}_{j+1}}vv_{xx}\phi_{i}-12\int_{\eta^{i}_{j}}^{\xi^{i}_{j+1}}vv_{x}\phi_{i}
-6\int_{\eta^{i}_{j}}^{\xi^{i}_{j+1}}v_{x}v_{xx}\phi_{i}\nonumber\\
&=\int_{\eta^{i}_{j}}^{\xi^{i}_{j+1}}\left(4v^{2}+v^{2}_{xx}+9v^{2}_{x}\right)\phi_{i}+I_{1}+I_{2}
+I_{3}
\label{4.81}
\end{align}
with
\begin{align}
I_{1}&=-4\int_{\eta^{i}_{j}}^{\xi^{i}_{j+1}}\partial_{x}(v\phi_{i})v_{x}
=-4\int_{\eta^{i}_{j}}^{\xi^{i}_{j+1}}v^{2}_{x}\phi_{i}
-4\int_{\eta^{i}_{j}}^{\xi^{i}_{j+1}}vv_{x}\phi'_{i}=-4\int_{\eta^{i}_{j}}^{\xi^{i}_{j+1}}v^{2}_{x}\phi_{i}
-2\int_{\eta^{i}_{j}}^{\xi^{i}_{j+1}}(v^{2})_{x}\phi'_{i}\nonumber\\
&=-2v^{2}(\xi^{i}_{j+1})\phi'_{i}(\xi^{i}_{j+1})+2v^{2}(\eta^{i}_{j})\phi'_{i}(\eta^{i}_{j})-4\int_{\eta^{i}_{j}}^{\xi^{i}_{j+1}}v^{2}_{x}\phi_{i}
+2\int_{\eta^{i}_{j}}^{\xi^{i}_{j+1}}v^{2}\phi''_{i}
\label{4.82}
\end{align}
\begin{align}
I_{2}&=-6\int_{\eta^{i}_{j}}^{\xi^{i}_{j+1}}\partial_{x}(v^{2})\phi_{i}
=-6v^{2}(\xi^{i}_{j+1})\phi_{i}(\xi^{i}_{j+1})+6v^{2}(\eta^{i}_{j})\phi_{i}(\eta^{i}_{j})+6\int_{\eta^{i}_{j}}^{\xi^{i}_{j+1}}v^{2}\phi'_{i}
\label{4.83}
\end{align}
and
\begin{align}
I_{3}&=-3\int_{\eta^{i}_{j}}^{\xi^{i}_{j+1}}\partial_{x}(v^{2}_{x})\phi_{i}=3\int_{\eta^{i}_{j}}^{\xi^{i}_{j+1}}v^{2}_{x}\phi'_{i}.
\label{4.84}
\end{align}
Adding \eqref{4.81}-\eqref{4.84}, we get
\begin{align}
I&=\int_{\eta^{i}_{j}}^{\xi^{i}_{j+1}}\left(4v^{2}+5v^{2}_{x}+v^{2}_{xx}\right)\phi_{i}-6v^{2}(\xi^{i}_{j+1})\phi_{i}(\xi^{i}_{j+1})+6v^{2}(\eta^{i}_{j})\phi_{i}(\eta^{i}_{j})\nonumber\\
&\hspace{1cm}-2v^{2}(\xi^{i}_{j+1})\phi'_{i}(\xi^{i}_{j+1})+2v^{2}(\eta^{i}_{j})\phi'_{i}(\eta^{i}_{j})
+R_1,
\label{4.85}
\end{align}
where using that $K=\sqrt{L}/8$, we have 
$$
|R_{1}|\le 6(\|\phi'_{i}\|_{L^{\infty}(\mathbb{R})}+\|\phi''_{i}\|_{L^{\infty}(\mathbb{R})})
\int_{\eta^{i}_{j}}^{\xi^{i}_{j+1}}(v^{2}+v^{2}_{x})\le O(L^{-1/2})\int_{\eta^{i}_{j}}^{\xi^{i}_{j+1}}(v^{2}+v^{2}_{x}).
$$
Similar computations lead to 
\begin{align}
J&=\int_{\xi^{i}_{j}}^{\eta^{i}_{j}}\left(4v^{2}+5v^{2}_{x}+v^{2}_{xx}\right)\phi_{i}-6v^{2}(\xi^{i}_{j})\phi_{i}(\xi^{i}_{j})+6v^{2}(\eta^{i}_{j})\phi_{i}(\eta^{i}_{j})\nonumber\\
&\hspace{1cm}+2v^{2}(\xi^{i}_{j})\phi'_{i}(\xi^{i}_{j})-2v^{2}(\eta^{i}_{j})\phi'_{i}(\eta^{i}_{j})
+R_2,
\label{4.86}
\end{align}
\begin{equation}
\int_{-\infty}^{\xi_1^{i} }g^2(x)\phi_i (x)  \, dx =\int_{-\infty}^{\xi_1^{i}}\left(4v^{2}+5v^{2}_{x}+v^{2}_{xx}\right) \phi_i 
-6v^{2}(\xi^{i}_{1})\phi_{i}(\xi^{i}_{1})-2v^{2}(\xi^{i}_{1})\phi'_{i}(\xi^{i}_{1})+R_3 \label{4.855}
\end{equation}
and 
\begin{equation}
\int_{\xi_{k_i+1}^{i}}^{+\infty }g^2(x)\phi_i (x)  \, dx =\int_{\xi_{k_i+1}^{i}}^{+\infty}\left(4v^{2}+5v^{2}_{x}+v^{2}_{xx}\right) \phi_i 
-6v^{2}(\xi^{i}_{1})\phi_{i}(\xi_{k_i+1}^{i})+2v^{2}(\xi_{k_i+1}^{i})\phi'_{i}(\xi_{k_i+1}^{i})+R_4,\label{4.8555}
\end{equation}
with 
$$
|R_{2}|\le O(L^{-1/2})\int_{\xi^{i}_{j}}^{\eta^{i}_{j}}(v^{2}+v^{2}_{x}),
~~|R_{3}|\le O(L^{-1/2})\int_{-\infty}^{\xi^{i}_{1}}(v^{2}+v^{2}_{x})
~~\text{and}~~|R_{4}|\le O(L^{-1/2})\int_{\xi^{i}_{k_{i}+1}}^{+\infty}(v^{2}+v^{2}_{x}).
$$
Then, adding \eqref{4.85} and \eqref{4.86}, and summing over  $j\in\lbrace 1,\ldots,k_{i}\rbrace$, we infer that 
 \begin{align}
\int_{\xi^{i}_{1}}^{\xi^{i}_{k_{i}+1}} g^{2}_{i}(x)\phi_{i}(x)dx
&=\int_{\xi^{i}_{1}}^{\xi^{i}_{k_{i}+1}}\left(4v^{2}+5v^{2}_{x}+v^{2}_{xx}\right)\phi_{i}
-6\sum_{j=1}^{k_{i}}v^{2}(\xi^{i}_{j+1})\phi_{i}(\xi^{i}_{j+1})-6\sum_{j=1}^{k_{i}}v^{2}(\xi^{i}_{j})\phi_{i}(\xi^{i}_{j})\nonumber\\
&\hspace{0.5cm}+12\sum_{j=1}^{k_{i}}v^{2}(\eta^{i}_{j})\phi_{i}(\eta^{i}_{j})
-2\sum_{j=1}^{k_{i}}v^{2}(\xi^{i}_{j+1})\phi'_{i}(\xi^{i}_{j+1})
+2\sum_{j=1}^{k_{i}}v^{2}(\xi^{i}_{j})\phi'_{i}(\xi^{i}_{j})+R,
\label{4.87}
\end{align}
with
$$|R|\le O(L^{-1/2})\int_{\xi^{i}_{1}}^{\xi^{i}_{k_{i}+1}}(v^{2}+v^{2}_{x}).$$
Finally, adding \eqref{4.855}-\eqref{4.87}, and recalling that $ \|v\|_{H^{1}}\le\|u\|_{\mathcal H} $, we obtain the lemma.
\hfill $ \square $ \vspace*{2mm}

\begin{Lem}[Connection Between $F_{i}(\cdot)$ and the Local Extrema of $v$]\label{Lemma 4.4}
Let $u\in H^{1}(\mathbb{R})$ and $v=(4-\partial^{2}_{x})^{-1}u\in H^{3}(\mathbb{R})$. For $i=1,\ldots,N$, define the function $h_{i}$ by
\begin{equation}
  h_{i}(x)=\left\{
    \begin{aligned}
     &-v_{xx}-6v_{x}+16v,~~x<\xi^{i}_{1},\\
     &-v_{xx}+6v_{x}+16v,~~\xi^{i}_{j}<x<\eta^{i}_{j},\\
     &-v_{xx}-6v_{x}+16v,~~\eta^{i}_{j}<x<\xi^{i}_{j+1},\\
     &-v_{xx}+6v_{x}+16v,~~x>\xi^{i}_{k_{i}+1},\\
    \end{aligned}
  \right.~~j=1,\ldots,k_{i}.
  \label{4.88}
\end{equation}
Then it holds
\begin{align}
\int_{\mathbb{R}}h_{i}(x)g^{2}_{i}(x)\phi_{i}(x)&=F_{i}(u)-144\left(\sum_{j=0}^{k_{i}}(M^{i}_{j+1})^{3}\phi_{i}(\xi^{i}_{j+1})-\sum_{j=1}^{k_{i}}(m^{i}_{j})^{3}\phi_{i}(\eta^{i}_{j})\right)+\|u\|^{3}_{\mathcal{H}}O(L^{-1/2}).
\label{4.89}
\end{align}
\end{Lem}

\textbf{Proof.} We have 
\begin{align}
\int_{\mathbb{R}}h_{i}(x)g^{2}_{i}(x)\phi_{i}(x)dx&=\int_{-\infty}^{\xi^{i}_{1}}
h_{i}(x)g^{2}_{i}(x)\phi_{i}(x)dx
+\sum_{j=1}^{k_{i}}\int_{\xi^{i}_{j}}^{\xi^{i}_{j+1}}h_{i}(x)g^{2}_{i}(x)\phi_{i}(x)dx\nonumber\\
&\hspace{1cm}+\int_{\xi^{i}_{k_{i}+1}}^{+\infty}h_{i}(x)g^{2}_{i}(x)\phi_{i}(x)dx.
\label{4.90}
\end{align}
For $j=1,\ldots,k_{i}$, 
\begin{align*}
\int_{\xi^{i}_{j}}^{\xi^{i}_{j+1}}h_{i}(x)g^{2}_{i}(x)\phi_{i}(x)dx&=\int_{\xi^{i}_{j}}^{\eta^{i}_{j}}\left(-v_{xx}-6v_{x}+16v\right)\left(2v+v_{xx}
-3v_{x}\right)^{2}\phi_{i}\\
&\hspace{1cm}+\int_{\eta^{i}_{j}}^{\xi^{i}_{j+1}}\left(-v_{xx}+6v_{x}+16v\right)\left(2v+v_{xx}
+3v_{x}\right)^{2}\phi_{i}\\
&=J+I.
\end{align*}
Computing $I$, we obtain 
\begin{align}
I&=\int_{\eta^{i}_{j}}^{\xi^{i}_{j+1}}\left(-v^{3}_{xx}+12vv^{2}_{xx}
+64v^{3}+60v^{2}v_{xx}\right)\phi_{i}-54\int_{\eta^{i}_{j}}^{\xi^{i}_{j+1}}v^{3}_{x}\phi_{i}
+27\int_{\eta^{i}_{j}}^{\xi^{i}_{j+1}}v^{2}_{x}v_{xx}\phi_{i}\nonumber\\
&\hspace{1cm}-108\int_{\eta^{i}_{j}}^{\xi^{i}_{j+1}}vv_{x}v_{xx}\phi_{i}
-216\int_{\eta^{i}_{j}}^{\xi^{i}_{j+1}}v^{2}v_{x}\phi_{i}
+216\int_{\eta^{i}_{j}}^{\xi^{i}_{j+1}}vv^{2}_{x}\phi_{i}\nonumber\\
&=\int_{\eta^{i}_{j}}^{\xi^{i}_{j+1}}\left(-v^{3}_{xx}+12vv^{2}_{xx}
+64v^{3}+60v^{2}v_{xx}\right)\phi_{i}-54\int_{\eta^{i}_{j}}^{\xi^{i}_{j+1}}v^{3}_{x}\phi_{i}
+I_{1}+I_{2}+I_{3}+I_{4}
\label{4.91}
\end{align}
with 
\begin{align}
I_{1}&=9\int_{\eta^{i}_{j}}^{\xi^{i}_{j+1}}\partial_{x}(v^{3}_{x})\phi_{i}=-9\int_{\eta^{i}_{j}}^{\xi^{i}_{j+1}}v^{3}_{x}\phi'_{i}\, ,
\label{4.92}
\end{align}
\begin{align}
I_{2}&=-54\int_{\eta^{i}_{j}}^{\xi^{i}_{j+1}}v\partial_{x}(v^{2}_{x})\phi_{i}=54\int_{\eta^{i}_{j}}^{\xi^{i}_{j+1}}\partial_{x}(v\phi_{i})v^{2}_{x}=54\int_{\eta^{i}_{j}}^{\xi^{i}_{j+1}}v^{3}_{x}\phi_{i}
+54\int_{\eta^{i}_{j}}^{\xi^{i}_{j+1}}vv^{2}_{x}\phi'_{i}\, ,
\label{4.93}
\end{align}
\begin{align}
I_{3}&=-72\int_{\eta^{i}_{j}}^{\xi^{i}_{j+1}}\partial_{x}(v^{3})\phi_{i}=-72v^{3}(\xi^{i}_{j+1})\phi_{i}(\xi^{i}_{j+1})+72v^{3}(\eta^{i}_{j})\phi_{i}(\eta^{i}_{j})+72\int_{\eta^{i}_{j}}^{\xi^{i}_{j+1}}v^{3}\phi_{i}'\, , 
\label{4.94}
\end{align}
and 
\begin{align}
I_{4}&=108\int_{\eta^{i}_{j}}^{\xi^{i}_{j+1}}\partial_{x}(v^{2})v_{x}\phi_{i}=-108\int_{\eta^{i}_{j}}^{\xi^{i}_{j+1}}v^{2}\partial_{x}(v_{x}\phi_{i})=-108\int_{\eta^{i}_{j}}^{\xi^{i}_{j+1}}v^{2}v_{xx}\phi_{i}-108\int_{\eta^{i}_{j}}^{\xi^{i}_{j+1}}v^{2}
v_{x}\phi'_{i}\nonumber\\
&=-108\int_{\eta^{i}_{j}}^{\xi^{i}_{j+1}}v^{2}v_{xx}\phi_{i}-36\int_{\eta^{i}_{j}}^{\xi^{i}_{j+1}}
\partial_{x}(v^{3})\phi'_{i}\nonumber\\
&=-36v^{3}(\xi^{i}_{j+1})\phi'_{i}(\xi^{i}_{j+1})+36v^{3}(\eta^{i}_{j})\phi'_{i}(\eta^{i}_{j})-108\int_{\eta^{i}_{j}}^{\xi^{i}_{j+1}}v^{2}v_{xx}\phi_{i}
+36\int_{\eta^{i}_{j}}^{\xi^{i}_{j+1}}v^{3}\phi''_{i}.
\label{4.95}
\end{align}
Adding \eqref{4.91}-\eqref{4.95}, we get
\begin{align}
I&=\int_{\eta^{i}_{j}}^{\xi^{i}_{j+1}}\left(-v^{3}_{xx}+12vv^{2}_{xx}
+64v^{3}-48v^{2}v_{xx}\right)\phi_{i}-72v^{3}(\xi^{i}_{j+1})\phi_{i}(\xi^{i}_{j+1})+72v^{3}(\eta^{i}_{j})\phi_{i}(\eta^{i}_{j})\nonumber\\
&\hspace{1cm}-36v^{3}(\xi^{i}_{j+1})\phi'_{i}(\xi^{i}_{j+1})
+36v^{3}(\eta^{i}_{j})\phi'_{i}(\eta^{i}_{j})+R,
\label{4.96}
\end{align}
where using that $\|v\|_{C^{1}(\mathbb{R})}\le C'_{S}\|v\|_{H^{2}(\mathbb{R})}$ (with $C'_{S}$ the constant of Sobolev), and 
$\|v\|_{H^{2}(\mathbb{R})}\sim\|u\|_{\mathcal{H}}$, the estimate of $R$ leads to 
\begin{align*}
|R|&\le (\|\phi'_{i}\|_{L^{\infty}(\mathbb{R})}+\|\phi''_{i}\|_{L^{\infty}(\mathbb{R})})
(\|v\|_{L^{\infty}(\mathbb{R})}+\|v_{x}\|_{L^{\infty}(\mathbb{R})})
\int_{\eta^{i}_{j}}^{\xi^{i}_{j+1}}(v^{2}+v^{2}_{x})\\
&\le O(L^{-1/2})\|u\|_{\mathcal{H}}\int_{\eta^{i}_{j}}^{\xi^{i}_{j+1}}(v^{2}+v^{2}_{x}).
\end{align*}
Similar computations lead to 
\begin{align}
J&=\int_{\xi^{i}_{j}}^{\eta^{i}_{j}}\left(-v^{3}_{x}+12vv^{2}_{xx}
+64v^{3}-48v^{2}v_{xx}\right)\phi_{i}-72v^{3}(\xi^{i}_{j})\phi_{i}(\xi^{i}_{j})
+72v^{3}(\eta^{i}_{j})\phi_{i}(\eta^{i}_{j})\nonumber\\
&\hspace{1cm}+36v^{3}(\xi^{i}_{j})\phi'_{i}(\xi^{i}_{j})
-36v^{3}(\eta^{i}_{j})\phi'_{i}(\eta^{i}_{j})
+ O(L^{-1/2})\|u\|_{\mathcal{H}}\int_{\xi^{i}_{j}}^{\eta^{i}_{j}}(v^{2}+v^{2}_{x}), \label{4.97}
\end{align}
\begin{align}
\int_{-\infty}^{\xi^{i}_{1}}h_{i}(x)g^{2}_{i}(x)\phi_{i}(x)
&=\int_{-\infty}^{\xi^{i}_{1}}\left(-v^{3}_{xx}+12vv^{2}_{xx}
+64v^{3}-48v^{2}v_{xx}\right)\phi_{i}-72v^{3}(\xi^{i}_{1})\phi_{i}(\xi^{i}_{1})\nonumber\\
&\hspace{1cm}-36v^{3}(\xi^{i}_{1})\phi'_{i}(\xi^{i}_{1}) 
+O(L^{-1/2})\|u\|_{\mathcal{H}}\int_{-\infty}^{\xi^{i}_{1}}(v^{2}+v^{2}_{x})
\label{4.966}
\end{align}
and 
\begin{align}
\int_{\xi^{i}_{k_{i}+1}}^{+\infty}h_{i}(x)g^{2}_{i}(x)\phi_{i}(x)
&=\int_{\xi^{i}_{k_{i}+1}}^{+\infty}\left(-v^{3}_{xx}+12vv^{2}_{xx}
+64v^{3}-48v^{2}v_{xx}\right)\phi_{i}-72v^{3}(\xi^{i}_{k_{i}+1})\phi_{i}(\xi^{i}_{k_{i}+1})\nonumber\\
&\hspace{1cm}  +36v^{3}(\xi^{i}_{k_{i}+1})\phi'_{i}(\xi^{i}_{k_{i}+1})
+O(L^{-1/2})\|u\|_{\mathcal{H}}\int_{\xi^{i}_{k_{i}+1}}^{+\infty}(v^{2}+v^{2}_{x}). 
\label{4.977}
\end{align}
Adding \eqref{4.96} and \eqref{4.97}, and summing over $j\in\lbrace 1,\ldots,k_{i}\rbrace$, we get 
\begin{align}
&\int_{\xi^{i}_{1}}^{\xi^{i}_{k_{i}+1}}h_{i}(x)g^{2}_{i}(x)\phi_{i}(x)dx\nonumber\\
&\hspace{1cm}=\int_{\xi^{i}_{1}}^{\xi^{i}_{k_{i}+1}}\left(-v^{3}_{x}+12vv^{2}_{xx}
+64v^{3}-48v^{2}v_{xx}\right)\phi_{i}
-72\sum_{j=1}^{k_{i}}v^{3}(\xi^{i}_{j+1})\phi_{i}(\xi^{i}_{j+1})\nonumber\\
&\hspace{2cm}-72\sum_{j=1}^{k_{i}}v^{3}(\xi^{i}_{j})\phi_{i}(\xi^{i}_{j})+144\sum_{j=1}^{k_{i}}v^{3}(\eta^{i}_{j})\phi_{i}(\eta^{i}_{j})
-36\sum_{j=1}^{k_{i}}v^{3}(\xi^{i}_{j+1})\phi'_{i}(\xi^{i}_{j+1})\nonumber\\
&\hspace{3cm}+36\sum_{j=1}^{k_{i}}v^{3}(\xi^{i}_{j})\phi'_{i}(\xi^{i}_{j})
+O(L^{-1/2})\|u\|_{\mathcal{H}}\int_{\xi^{i}_{1}}^{\xi^{i}_{k_{i}+1}}(v^{2}+v^{2}_{x}).
\label{4.98}
\end{align}
Finally, adding \eqref{4.966}-\eqref{4.98}, we obtain the lemma.
\hfill $ \square $ \vspace*{2mm}

\begin{Lem}[Connection Between $E_{i}(\cdot)$ and $F_{i}(\cdot)$]\label{Lemma 4.5}
Let $u\in H^{1}(\mathbb{R})$, with $y=(1-\partial^{2}_{x})u\in\mathcal{M}^{+}(\mathbb{R})$, that satisfies \eqref{4.2}. Let be given $N-1$ real numbers $-\infty=y_{1}<y_{2}<\ldots<y_{N}<y_{N+1}=+\infty$ with $y_{i}-y_{i-1}\ge 2L/3$. For $i=1,\ldots,N$, assume that $v=(4-\partial^{2}_{x})^{-1}u$ satisfies \eqref{4.70}-\eqref{4.77}, with local extrema on $[\eta^{i}_{0},\eta^{i}_{k_{i}+1}]$ arranged in decreasing order in the following way: 
\begin{equation}
M^{i}_{1}\ge M^{i}_{2}\ge\ldots\ge M^{i}_{k_{i}+1}\ge 0,~~m^{i}_{1}\ge m^{i}_{2}\ge\ldots\ge m^{i}_{k_{i}}\ge 0,~~M^{i}_{j+1}\ge m^{i}_{j},~~j=1,\ldots,k_{i}.
\label{4.99}
\end{equation}
There exist $\gamma_{0}>0$ and $L_{0}>0$ only depending on the speeds $(c_{i})_{i=1}^{N}$, such that if 
$0<\gamma<\gamma_{0}$ and $L>L_{0}>0$, then defining the functional $E_{i}(\cdot)$'s and $F_{i}(\cdot)$'s as in \eqref{4.a1}-\eqref{4.69}, it holds
\begin{equation}
F_{i}(u)\le 18M^{i}_{1}E_{i}(u)-72(M^{i}_{1})^{3}
+\|u\|^{3}_{\mathcal{H}}O(L^{-1/2}),~~i=1,\ldots,N.
\label{4.100}
\end{equation}
\end{Lem}

\textbf{Proof.} 
Combining \eqref{4.67}, \eqref{4.77} and \eqref{4.79} with $K=\sqrt{L}/8$, we get
\begin{equation}
\int_{\mathbb{R}}g^{2}_{i}(x)\phi_{i}(x)
=E_{i}(u)-12\left(\sum_{j=0}^{k_{i}}(M^{i}_{j+1})^{2}-\sum_{j=1}^{k_{i}}(m^{i}_{j})^{2}\right)
+\|u\|^{2}_{\mathcal{H}}O(L^{-1/2}).
\label{4.101}
\end{equation}
Similarly, combining \eqref{4.67}, \eqref{4.77} and \eqref{4.89}, we get
\begin{equation}
\int_{\mathbb{R}}h_{i}(x)g^{2}_{i}(x)\phi_{i}(x)dx =F_{i}(u)-144\left(\sum_{j=0}^{k_{i}}(M^{i}_{j+1})^{3}-\sum_{j=1}^{k_{i}}(m^{i}_{j})^{3}\right)
+\|u\|^{3}_{\mathcal{H}}O(L^{-1/2}).
\label{4.104}
\end{equation}
Now, let us show that $h_{i}\le 18 M^{i}_{1}$ on $\Omega_{i}$. Note that by \eqref{4.6} and \eqref{A1}, one can check that 
$18M^{i}_{1}\ge c_{i}/4$. We rewrite the function $h_{i}$ as
\begin{equation*}
  h_{i}(x)=\left\{
    \begin{aligned}
     &-v_{xx}-6v_{x}+16v,~~x<\eta^{i}_{0},\\
     &-\left(\partial^{2}_{x}+3\partial_{x}+2\right)v-3v_{x}+18v,~~\eta^{i}_{0}<x<\xi^{i}_{1},\\
     &-\left(\partial^{2}_{x}-3\partial_{x}+2\right)v+3v_{x}+18v,~~\xi^{i}_{j}<x<\eta^{i}_{j},\\
     &-\left(\partial^{2}_{x}+3\partial_{x}+2\right)v-3v_{x}+18v,~~\eta^{i}_{j}<x<\xi^{i}_{j+1},\\
     &-\left(\partial^{2}_{x}-3\partial_{x}+2\right)v+3v_{x}+18v,~~\xi^{i}_{k_{i}+1}<x<\eta^{i}_{k_{i}+1},\\
     &-v_{xx}+6v_{x}+16v,~~x>\eta^{i}_{k_{i}+1},\\
    \end{aligned}
  \right.~~j=1,\ldots,k_{i}.
\end{equation*}
Then, if $x\in\Omega_{i}\setminus[\eta^{i}_{0},\eta^{i}_{k_{i}+1}]$, using that $v_{xx}=4v-u$, \eqref{3.36}, \eqref{4.77} and \eqref{4.78}, it holds 
$$h_{i}\le|v_{xx}|+6|v_{x}|+16v
\le u+32v\\
\le\frac{c_{i}}{9}.$$
If $\eta^{i}_{0}<x<\xi^{i}_{1}$, then $v_{x}\ge 0$, and using that $y=(1-\partial^{2}_{x})u\ge 0$, it follows from Lemma \ref{Lemma 2.2} that 
\begin{align*}
h_{i}&=-(\partial^{2}_{x}+3\partial_{x}+2)v-3v_{x}+18v\\
&=-(2+\partial_{x})(4-\partial^{2}_{x})^{-1}(1+\partial_{x})u-3v_{x}+18v\\
&\le 18v.
\end{align*}
If $\xi^{i}_{j}<x<\eta^{i}_{j}$, then $v_{x}\le 0$, and similarly using that $y=(1-\partial^{2}_{x})u\ge 0$, it follows from Lemma \ref{Lemma 2.2} that 
\begin{align*}
h_{i}&=-(\partial^{2}_{x}-3\partial_{x}+2)v+3v_{x}+18v\\
&=-(2-\partial_{x})(4-\partial^{2}_{x})^{-1}(1-\partial_{x})u+3v_{x}+18v\\
&\le 18v.
\end{align*}
Therefore, it holds 
\begin{equation}
h_{i}(x)\le 18\max_{x\in\Omega_{i}}v(x)=18M^{i}_{1},~~\forall x\in\Omega_{i}.
\label{4.105}
\end{equation}
Now, taking $\phi_{i}\equiv 1$ on $\mathbb{R}$ in \eqref{4.79}, we have 
$\|g_{i}\|_{L^{2}(\mathbb{R})}\le\|u\|_{\mathcal{H}}$. Also, from the definition of $h_{i}$, and using \eqref{2.3} and
Remark \ref{Remark 4.1}, we have 
$\|h_{i}\|_{L^{\infty}(\mathbb{R})}\le \|u\|_{L^{\infty}(\mathbb{R})}+32\|v\|_{L^{\infty}(\mathbb{R})}
\le O(\|u\|_{\mathcal{H}})$. Then, combining \eqref{4.101}-\eqref{4.105}, we obtain
\begin{align*}
F_{i}(u)&-144\left(\sum_{j=0}^{k_{i}}(M^{i}_{j+1})^{3}-\sum_{j=1}^{k_{i}}(m^{i}_{j})^{3}\right)\\
&=\int_{\mathbb{R}}h_{i}(x)g^{2}_{i}(x)\phi_{i}(x)dx
+\|u\|^{3}_{\mathcal{H}}O(L^{-1/2})\\
&=\int_{\Omega_{i}}h_{i}(x)g^{2}_{i}(x)\phi_{i}(x)dx+\int_{\Omega^{c}_{i}}h_{i}(x)g^{2}_{i}(x)\phi_{i}(x)dx
+\|u\|^{3}_{\mathcal{H}}O(L^{-1/2})\\
&\le 18M^{i}_{1}\int_{\Omega_{i}}g^{2}_{i}(x)\phi_{i}(x)dx+\|h_{i}\|_{L^{\infty}(\mathbb{R})}\|g_{i}\|^{2}_{L^{2}(\mathbb{R})}\|\phi_{i}\|_{L^{\infty}(\Omega^{c}_{i})}
+\|u\|^{3}_{\mathcal{H}}O(L^{-1/2})\\
&\le 18M^{i}_{1}\left[E_{i}(u)-12\left(\sum_{j=0}^{k_{i}}(M^{i}_{j+1})^{2}-\sum_{j=1}^{k_{i}}(m^{i}_{j})^{2}\right)\right]
+\|u\|^{3}_{\mathcal{H}}O(L^{-1/2}).
\end{align*}
Therefore, using that $M^{i}_{j+1}\ge m^{i}_{j}$ and proceeding as in Lemma \ref{Lemma 3.6} (see \eqref{3.39}), we infer that
\begin{align*}
F_{i}(u)
&\le 18M^{i}_{1}E_{i}(u)-72(M^{i}_{1})^{3}
+144\sum_{j=1}^{k_{i}}\left\lbrace\left[(M^{i}_{j+1})^{3}-(m^{i}_{j})^{3}\right]
-\frac{3}{2}M^{i}_{1}\left[(M^{i}_{j+1})^{2}-(m^{i}_{j})^{2}\right]\right\rbrace \nonumber\\
&\hspace{1cm}+\|u\|^{3}_{\mathcal{H}}O(L^{-1/2})\nonumber\\
&\le 18M^{i}_{1}E_{i}(u)-72(M^{i}_{1})^{3}
+\|u\|^{3}_{\mathcal{H}}O(L^{-1/2}).
\label{4.106}
\end{align*}
This proves the lemma.
\hfill $ \square $ \vspace*{2mm} 

The lemma below is the generalization of Lemma \ref{Lemma 3.3}.

\begin{Lem}[General Quadratic Identity]
Let $Z=(z_{i})_{i=1}^{N}\in\mathbb{R}^{N}$ with $|z_{i}-z_{i-1}|\ge L/2$, and $u\in L^{2}(\mathbb{R})$. It holds
\begin{equation}
E(u)-\sum_{i=1}^{N}E(\varphi_{c_{i}})=\|u-S_{Z}\|^{2}_{\mathcal{H}}+4\sum_{i=1}^{N}c_{i}\left(v(z_{i})-\frac{c_{i}}{6}\right)+O(e^{-L/4}),
\label{4.107}
\end{equation}
where $S_{Z}$ is defined in \eqref{4.12} and $O(\cdot)$ only depends on $(c_{i})_{i=1}^{N}$.
\end{Lem}

\textbf{Proof.}
Let us compute 
\begin{align}
\|u-S_{Z}\|^{2}_{\mathcal{H}} &=\int_{\mathbb{R}}[(1-\partial^{2}_{x})(u-S_{Z})]
[(4-\partial^{2}_{x})^{-1}(u-S_{Z})]\nonumber\nonumber\\
&=\|u\|^{2}_{\mathcal{H}}+\|S_{Z}\|^{2}_{\mathcal{H}}-2\int_{\mathbb{R}}
[(1-\partial^{2}_{x})S_{Z}][(4-\partial^{2}_{x})^{-1}u]\nonumber\\
&=\|u\|^{2}_{\mathcal{H}}+\|S_{Z}\|^{2}_{\mathcal{H}}-2\sum_{i=1}^{N}c_{i}\int_{\mathbb{R}}
[(1-\partial^{2}_{x})\varphi_{c_{i}}(x-z_{i})]v\nonumber\\
&=\|u\|^{2}_{\mathcal{H}}+\|S_{Z}\|^{2}_{\mathcal{H}}-4\sum_{i=1}^{N}c_{i}v(z_{i}),
\label{4.108}
\end{align}
where we use that 
$$(1-\partial^{2}_{x})\varphi_{c_{i}}(\cdot-z_{i})=2c_{i}\delta_{z_{i}}$$
with $\delta_{z_{i}}$ the Dirac mass applied at point $z_{i}$. We also have 
\begin{align}
\|S_{Z}\|^{2}_{\mathcal{H}}&=\int_{\mathbb{R}}[(1-\partial^{2}_{x})S_{Z}][(4-\partial^{2}_{x})^{-1}S_{Z}]\nonumber\\
&=\sum_{i=1}^{N}\int_{\mathbb{R}}[(1-\partial^{2}_{x})\varphi_{c_{i}}(x-z_{i})][(4-\partial^{2}_{x})^{-1}S_{Z}]\nonumber\\
&=2\sum_{i=1}^{N}c_{i}\left\langle\delta_{z_{i}},(4-\partial^{2}_{x})^{-1}S_{Z}\right\rangle_{H^{-1},H^{1}}\nonumber\\
&=2\sum_{1\le i,j\le N}c_{i}\left\langle\delta_{z_{i}},
(4-\partial^{2}_{x})^{-1}\varphi_{c_{j}}(\cdot-z_{j})\right\rangle_{H^{-1},H^{1}}\nonumber\\
&=2\sum_{i=1}^{N}c_{i}(4-\partial^{2}_{x})^{-1}\varphi_{c_{i}}(0)+2\sum\limits_{\underset{i \neq j}{1\le i,j\le N}}c_{i}c_{j}(4-\partial^{2}_{x})^{-1}e^{-|z_{i}-z_{j}|},
\label{4.109}
\end{align}
where $\langle\cdot,\cdot\rangle_{H^{-1},H^{1}}$ denote the duality $H^{-1}/H^{1}$, and we recall that 
$\delta_{z_{i}}\in H^{-1}(\mathbb{R})$ since $\|\delta_{z_{i}}\|_{H^{-1}(\mathbb{R})}\le C_{S}$, with $C_{S}$ the constant appearing in \eqref{2.3}. Now, using that $|z_{i}-z_{i-1}|\ge L/2$,
\begin{align}
(4-\partial^{2}_{x})^{-1}e^{-|z_{i}-z_{j}|}&=\frac{1}{4}\int_{\mathbb{R}}e^{-2|x'-(z_{i}-z_{j})|}
e^{-|x'|}dx'\nonumber\\
&=\frac{1}{3}e^{-|z_{i}-z_{j}|}-\frac{1}{6}e^{-2|z_{i}-z_{j}|}\nonumber\\
&=O(e^{-L/4}),
\label{4.110}
\end{align} 
and combining \eqref{4.109} and \eqref{4.110}, for $L>L_{0}>0$ with $L_{0}\gg 1$, we get
\begin{equation}
\|S_{Z}\|^{2}_{\mathcal{H}}=\sum_{i=1}^{N}E(\varphi_{c_{i}})+O(e^{-L/4}).
\label{4.111}
\end{equation}
Finally, combining \eqref{4.108} and \eqref{4.111}, we obtain the lemma.
\hfill $ \square $ \vspace*{2mm} 

The last lemma is the localized version of Lemma \ref{Lemma 3.1}.

\begin{Lem}[Control of the Distances Between Local and Global Energies at $t=0$]\label{Lemma 4.7}
Let $u_0\in H^{1}(\mathbb{R})$ satisfying \eqref{1.11}-\eqref{1.13}. Then it holds
\begin{equation}
\left|E(u_{0})-\sum_{i=1}^{N}E(\varphi_{c_{i}})\right|\le O(\varepsilon^{2}) +O(e^{-L/4}),
\label{4.112}
\end{equation}
\begin{equation}
|E_{i}(u_{0})-E(\varphi_{c_{i}})|\le O(\varepsilon^{2})+O(e^{-\sqrt{L}}),~~i=1,\ldots,N,
\label{4.113}
\end{equation}
and 
\begin{equation}
|F_{i}(u_{0})-F(\varphi_{c_{i}})|\le O(\varepsilon^{2})+O(e^{-\sqrt{L}}),~~i=1,\ldots,N,
\label{4.114}
\end{equation}
where $O(\cdot)$ only depend on $(c_{i})_{i=1}^{N}$.
\end{Lem}

\textbf{Proof.} 
For the first estimate, applying triangular inequality and \eqref{1.13}, we have
\begin{align}
|E(u_{0})-E(S_{Z^{0}})|&=\left|\|u_{0}\|_{\mathcal{H}}-\|S_{Z^{0}}\|_{\mathcal{H}}\right|\left(\|u_{0}\|_{\mathcal{H}}+\|S_{Z^{0}}\|_{\mathcal{H}}\right)\nonumber\\
&\le\|u_{0}-S_{Z^{0}}\|_{\mathcal{H}}\left(\|u_{0}-S_{Z^{0}}\|_{\mathcal{H}}+2\|S_{Z^{0}}\|_{\mathcal{H}}\right)\nonumber\\
&\le\varepsilon^{2}\left(\varepsilon^{2}+\frac{2}{\sqrt{3}}\sum_{i=1}^{N}c_{i}\right).
\label{4.115}
\end{align}
Thus, combining \eqref{4.111} and \eqref{4.115}, it holds
\begin{align*}
\left|E(u_{0})-\sum_{i=1}^{N}E(\varphi_{c_{i}})\right|& \le\left|E(u_{0})-E(S_{Z^{0}})\right|
+\left|E(S_{Z^{0}})-\sum_{i=1}^{N}E(\varphi_{c_{i}})\right|\\
&\le\varepsilon^{2}(\varepsilon^{2}+O(1))+O(e^{-L/4})\\
&\le O(\varepsilon^{2})+O(e^{-L/4}).
\end{align*}

For the second estimate, using the exponential decay of $\varphi_{c_{i}}$'s and the $\phi_{i}$'s, and the definition of $E_{i}(\cdot)$, we have
\begin{align*}
|E_{i}(u_{0})&-E(\varphi_{c_{i}})|\\
&\le\left|\|u_{0}\|^{2}_{\mathcal{H}(\Omega_{i})}-\|\varphi_{c_{i}}\|^{2}_{\mathcal{H}(\Omega_{i})}\right|+O(e^{-\sqrt{L}})\\
&=\left|\|u_{0}\|_{\mathcal{H}(\Omega_{i})}-\|\varphi_{c_{i}}\|_{\mathcal{H}(\Omega_{i})}\right|\left(\|u_{0}\|_{\mathcal{H}(\Omega_{i})}+\|\varphi_{c_{i}}\|_{\mathcal{H}(\Omega_{i})}\right)+O(e^{-\sqrt{L}})\\
&\le\left(\|u_{0}-S_{Z^{0}}\|_{\mathcal{H}(\Omega_{i})}+\sum\limits_{\underset{j \neq i}{1\le j\le N}}\|\varphi_{c_{j}}\|_{\mathcal{H}(\Omega_{i})}\right)
\left(\|u_{0}-S_{Z^{0}}\|_{\mathcal{H}}+\frac{2}{\sqrt{3}}\sum_{j=1}^{N}c_{j}\right)+O(e^{-\sqrt{L}})\\
&\le\left(\varepsilon^{2}+O(e^{-L/8})\right)\left(\varepsilon^{2}+O(1)\right)+O(e^{-\sqrt{L}})\\
&\le O(\varepsilon^{2})+O(e^{-\sqrt{L}}).
\end{align*}

Similarly, for the third estimate, using the exponential decay of $\varphi_{c_{i}}$'s and the $\phi_{i}$'s, and the definition of $F_{i}(\cdot)$, we have 
\begin{align*}
|F_{i}(u_{0})&-F(\varphi_{c_{i}})|\\
&\le\left|\int_{\Omega_{i}}\left(u^{3}_{0}-\varphi^{3}_{c_{i}}\right)\right|+O(e^{-\sqrt{L}})\\
&\le\int_{\Omega_{i}}|u_{0}-\varphi_{c_{i}}|\left(u^{2}_{0}+u_{0}\varphi_{c_{i}}+\varphi^{2}_{c_{i}}\right)+O(e^{-\sqrt{L}})\\
&\le\|u_{0}-\varphi_{c_{i}}\|_{L^{2}(\Omega_{i})}\left(\int_{\Omega_{i}}\left
(u^{2}_{0}+u_{0}\varphi_{c_{i}}+\varphi^{2}_{c_{i}}\right)^{2}\right)^{1/2}+O(e^{-\sqrt{L}})\\
&\le\left(\|u_{0}-S_{Z^{0}}\|_{L^{2}(\Omega_{i})}+\sum\limits_{\underset{j \neq i}{1\le j\le N}}\|\varphi_{c_{j}}\|_{L^{2}(\Omega_{i})}\right)\cdot O(1)+O(e^{-\sqrt{L}})\\
&\le\left(\varepsilon^{2}+O(e^{-L/8})\right)\cdot O(1)+O(e^{-\sqrt{L}})\\
&\le  O(\varepsilon^{2})+O(e^{-\sqrt{L}}).
\end{align*}
This proves the lemma.
\hfill $ \square $ \vspace*{2mm}

\subsection{End of the proof of Theorem 1.1}
Let $u\in\mathcal{X}([0,T[)$, with $0<T\le +\infty$, be a solution of the DP equation satisfying \eqref{1.11}-\eqref{1.13} and \eqref{4.2} for some $ t_0\in ]0,T[ $. Let $M^{i}_{1}=v(t_{0},\xi^{i}_{1}(t_{0}))=\max_{x\in J_{i}}v(t_{0},x)$, with $J_{i}$'s as in \eqref{4.9}, and $\delta_{i}=c_{i}/6-M^{i}_{1}$. First, from \eqref{4.8} and \eqref{4.10}, we know that for $i=2,\ldots,N$, 
 $$\xi^{i}_{1}(t_{0})-\xi^{i-1}_{1}(t_{0})\ge \frac{2L}{3}>\frac{L}{2} .
 $$
  Applying \eqref{4.107} and \eqref{4.112} with $ u(t_0) $, we get 
\begin{equation}
\left\|u(t_{0})-\sum_{i=1}^{N}\varphi_{c_{i}}(\cdot-\xi^{i}_{1}(t_{0}))\right\|^{2}_{\mathcal{H}}\le 4\sum_{i=1}^{N}c_{i}\delta_{i}+O(\varepsilon^{2})+O(e^{-L/4}).
\label{4.116}
\end{equation}
In the same way, from  \eqref{4.100}  we get 
$$F_{i}(u(t_0))\le 18M^{i}_{1}E_{i}(u(t_0))-72(M^{i}_{1})^{3}+O(L^{-1/2}),$$
which leads to 
\begin{equation}\label{ww}
F(u(t_0))=\sum_{i=1}^{N}F_{i}(u(t_0))\le 18\sum_{i=1}^{N}M^{i}_{1}E_{i}(u(t_0))-72\sum_{i=1}^{N}(M^{i}_{1})^{3} + O(L^{-1/2})\, , 
\end{equation}
by summing over $i\in \{1,\ldots,N\} $.

 Now, we will use the following notation: for a function  $f :\mathbb{R}_{+}\mapsto\mathbb{R}$, we set 
\begin{equation}
\Delta^{t_0}_{0}f=f(t_0)-f(0).
\label{4.117}
\end{equation}
From \eqref{ww} and the fact that $ E(\cdot)$ and $ F(\cdot) $ are conservation laws for $ u $, we obtain
\begin{align}
0=\Delta^{t_0}_{0}F(u)&=\sum_{i=1}^{N}\Delta^{t_0}_{0}F_{i}(u)\le 18\sum_{i=1}^{N}M^{i}_{1}\Delta^{t_0}_{0}E_{i}(u)\nonumber\\
&+\sum_{i=1}^{N}\left[-72(M^{i}_{1})^{3}+18M^{i}_{1}E_{i}(u_{0})-F_{i}(u_{0})\right]+O(L^{-1/2}).
\label{4.118}
\end{align}
Note that, from \eqref{4.113} and \eqref{4.114}, for $0<\varepsilon<\varepsilon_{0}$ and $L>L_{0}>0$ with $\varepsilon_{0}\ll 1$ and $L_{0}\gg 1$, it holds
\begin{equation}
\sum_{i=1}^{N}\left[-72(M^{i}_{1})^{3}+18M^{i}_{1}E_{i}(u_{0})-F_{i}(u_{0})\right]=-72\sum_{i=1}^{N}\delta^{2}_{i}\left(M^{i}_{1}+\frac{c_{i}}{3}\right)+O(\varepsilon^{2})+O(e^{-\sqrt{L}}).
\label{4.119}
\end{equation}
Combining \eqref{4.118} and \eqref{4.119}, we get
$$\sum_{i=1}^{N}\delta^{2}_{i}\left(M^{i}_{1}+\frac{c_{i}}{3}\right)\le\frac{1}{4}\sum_{i=1}^{N}M^{i}_{1}\Delta^{t_0}_{0}E_{i}(u)+O(\varepsilon^{2})+O(L^{-1/2}),$$
and using the Abel transformation with $M^{0}_{1}=0$, we obtain  
\begin{equation}
\sum_{i=1}^{N}\delta^{2}_{i}\left(M^{i}_{1}+\frac{c_{i}}{3}\right)\le\frac{1}{4}\sum_{i=2}^{N}(M^{i}_{1}-M^{i-1}_{1})\Delta^{t}_{0}\mathcal{J}_{i,K}+O(\varepsilon^{2})+O(L^{-1/2})\, , 
\label{4.120}
\end{equation}
where $\mathcal{J}_{i,K}(t)$ is defined in \eqref{4.26}. 
From \eqref{4.2} we know that $u(t_{0})\in U(\gamma,L/2)$, on account of Lemma \ref{Lemma 4.1} there exists $\tilde{X}=(\tilde{x}_{1},\ldots,\tilde{x}_{N})$ with $\tilde{x}_{i}\in J_{i}$ such that $E\left(u(t_{0})-S_{\tilde{X}}\right)\le O(\gamma^{2})$, where $S_{\tilde{X}}$ is defined in \eqref{4.12}. Recalling that $v(t_{0},\xi^{i}_{1}(t_{0}))=\max_{x\in J_{i}}v(t_{0},x)$ and using \eqref{4.107}, we obtain $E\left(u(t_{0})-S_{\xi_{1}}\right)\le O(\gamma^{2})+O(e^{-L/4})$, with $\xi_{1}=(\xi^{1}_{1},\ldots,\xi^{N}_{1})$. From \eqref{4.6}, we deduce that 
$$\left\|v(t_{0})-\sum_{j=1}^{N}\rho_{c_{j}}(\cdot-\xi^{j}_{1}(t_{0}))\right\|_{L^{\infty}(\mathbb{R})}
\le O(\gamma)+O(e^{-L/8}).$$
Thus, we infer that 
$$v(x)=\sum_{j=1}^{N}\rho_{c_{j}}(\cdot-\xi^{j}_{1}(t_{0}))+O(\gamma)+O(e^{-L/8}),~~\forall x\in\mathbb{R}, $$
and applying this formula with $x=\xi^{i}_{1}(t_{0})$ and using that $\xi^{j}_{1}(t_{0})-\xi^{j-1}_{1}(t_{0})>L/2$, we get 
\begin{align*}
v(\xi^{i}_{1}(t_{0}))&=\sum_{j=1}^{N}\rho_{c_{j}}
(\xi^{i}_{1}(t_{0})-\xi^{j}_{1}(t_{0}))+O(\gamma)+O(e^{-L/8})\\
&=\frac{c_{i}}{6}+\sum\limits_{\underset{j \neq i}{1\le j\le N}}\rho_{c_{j}}
(\xi^{i}_{1}(t_{0})-\xi^{j}_{1}(t_{0}))+O(\gamma)+O(e^{-L/8})\\
&=\frac{c_{i}}{6}+O(\gamma)+O(e^{-L/8}).\\
\end{align*}
We take $\gamma=A(\sqrt{\varepsilon}+L^{-1/8})$, then $M^{i}_{1}=c_{i}/6+O(\sqrt{\varepsilon})+O(L^{-1/8})$. Therefore, for $0<\varepsilon<\varepsilon_{0}$
and $L>L_{0}>0$, with $\varepsilon_{0}\ll 1$ and $L_{0}\gg 1$, it holds
\begin{equation}
0<M^{1}_{1}<M^{2}_{1}<\ldots<M^{N}_{1}.
\label{4.121}
\end{equation}
Combining \eqref{4.120}, \eqref{4.121} and using the monotonicity estimate \eqref{4.28}, it holds
$$\sum_{i=1}^{N}\delta^{2}_{i}\left(M^{i}_{1}+\frac{c_{i}}{3}\right)\le O(\varepsilon^{2})
+O(L^{-1/8}).$$
Therefore, using that $(M^{i}_{1}+c_{i}/3)^{-1}<3/c_{i}$, there exists $C>0$ only depending on $(c_{i})^{N}_{i=1}$ such that
\begin{equation}
\delta_{i}\le C(\varepsilon+L^{-1/4}),~~i=1,\ldots,N.
\label{4.122}
\end{equation}
Now, combining \eqref{4.116} and \eqref{4.122}, we obtain 
$$\left\|u(t_{0})-\sum_{i=1}^{N}\varphi_{c_{i}}(\cdot-\xi^{i}_{1}(t_{0}))\right\|_{\mathcal{H}}\le C(\sqrt{\varepsilon}+L^{-1/8}),$$
and the  theorem follows by choosing $A=2C$.

\begin{Rem}[The Role of the Number of Extrema]\label{x1}
\normalfont
In the case where $v=(4-\partial^{2}_{x})^{-1}u$ admits a countable infinite number of local maximal values on some  $ [\alpha_i,\beta_i] $ (see \eqref{4.70}), with $ i\in \{1,\ldots,N\}$, it suffices to change the finite sums over $ j$ by  infinite sums in Lemmas \ref{Lemma 3.4}-\ref{Lemma 3.5} and Lemmas \ref{Lemma 4.3}-\ref{Lemma 4.4}.
\end{Rem}

\subsection*{Appendix. Proof of Lemma \ref{Lemma 4.2}}\label{Appendix}
The aim of this subsection is to prove Lemma \ref{Lemma 4.2}. Let us first  assume that $u$ is smooth solution. The case $u\in \mathcal{X}([0,T[)$ will follow by a density argument.

We compute the time variation of the following energy:
\begin{align*}
\frac{d}{dt}\int_{\mathbb{R}}yvg &=\int_{\mathbb{R}}y_{t}vg+\int_{\mathbb{R}}yv_{t}g\\
&=I+J.
\end{align*}
Applying the operator $(1-\partial^{2}_{x})(\cdot)$ on both sides of equation \eqref{1.6}, we get
$$y_{t}=-\frac{1}{2}(1-\partial^{2}_{x})\partial_{x}u^{2}-\frac{3}{2}\partial_{x}u^{2}$$
and substituting $y_{t}$ by this value, $I$ becomes
\begin{align*}
I&=-\frac{1}{2}\int_{\mathbb{R}}\left[(1-\partial^{2}_{x})\partial_{x}u^{2}\right]vg-\frac{3}{2}\int_{\mathbb{R}}\left[\partial_{x}(u^{2})\right]vg\\
&=I_{1}+I_{2}.
\end{align*}
By computing 
\begin{equation}
I_{2}=\frac{3}{2}\int_{\mathbb{R}}u^{2}\partial_{x}(vg)
=\frac{3}{2}\int_{\mathbb{R}}u^{2}v_{x}g+\frac{3}{2}\int_{\mathbb{R}}u^{2}vg'
\label{4.30}
\end{equation}
and
\begin{align*}
I_{1}&=\frac{1}{2}\int_{\mathbb{R}}\left[(1-\partial^{2}_{x})u^{2}\right]\partial_{x}(vg) \\
&=\frac{1}{2}\int_{\mathbb{R}}\left[(1-\partial^{2}_{x})u^{2}\right]v_{x}g+\frac{1}{2}\int_{\mathbb{R}}\left[(1-\partial^{2}_{x})u^{2}\right]vg'\\
&=I_{3}+I_{4}
\end{align*}
with
\begin{align}
I_{3}&=\frac{1}{2}\int_{\mathbb{R}}u^{2}v_{x}g-\frac{1}{2}\int_{\mathbb{R}}\partial^{2}_{x}(u^{2})v_{x}g\nonumber\\
&=\frac{1}{2}\int_{\mathbb{R}}u^{2}v_{x}g+\frac{1}{2}\int_{\mathbb{R}}\partial_{x}(u^{2})\partial_{x}(v_{x}g)\nonumber\\
&=\frac{1}{2}\int_{\mathbb{R}}u^{2}v_{x}g+\frac{1}{2}\int_{\mathbb{R}}\partial_{x}(u^{2})v_{xx}g+\frac{1}{2}\int_{\mathbb{R}}\partial_{x}(u^{2})v_{x}g'\nonumber\\
&=\frac{1}{2}\int_{\mathbb{R}}u^{2}v_{x}g-\frac{1}{2}\int_{\mathbb{R}}u^{2}\partial_{x}(v_{xx}g)-\frac{1}{2}\int_{\mathbb{R}}u^{2}\partial_{x}(v_{x}g')\nonumber\\
&=\frac{1}{2}\int_{\mathbb{R}}u^{2}v_{x}g-\frac{1}{2}\int_{\mathbb{R}}u^{2}v_{xxx}g-\int_{\mathbb{R}}u^{2}v_{xx}g'-\frac{1}{2}\int_{\mathbb{R}}u^{2}v_{x}g''
\label{4.31}
\end{align}
and
\begin{align}
I_{4}&=\frac{1}{2}\int_{\mathbb{R}}u^{2}vg'-\frac{1}{2}\int_{\mathbb{R}}\partial^{2}_{x}(u^{2})vg'\nonumber\\
&=\frac{1}{2}\int_{\mathbb{R}}u^{2}vg'+\frac{1}{2}\int_{\mathbb{R}}\partial_{x}(u^{2})\partial_{x}(vg')\nonumber\\
&=\frac{1}{2}\int_{\mathbb{R}}u^{2}vg'+\frac{1}{2}\int_{\mathbb{R}}\partial_{x}(u^{2})v_{x}g'+\frac{1}{2}\int_{\mathbb{R}}\partial_{x}(u^{2})vg''\nonumber\\
&=\frac{1}{2}\int_{\mathbb{R}}u^{2}vg'-\frac{1}{2}\int_{\mathbb{R}}u^{2}\partial_{x}(v_{x}g')-\frac{1}{2}\int_{\mathbb{R}}u^{2}\partial_{x}(vg'')\nonumber\\
&=\frac{1}{2}\int_{\mathbb{R}}u^{2}vg'-\frac{1}{2}\int_{\mathbb{R}}u^{2}v_{xx}g'-\int_{\mathbb{R}}u^{2}v_{x}g''-\frac{1}{2}\int_{\mathbb{R}}u^{2}vg'''.
\label{4.32}
\end{align}
Adding \eqref{4.31} and \eqref{4.32}, we get
\begin{align}
I_{1}&=\frac{1}{2}\int_{\mathbb{R}}u^{2}v_{x}g-\frac{1}{2}\int_{\mathbb{R}}u^{2}v_{xxx}g-\int_{\mathbb{R}}u^{2}v_{xx}g'-\frac{1}{2}\int_{\mathbb{R}}u^{2}v_{x}g''\nonumber\\
&\hspace{1cm}+\frac{1}{2}\int_{\mathbb{R}}u^{2}vg'-\frac{1}{2}\int_{\mathbb{R}}u^{2}v_{xx}g'-\int_{\mathbb{R}}u^{2}v_{x}g''-\frac{1}{2}\int_{\mathbb{R}} u^{2}vg'''\nonumber\\
&=\frac{1}{2}\int_{\mathbb{R}}u^{2}v_{x}g-\frac{1}{2}\int_{\mathbb{R}}u^{2}v_{xxx}g-\frac{3}{2}\int_{\mathbb{R}}u^{2}v_{xx}g'\nonumber\\
&\hspace{1cm}+\frac{1}{2}\int_{\mathbb{R}}u^{2}vg'-\frac{3}{2}\int_{\mathbb{R}}u^{2}v_{x}g''-\frac{1}{2}\int_{\mathbb{R}}u^{2}vg'''
\label{4.33}
\end{align}
and adding \eqref{4.30} and \eqref{4.33}, we get
\begin{align*}
I&=2\int_{\mathbb{R}}u^{2}v_{x}g-\frac{1}{2}\int_{\mathbb{R}}u^{2}v_{xxx}g-\frac{3}{2}\int_{\mathbb{R}}u^{2}v_{xx}g'\\
&\hspace{1cm}+2\int_{\mathbb{R}}u^{2}vg'-\frac{3}{2}\int_{\mathbb{R}}u^{2}v_{x}g''-\frac{1}{2}\int_{\mathbb{R}}u^{2}vg'''\\
&=\frac{1}{2}\int_{\mathbb{R}}u^{2}\left[(4-\partial^{2}_{x})v_{x}\right]g+\frac{3}{2}\int_{\mathbb{R}}u^{2}\left[\left(\frac{4}{3}-\partial^{2}_{x}\right)v\right]g'\\
&\hspace{1cm}-\frac{3}{2}\int_{\mathbb{R}}u^{2}v_{x}g''-\frac{1}{2}\int_{\mathbb{R}}u^{2}vg'''.
\end{align*}
The first two integrals give us
$$\frac{1}{2}\int_{\mathbb{R}}u^{2}\left[(4-\partial^{2}_{x})v_{x}\right]g =\frac{1}{2}\int_{\mathbb{R}}u^{2}u_{x}g
=\frac{1}{6}\int_{\mathbb{R}}\partial_{x}(u^{3})g\\
=-\frac{1}{6}\int_{\mathbb{R}}u^{3}g'$$
and
$$\frac{3}{2}\int_{\mathbb{R}}u^{2}\left[\left(\frac{4}{3}-\partial^{2}_{x}\right)v\right]g'=\frac{3}{2}\int_{\mathbb{R}}u^{2}\left[(4-\partial^{2}_{x})v-\frac{8}{3}v\right]g'
=\frac{3}{2}\int_{\mathbb{R}}u^{3}g'-4\int_{\mathbb{R}}u^{2}vg'.$$
Finally, we obtain 
\begin{equation}
I=\frac{4}{3}\int_{\mathbb{R}}u^{3}g'-4\int_{\mathbb{R}}u^{2}vg'-\frac{3}{2}\int_{\mathbb{R}}u^{2}v_{x}g''-\frac{1}{2}\int_{\mathbb{R}}u^{2}vg'''.
\label{4.34}
\end{equation}
We set $h=(1-\partial^{2}_{x})^{-1}u^{2}$. Applying the operator $(4-\partial^{2}_{x})^{-1}(\cdot)$ on both sides of equation \eqref{1.6} and using \eqref{4.17}, we get
\begin{equation}
v_{t}=-\frac{1}{2}(1-\partial^{2}_{x})^{-1}\partial_{x}u^{2}
=-\frac{1}{2}h_{x}.
\label{4.35}
\end{equation}
Substituting $v_{t}$ by this value, $J$ becomes
\begin{align*}
J&=-\frac{1}{2}\int_{\mathbb{R}}yh_{x}g\\
&=\frac{1}{2}\int_{\mathbb{R}}y_{x}hg+\frac{1}{2}\int_{\mathbb{R}}yhg'\\
&=J_{1}+J_{2}.
\end{align*}
By computing
\begin{align}
J_{2}&=\frac{1}{2}\int_{\mathbb{R}}(u-u_{xx})hg'\nonumber\\
&=\frac{1}{2}\int_{\mathbb{R}}uhg'-\frac{1}{2}\int_{\mathbb{R}}u_{xx}hg'\nonumber\\
&=\frac{1}{2}\int_{\mathbb{R}}uhg'+\frac{1}{2}\int_{\mathbb{R}}u_{x}\partial_{x}(hg')\nonumber\\
&=\frac{1}{2}\int_{\mathbb{R}}uhg'+\frac{1}{2}\int_{\mathbb{R}}u_{x}h_{x}g'+\frac{1}{2}\int_{\mathbb{R}}u_{x}hg''\nonumber\\
&=\frac{1}{2}\int_{\mathbb{R}}uhg'-\frac{1}{2}\int_{\mathbb{R}}u\partial_{x}(h_{x}g')-\frac{1}{2}\int_{\mathbb{R}}u\partial_{x}(hg'')\nonumber\\
&=\frac{1}{2}\int_{\mathbb{R}}uhg'-\frac{1}{2}\int_{\mathbb{R}}uh_{xx}g'-\int_{\mathbb{R}}uh_{x}g''-\frac{1}{2}\int_{\mathbb{R}}uhg'''
\label{4.36}
\end{align}
and 
\begin{align*}
J_{1}&=\frac{1}{2}\int_{\mathbb{R}}(u_{x}-u_{xxx})hg\\
&=\frac{1}{2}\int_{\mathbb{R}}u_{x}hg-\frac{1}{2}\int_{\mathbb{R}}u_{xxx}hg\\
&=\frac{1}{2}\int_{\mathbb{R}}u_{x}hg+\frac{1}{2}\int_{\mathbb{R}}u_{xx}\partial_{x}(hg)\\
&=\frac{1}{2}\int_{\mathbb{R}}u_{x}hg+\frac{1}{2}\int_{\mathbb{R}}u_{xx}h_{x}g+\frac{1}{2}\int_{\mathbb{R}}u_{xx}hg'\\
&=J_{3}+J_{4}+J_{5}
\end{align*}
with
\begin{equation}
J_{3}=-\frac{1}{2}\int_{\mathbb{R}}u\partial_{x}(hg)
=-\frac{1}{2}\int_{\mathbb{R}}uh_{x}g-\frac{1}{2}\int_{\mathbb{R}}uhg'\, ,
\label{4.37}
\end{equation}
\begin{align}
J_{5}&=-\frac{1}{2}\int_{\mathbb{R}}u_{x}\partial_{x}(hg')\nonumber\\
&=-\frac{1}{2}\int_{\mathbb{R}}u_{x}h_{x}g'-\frac{1}{2}\int_{\mathbb{R}}u_{x}hg''\nonumber\\
&=\frac{1}{2}\int_{\mathbb{R}}u\partial_{x}(h_{x}g')+\frac{1}{2}\int_{\mathbb{R}}u\partial_{x}(hg'')\nonumber\\
&=\frac{1}{2}\int_{\mathbb{R}}uh_{xx}g'+\int_{\mathbb{R}}uh_{x}g''+\frac{1}{2}\int_{\mathbb{R}}uhg'''
\label{4.38}
\end{align}
and 
\begin{align}
J_{4}&=-\frac{1}{2}\int_{\mathbb{R}}u_{x}\partial_{x}(h_{x}g)\nonumber\\
&=-\frac{1}{2}\int_{\mathbb{R}}u_{x}h_{xx}g-\frac{1}{2}\int_{\mathbb{R}}u_{x}h_{x}g'\nonumber\\
&=\frac{1}{2}\int_{\mathbb{R}}u\partial_{x}(h_{xx}g)+\frac{1}{2}\int_{\mathbb{R}}u\partial_{x}(h_{x}g')\nonumber\\
&=\frac{1}{2}\int_{\mathbb{R}}uh_{xxx}g+\int_{\mathbb{R}}uh_{xx}g'+\frac{1}{2}\int_{\mathbb{R}}uh_{x}g''.
\label{4.39}
\end{align}
Adding \eqref{4.36}-\eqref{4.39}, we get  
$$J=-\frac{1}{2}\int_{\mathbb{R}}uh_{x}g+\frac{1}{2}\int_{\mathbb{R}}uh_{xxx}g+\int_{\mathbb{R}}uh_{xx}g'+\frac{1}{2}\int_{\mathbb{R}}uh_{x}g''.$$
Using that $h_{xx}=-u^{2}+h$ and $h_{xxx}=-2uu_{x}+h_{x}$, we have
$$\int_{\mathbb{R}}uh_{xx}g'=\int_{\mathbb{R}}u(-u^{2}+h)g'
=-\int_{\mathbb{R}}u^{3}g'+\int_{\mathbb{R}}uhg'$$
and
\begin{align*}
\frac{1}{2}\int_{\mathbb{R}}uh_{xxx}g&=\frac{1}{2}\int_{\mathbb{R}}u(-2uu_{x}+h_{x})g\\
&=-\int_{\mathbb{R}}u^{2}u_{x}g+\frac{1}{2}\int_{\mathbb{R}}uh_{x}g\\
&=-\frac{1}{3}\int_{\mathbb{R}}\partial_{x}(u^{3})g+\frac{1}{2}\int_{\mathbb{R}}uh_{x}g\\
&=\frac{1}{3}\int_{\mathbb{R}}u^{3}g'+\frac{1}{2}\int_{\mathbb{R}}uh_{x}g.
\end{align*}
At this stage it is worth noticing that the term $\int_{\mathbb{R}}uh_{x}g$ cancels with the one in $J$.
Finally, we obtain
\begin{equation}
J=-\frac{2}{3}\int_{\mathbb{R}}u^{3}g'+\int_{\mathbb{R}}uhg'+\frac{1}{2}\int_{\mathbb{R}}uh_{x}g''.
\label{4.40}
\end{equation}
Combining \eqref{4.34} and \eqref{4.40}, we get 
\begin{equation}
\frac{d}{dt}\int_{\mathbb{R}}yvg=\frac{2}{3}\int_{\mathbb{R}}u^{3}g'-4\int_{\mathbb{R}}u^{2}vg'-\frac{3}{2}\int_{\mathbb{R}}u^{2}v_{x}g''-\frac{1}{2}\int_{\mathbb{R}}u^{2}vg'''+\int_{\mathbb{R}}uhg'+\frac{1}{2}\int_{\mathbb{R}}uh_{x}g''.
\label{4.41}
\end{equation}
Now, substituting $u$ by $4v-v_{xx}$ and using integration by parts, we rewrite the energy as
\begin{align*}
\int_{\mathbb{R}}yvg&=\int_{\mathbb{R}}v\left[(1-\partial^{2}_{x})(4v-v_{xx})\right]g\\
&=4\int_{\mathbb{R}}v^{2}g-5\int_{\mathbb{R}}vv_{xx}g
+\int_{\mathbb{R}}v(\partial^{4}_{x}v)g\\
&=4\int_{\mathbb{R}}v^{2}g+K_{1}+K_{2}.\\
\end{align*}
By computing
\begin{align}
K_{1}&=5\int_{\mathbb{R}}\partial_{x}(vg)v_{x}\nonumber\\
&=5\int_{\mathbb{R}}v^{2}_{x}g+5\int_{\mathbb{R}}vv_{x}g'\nonumber\\
&=5\int_{\mathbb{R}}v^{2}_{x}g+\frac{5}{2}\int_{\mathbb{R}}\partial_{x}(v^{2})g'\nonumber\\
&=5\int_{\mathbb{R}}v^{2}_{x}g-\frac{5}{2}\int_{\mathbb{R}}v^{2}g''
\label{4.42}
\end{align}
and
\begin{align}
K_{2}&=-\int_{\mathbb{R}}\partial_{x}(vg)v_{xxx}\nonumber\\
&=-\int_{\mathbb{R}}v_{x}v_{xxx}g-\int_{\mathbb{R}}vv_{xxx}g'\nonumber\\
&=K_{3}+K_{4}
\label{4.43}
\end{align}
with
\begin{align}
K_{3}&=\int_{\mathbb{R}}\partial_{x}(v_{x}g)v_{xx}\nonumber\\
&=\int_{\mathbb{R}}v^{2}_{xx}g+\int_{\mathbb{R}}v_{x}v_{xx}g'\nonumber\\
&=\int_{\mathbb{R}}v^{2}_{xx}g+\frac{1}{2}\int_{\mathbb{R}}\partial_{x}(v^{2}_{x})g'\nonumber\\
&=\int_{\mathbb{R}}v^{2}_{xx}g-\frac{1}{2}\int_{\mathbb{R}}v^{2}_{x}g''
\label{4.44}
\end{align}
and
\begin{align}
K_{4}&=\int_{\mathbb{R}}\partial_{x}(vg')v_{xx}\nonumber\\
&=\int_{\mathbb{R}}v_{x}v_{xx}g'+\int_{\mathbb{R}}vv_{xx}g''\nonumber\\
&=\frac{1}{2}\int_{\mathbb{R}}\partial_{x}(v^{2}_{x})g'-\int_{\mathbb{R}}\partial_{x}(vg'')v_{x}\nonumber\\
&=-\frac{1}{2}\int_{\mathbb{R}}v^{2}_{x}g''-\int_{\mathbb{R}}v^{2}_{x}g''-\int_{\mathbb{R}}vv_{x}g'''\nonumber\\
&=-\frac{3}{2}\int_{\mathbb{R}}v^{2}_{x}g''-\frac{1}{2}\int_{\mathbb{R}}\partial_{x}(v^{2})g'''\nonumber\\
&=-\frac{3}{2}\int_{\mathbb{R}}v^{2}_{x}g''+\frac{1}{2}\int_{\mathbb{R}}v^{2}g^{(4)}.
\label{4.45}
\end{align}
Combining \eqref{4.42}-\eqref{4.45}, we get 
$$\int_{\mathbb{R}}yvg=\int_{\mathbb{R}}\left(4v^{2}+5v^{2}_{x}+v^{2}_{xx}\right)g
+\frac{1}{2}\int_{\mathbb{R}}v^{2}(g^{(4)}-5g'')-2\int_{\mathbb{R}}v^{2}_{x}g''$$
and differentiating with respect to time
\begin{equation}
\frac{d}{dt}\int_{\mathbb{R}}yvg=\frac{d}{dt}\int_{\mathbb{R}}\left(4v^{2}+5v^{2}_{x}+v^{2}_{xx}\right)g+L_{1}+L_{2}.
\label{4.46}
\end{equation}
Using \eqref{4.35}, we have 
\begin{equation}
L_{1}=\int_{\mathbb{R}}vv_{t}(g^{(4)}-5g'')
=\frac{5}{2}\int_{\mathbb{R}}vh_{x}g''-\frac{1}{2}\int_{\mathbb{R}}vh_{x}g^{(4)}
\label{4.47}
\end{equation}
and 
\begin{equation}
L_{2}=-4\int_{\mathbb{R}}v_{x}v_{tx}g''
=2\int_{\mathbb{R}}v_{x}h_{xx}g''
=-2\int_{\mathbb{R}}u^{2}v_{x}g''+2\int_{\mathbb{R}}v_{x}hg''.
\label{4.48}
\end{equation}
Lemma \ref{Lemma 4.2} follows by combining \eqref{4.41} and \eqref{4.46}-\eqref{4.48}.
\hfill $ \square $ \vspace*{2mm}

\begin{Ack}
\normalfont
The author would like to thank his PhD advisor Luc Molinet for his help and his careful reading of this manuscript.
\end{Ack}

\bibliographystyle{plain}
\bibliography{mabiblio}
\end{document}